\UseRawInputEncoding

\documentclass[12pt,a4paper]{amsart}
\usepackage{lmodern}
\usepackage{amsmath}
\usepackage[utf8]{inputenc}
\usepackage{bbm}
\usepackage{tikz}
\usepackage{tikz-cd}
\usepackage{utfsym}
\usepackage{mathrsfs}
\usepackage{fix-cm}
\usepackage{cite}
\usepackage{mathtools}
\usetikzlibrary{matrix}
\usepackage{tikz-cd}
\usepackage{fontawesome5}

\usepackage[a4paper,left=1.8cm,right=1.8cm,top=2.5cm,bottom=2.5cm]{geometry}

\usepackage{dsfont}         
                          
\usepackage{float}
\usepackage{xurl}
\usepackage{ragged2e}
\usepackage{xcolor}
 
\usepackage{enumitem}
\usepackage{float}

\setlist{nosep}
\setlist[enumerate]{label=(\roman*)}
 
\usepackage{etoolbox}
\AtBeginEnvironment{thebibliography}{%
  \setlength{\itemsep}{0pt}%
  \setlength{\parskip}{0pt}%
  \setlength{\topsep}{0pt}%
  \setlength{\partopsep}{0pt}%
}
 
\newlist{defenum}{enumerate}{3}
\setlist[defenum]{label=(\alph*),leftmargin=*,align=left}
\newlist{defsubenum}{enumerate}{1}
\setlist[defsubenum]{label=(\roman*),leftmargin=*,align=left}

\usepackage{mathtools}
\usepackage{amssymb,amsthm}
\usepackage{graphicx}
\usepackage[cal=rsfso,scaled=1.05]{mathalfa}

\usepackage[english]{babel}
\usepackage[T1]{fontenc}
\usepackage{tikz-cd}
\usepackage{bm}

\usepackage[most]{tcolorbox}

\usepackage{array,booktabs,tabularx,multirow}

\usepackage[colorlinks=true,pagebackref,hypertexnames=false,linkcolor=red,citecolor=blue,filecolor=blue,urlcolor=blue]{hyperref}

\usepackage[nameinlink,noabbrev]{cleveref}
\usepackage{autonum}
\usepackage[alphabetic,backrefs]{amsrefs}

\makeatletter
\renewcommand\normalsize{%
    \@setfontsize\normalsize{11.7}{14pt plus .3pt minus .3pt}%
    \abovedisplayskip 10\p@ \@plus4\p@ \@minus4\p@
    \abovedisplayshortskip 6\p@ \@plus2\p@
    \belowdisplayshortskip 6\p@ \@plus2\p@
    \belowdisplayskip \abovedisplayskip}
\renewcommand\small{%
    \@setfontsize\small{9.5}{12\p@ plus .2\p@ minus .2\p@}%
    \abovedisplayskip 8.5\p@ \@plus4\p@ \@minus1\p@
    \belowdisplayskip \abovedisplayskip
    \abovedisplayshortskip \abovedisplayskip
    \belowdisplayshortskip \abovedisplayskip}
\renewcommand\footnotesize{%
    \@setfontsize\footnotesize{8.5}{9.25\p@ plus .1pt minus .1pt}
    \abovedisplayskip 6\p@ \@plus4\p@ \@minus1\p@
    \belowdisplayskip \abovedisplayskip
    \abovedisplayshortskip \abovedisplayskip
    \belowdisplayshortskip \abovedisplayskip}

\ifdefined\pdfpagewidth
\else
\fi
\calclayout
\makeatother

\newtheorem{theorem}{Theorem}[section]
\newtheorem{Mtheorem}{Theorem}

\newtheorem{corollary}[theorem]{Corollary}

\newtheorem{lemma}[theorem]{Lemma}
\newtheorem{conjecture}[theorem]{Conjecture}
\newtheorem{proposition}[theorem]{Proposition}
\theoremstyle{definition}
\newtheorem{definition}[theorem]{Definition}
\newtheorem{example}[theorem]{Example}   
\theoremstyle{remark}
\newtheorem{remark}[theorem]{Remark}
\newcommand{\ms}[1]{\{\!\{#1\}\!\}}

\newcommand{\cl}[1]{\overline{#1}}

\crefname{theorem}{Theorem}{Theorems}
\crefname{Mtheorem}{Main Theorem}{Main Theorems}
\crefname{lemma}{Lemma}{Lemmata}
\crefname{corollary}{Corollary}{Corollaries}
\crefname{proposition}{Proposition}{Propositions}
\crefname{definition}{Definition}{Definitions}
\crefname{example}{Example}{Examples}    

\AtEndEnvironment{example}{\exend}

\newcommand{\fm}{\mathfrak{m}}

\newcommand{\Q}{\mathbb{Q}}
\newcommand{\R}{\mathbb{R}}
\newcommand{\Z}{\mathbb{Z}}
\newcommand{\C}{\mathbb{C}}

\newcommand{\PP}{\mathbb{P}}
\newcommand{\T}{\mathbb{T}}
{\End}{\mathrm{End}}
\newcommand{\Hom}{\mathrm{Hom}}

\newcommand{\rank}{\mathrm{rk}}

\newcommand{\rH}{\mathrm{H}}

\newcommand{\cE}{\mathcal{E}}

\newcommand{\cO}{\mathcal{O}}

\newcommand{\Span}{\mathrm{Span}}

\DeclareMathOperator{\Spec}{\mathrm{Spec}}
\DeclareMathOperator{\Cone}{Cone}

\newcolumntype{L}{>{$}l<{$}}
\newcolumntype{C}{>{$}c<{$}}
\newcolumntype{R}{>{$}r<{$}}

\newcommand{\Pic}{\text{Pic}}

\lstdefinestyle{gwflags}{
  basicstyle=\ttfamily\small,
  columns=fullflexible,
  keepspaces=true,
  showstringspaces=false,
  morecomment=[l]{\#},
  commentstyle=\itshape\color{black!55},
  aboveskip=0pt, belowskip=0pt,
  frame=none
}

\NewTCBListing{inputbox}{ o }{
  listing only,
  listing style=gwflags,
  colback=black!3, colframe=black!80,
  boxrule=0.8pt, arc=0pt, outer arc=0pt,
  left=8pt, right=8pt, top=6pt, bottom=6pt,
  fonttitle=\bfseries, coltitle=black,
  attach boxed title to top left={yshift=-2mm, xshift=4mm},
  boxed title style={colback=white, colframe=white},
  IfValueT={#1}{title={#1}}
}

\NewTColorBox{outputbox}{ o }{
  breakable,
  colback=white, colframe=black!80,
  boxrule=0.8pt, arc=0pt, outer arc=0pt,
  left=8pt, right=8pt, top=8pt, bottom=8pt,
  fonttitle=\bfseries, coltitle=black,
  attach boxed title to top left={yshift=-2mm, xshift=4mm},
  boxed title style={colback=white, colframe=white},
  IfValueT={#1}{title={#1}}
}
 
\newcommand{\exend}{\hfill$\diamond$}
 
\begin{document}
\title[On some atomic decompositions in toric varieties]{On the atomic decomposition of complete intersections in toric varieties}

 \author[Brockhof]{Milena Brockhof}
 \address{Instituto de Matemática, Estatística e Computação Científica (IMECC) da Universidade Estadual de Campinas (Unicamp), Cidade Universitária, Campinas - SP, 13083-856, Brazil.} 
 \email{mbrockhof2001@gmail.com}

 \author[Galindo]{Giovane Galindo}
 \address{Institute of Mathematics and Informatics, Bulgarian Academy of Sciences}
 \email{giovanepgneto@gmail.com}

 \author[Katzarkov]{Ludmil Katzarkov}
 \address{Simons Center for Geometry and Physics, Stony Brook University, Stony Brook, NY 11794, USA. \& International Center for Mathematical Sciences (ICMS), Sofia, Bulgaria. \& Institute of Mathematics and Informatics, Bulgarian Academy of Sciences, Sofia, Bulgaria.} 
 \email{lkatzarkov@gmail.com}

 \author[Martins]{Pedro Antonio Muniz Martins}
 \address{Instituto de Matemática, Estatística e Computação Científica (IMECC) da Universidade Estadual de Campinas (Unicamp), Cidade Universitária, Campinas - SP, 13083-856, Brazil.} 
 \email{pedroa.muniz9@gmail.com}

\begin{abstract}
Building on Liu's and Liu--Sheshmani's systematized description of Kontsevich's localization methods in toric varieties, we develop a computational framework for genus-zero Gromov--Witten invariants of smooth projective toric varieties and their smooth complete intersections. This is grounded in equivariant localization and the combinatorial structure of toric fans. We describe how the geometric data required by the localization formulas can be recovered directly from the fan and extend the existing software for flag varieties to compute the small quantum multiplication by the first Chern class. The implementation also incorporates cyclic determinants and regularized-exponent multisets, providing computational tools for investigating the atomic decomposition of quantum cohomology.

Using this framework, we obtain the nonrationality of every smooth complete intersection of a quadric and a cubic in $\PP^5$. Also applying new methods to study Hodge atoms in Fano fivefolds, we recover the nonrationality of every smooth complete intersection of a quadric and a cubic in $\PP^7$, due to Cavenaghi--Katzarkov--Kontsevich, and we prove the nonrationality of every smooth quartic double fivefold and every smooth quartic fivefold.
\end{abstract}

\maketitle
\tableofcontents

\section{Introduction}
\label{sec:intro}

The theory of Hodge atoms, introduced by Katzarkov--Kontsevich--Pantev--Yu~\cite{KKPY}, offers a new approach to nonrationality problems by combining Gromov--Witten theory with Hodge theory. The first major application is the proof of the non-rationality of the very general cubic fourfold~\cite{KKPY}*{Theorem 6.8}.

Recent developments of the theory include~\cites{alexeev2026atomic, Benedetti_Fay_Guere_Manivel_Perrin, JeremyCubic, Leo_Luiz_Pedro_Lud_Max, CKK-lectures, Andrew, AtomsMeetSymbols, Benedetti_Vladimiro_Manivel, cavenaghi2026atomicapproachshokurovsconicbundle}. We recall some of these results, without attempting an exhaustive overview:

\begin{itemize}
    \item[1)] In \cite{AtomsMeetSymbols}, the authors initiate the theory of $G$-equivariant Hodge atoms. For a smooth hypersurface $X(1,1,1,1)\subset (\PP^1)^4$ of multidegree $(1,1,1,1)$ invariant under the involution interchanging the first two and the last two factors, they prove that this $\Z/2$-action is not equivariantly birational to a linearizable $\Z/2$-action on $\PP^3$~\cite{AtomsMeetSymbols}*{Theorem 2.16}. The same paper initiates the theory of Chen--Ruan atoms; a full orbifold birational theory would require an appropriate analog of Iritani's blowup decomposition for ordinary quantum cohomology~\cite{iritani2025quantumcohomologyblowups}*{Theorem 5.18}. Building on the equivariant theory, Fay proves that a very general symmetric Verra fourfold is not $\Z/2$-birational to $\PP^4$~\cite{Fay}.

    \item[2)] In \cite{Benedetti_Vladimiro_Manivel}, the authors use Hodge atoms to prove the nonrationality of the very general Gushel--Mukai fourfold~\cite{Benedetti_Vladimiro_Manivel}*{Theorem 4}. They compute the matrix of the small quantum multiplication by the first Chern class restricted to the ambient cohomology and use the atomic irrationality criterion of~\cite{KKPY}*{Proposition 5.30}.

    \item[3)] In \cite{Benedetti_Fay_Guere_Manivel_Perrin}, the authors use Hodge atoms to establish a numerical irrationality criterion for fourfolds in terms of their Hodge numbers~\cite{Benedetti_Fay_Guere_Manivel_Perrin}*{Theorem 4.1}. They apply it to the Hodge-general Küchle fourfold of type $(c5)$~\cite{Benedetti_Fay_Guere_Manivel_Perrin}*{Proposition 5.1}.

    \item[4)] In \cite{alexeev2026atomic}, the authors use localization methods originating with Kontsevich~\cite{Kontsevich1995} to develop software computing small quantum multiplication by the first Chern class for complete intersections in flag varieties. They make computed matrices available through \href{https://leonardofcavenaghi.github.io/atomicdecompositions/}{GEOATOMS} and prove new irrationality results, including the nonrationality of an intersection of five hyperplanes and a quadric in the spinor tenfold~\cite{alexeev2026atomic}*{Theorem 6.18}.

    \item[5)] In \cite{Leo_Luiz_Pedro_Lud_Max}, the authors obtain a numerical obstruction to rationality from cyclic factors of the formal quantum connection and their regularized exponents~\cite{Leo_Luiz_Pedro_Lud_Max}*{Theorem 4.5}. Related formal-connection methods are used in~\cite{Andrew} to prove nonrationality for families of Fano complete-intersection fourfolds in weighted projective spaces~\cite{Andrew}*{Theorem 1.1}.
\end{itemize}

Roughly speaking, the computation of Hodge atoms involves studying the spectral decomposition of the operator
\[
\boldsymbol{\kappa}_b
=\mathsf{Eu}\star_b(-)
\]
where $\mathsf{Eu}$ is the Euler vector field, $\star_b$ denotes big quantum multiplication at $b \in B_X^{\mathsf{Hod}}$, that is, $b$ is in the fixed locus of the Hodge group acting on the base of the A-model F-bundle. More precisely, the generic atomic decomposition is controlled by this operator on the locus in $B_X^{\mathsf{Hod}}$ where the number of distinct eigenvalues is maximal. Since computing big quantum multiplication is generally difficult, we instead study the operator on the small quantum locus.

To describe the small quantum product, let $\{T_i\}$ be a homogeneous basis of $\rH^*(X)$, and let $(g^{ab})$ denote the inverse of the Poincaré pairing matrix $g_{ab}=\int_X T_a\cup T_b$. Then
\[
T_i\star_{\mathrm{small}}T_j
=
\sum_{a,b}
\left(
\sum_{\beta\in\mathsf{NE}(X)}
\big\langle T_i,T_j,T_a\big\rangle_{3,\beta}
q^\beta
\right)g^{ab}T_b,
\]
where $\mathsf{NE}(X)$ denotes the semigroup of effective integral curve classes, $\langle T_i,T_j,T_a\rangle_{3,\beta}$ denotes the genus-zero three-point Gromov--Witten invariant of class $\beta$, and $q^\beta$ is the corresponding Novikov monomial.

When $X$ is a smooth projective Fano variety, the ampleness of $c_1(TX)$ and the dimension constraint imply that only finitely many effective curve classes contribute to any fixed small quantum product. For a complete intersection computed through an ambient variety, however, turning this finiteness statement into an explicit stopping rule may require further bounds on the ambient images of effective curve classes, as discussed in Appendix~\ref{app:toricsoftware}.

Of particular interest is the operator
\[
c_1(TX)\star_{\mathrm{small}}(-),
\]
which corresponds to the Euler multiplication operator at the ``small locus'' in $B_X^{\mathsf{Hod}}$. Although its spectral decomposition need not coincide with the generic atomic decomposition, it can provide useful constraints on the possible Hodge atoms. In certain cases, this information is already sufficient to obtain obstructions to rationality, without explicitly computing the full big quantum product.

\bigskip

The localization methods proposed in \cite{Kontsevich1995} and further developed in \cite{Liu_Chiu}, particularly in \cite{liu2013localization}, which focuses on the toric case, provide a way to compute the Gromov–Witten invariants of \emph{GKM spaces}. Here we use the term for smooth projective varieties equipped with an algebraic torus action having finitely many fixed points and finitely many one-dimensional torus orbits. Our focus is on toric varieties, which are especially well-suited to these methods because the combinatorial structure of their defining fan determines their geometry entirely.

In \cite{alexeev2026atomic}, software was developed to carry out the computations of genus zero Gromov-Witten invariants in flag varieties and their complete intersections. The first main goal of the present paper is to provide the theoretical background for extending this software to toric varieties and their complete intersections. This concerns using the same machinery developed for decorated trees, but switching the combinatorial structure previously given from Lie-theoretic data to fan combinatorics. We also add new features discussed in \cite{Leo_Luiz_Pedro_Lud_Max}: the cyclic determinant and the regularized exponent multiset. These two inclusions are telling: the first allows us to understand the sufficient conditions for small quantum multiplication to yield the actual atomic splitting, while the second concerns a new type of obstruction to rationality arising from the quantum differential equation.

\ 

The second main goal of this paper, which justifies the software updates we propose, concerns birational-geometry applications arising from computations performed by the new software version. Here we present only some of the computations (and applications) we have obtained; a full list of computations, as well as the source code, may be found on the website \href{https://leonardofcavenaghi.github.io/atomicdecompositions/}{GEOATOMS}. We prove:

\begin{Mtheorem}[= Theorem~\ref{thm:quadric_cubic_P5}]\label{Mtheorem:quadric_cubic_P5}
    Every smooth intersection of a quadric and cubic in $\mathbb{P}^{5}$ is nonrational.
\end{Mtheorem}
Theorem~\ref{Mtheorem:quadric_cubic_P5} strengthens the classical result of \cite{IskovskikhPukhlikov1996}, which establishes nonrationality for a general smooth complete intersection of this type.

\ 

Although a complete classification of atoms appearing in the atomic decomposition of smooth projective threefolds is not yet available, it is still possible to prove the nonrationality of some smooth Fano fivefolds~\cite{IlievManivel2015}. The five-dimensional criterion used below is due to Cavenaghi--Katzarkov--Kontsevich~\cite{CKK-lectures}. We recall their nonrationality result for the $(2,3)$ complete intersection in $\PP^7$ and provide only a brief account of its application. We claim no originality for this result:
\begin{theorem}[Cavenaghi--Katzarkov--Kontsevich]\label{Mtheorem:CKK}
     Every smooth intersection of a quadric and a cubic in $\mathbb{P}^{7}$ is nonrational.
\end{theorem}

This example is one of the three families studied in \cite{IlievManivel2015}, where they classify all complete intersections of Calabi--Yau type in weighted projective spaces of dimension at least four:

\begin{proposition}[\cite{IlievManivel2015}, Section 3.1]
Let $X$ be a smooth complete intersection of Calabi--Yau
type in a weighted projective space. Assume that
$\dim X \geq 4$. Then $X$ is necessarily one of the
following:
\begin{enumerate}
    \item A smooth cubic sevenfold in $\mathbb{P}^8$;
    \item A smooth quartic double fivefold in
    $\mathbb{P}(1,1,1,1,1,1,2)$;
    \item A transverse intersection of a smooth cubic
    and a smooth quadric in $\mathbb{P}^7$.
\end{enumerate}
\end{proposition}

While the cubic sevenfold remains outside the scope of the fivefold criterion used here, the quartic double fivefold can be realized as a hypersurface in a smooth toric variety. The criterion of Cavenaghi--Katzarkov--Kontsevich then yields:
\begin{Mtheorem}[= Theorem~\ref{thm:non_rat_quartic_double_fivefold}]\label{thm:intro_double_fivefold}
    Every smooth quartic double fivefold is nonrational.
\end{Mtheorem}

Theorem~\ref{thm:intro_double_fivefold} verifies a conjecture formulated by Cheng--Perry--Zhao in~\cite{Cheng2025}, inspired by Kuznetsov's categorical approach to rationality. To properly contextualize it, recall that for a smooth projective Fano variety $X$ with $\mathrm{Pic}(X)\cong \Z[H]$, where $H$ is the ample generator and $-K_X=rH$, the semiorthogonal decomposition of the bounded derived category of coherent sheaves is
\[
D^{b}(\mathrm{Coh}(X)) = \langle \mathsf{Ku}(X) , \mathcal{O}_X ,\mathcal{O}_X(H) , \ldots, \mathcal{O}_X\left( (r-1)H  \right) \rangle.
\]
The component $\mathsf{Ku}(X)$ is called the \emph{Kuznetsov component}, and is expected to be related to birational geometry. One prominent example is Kuznetsov's original rationality prediction for cubic fourfolds, whose sufficiency direction has since been disproved.
\begin{conjecture}[Kuznetsov, original formulation]
   A smooth cubic fourfold $X$ is rational if and only if there exists a smooth projective K3 surface $K$ and an exact equivalence
\[
\mathsf{Ku}(X)\simeq D^b(\mathrm{Coh}(K)).
\]
\end{conjecture}
The recent work~\cite{OpenAI2026IrrationalCubicFourfolds} provides counterexamples to the ``if'' direction: it constructs irrational smooth cubic fourfolds whose Kuznetsov components are equivalent to $D^{b}(\mathrm{Coh}(K))$ for smooth projective K3 surfaces $K$. Gu\'er\'e~\cite{guere2026integralhodgestructuresatoms} interprets the resulting birational obstruction in terms of integral Hodge structures and atoms, using the compatibility of Iritani's blowup formula with the classical blowup decomposition. In particular, his result gives a necessary integral Hodge-theoretic condition for a cubic fourfold to be rational.

An analogous prediction to Kuznetsov's, now for quartic double fivefolds, was formulated by Cheng--Perry--Zhao.

\begin{conjecture}[\cite{Cheng2025}*{Conjecture 1.2}]
    If $X$ is a smooth quartic double fivefold which is rational, then there exists an equivalence $\mathsf{Ku}(X) \cong D^{b}(\mathrm{Coh}(W))$ for a smooth projective Calabi–Yau threefold $W$.
\end{conjecture}
The conjecture is equivalent to the nonrationality of every smooth quartic double fivefold~\cite{Cheng2025}*{Conjecture 1.2}. Indeed, the three-Calabi--Yau category $\mathsf{Ku}(X)$ satisfies~\cite{Perry2020}*{Lemma 6.9}
\[
HH^{\bullet}(\mathsf{Ku}(X))
\cong \C\oplus\C^{90}[-2]\oplus\C^2[-3]
\oplus\C^{90}[-4]\oplus\C[-6].
\]
In particular, $\dim_{\C}HH^3(\mathsf{Ku}(X))=2$. If $\mathsf{Ku}(X)\simeq D^b(\mathrm{Coh}(W))$ for a smooth projective variety $W$, comparison of Serre functors forces $W$ to be a Calabi--Yau threefold. The Hochschild--Kostant--Rosenberg isomorphism would then give
\[
\dim_{\C}HH^3(W)=2+2h^{1,1}(W)\geq4,
\]
because $W$ is projective.  Theorem~\ref{thm:intro_double_fivefold} establishes the irrationality predicted by the Cheng--Perry--Zhao conjecture. It also provides evidence for the proposed relationship between Hodge atoms and semiorthogonal decompositions discussed in~\cite[Section 2.7]{AtomsMeetSymbols}.

Our final result concerns quartic fivefolds. Pavi{\'c} and Schreieder~\cite{Pavic2023} proved that a very general quartic fivefold in $\PP^6$ does not admit a decomposition of the diagonal and is therefore not retract rational. Here, we prove:

\begin{Mtheorem}[= Theorem~\ref{thm:quartic_fivefold}]
\label{thm:quartic_fivefold_intro}
Every smooth quartic fivefold in $\PP^{6}$ is nonrational.
\end{Mtheorem}

\bigskip

\subsection*{This paper is organized as follows}
Section~\ref{sec:torics} aims to fix notation and introduce an unfamiliar reader to the concept of toric varieties. No original result is presented, and the decision to include it was to produce a self-contained text. Section~\ref{sec:Loc_methods} concerns localization methods in Gromov--Witten theory, with particular focus on toric varieties. We did not set out to recall the basics of Gromov--Witten invariants, which readers can find, for example, in \cites{Kontsevich1994, fulton1997notesstablemapsquantum, CKK-lectures}. The point of this section is to systematically convert localization methods in the toric setting, translating everything into the language of the previous section, and to familiarize the reader with the algorithm we have used to implement the software. Section~\ref{sec:atoms} produces a very concise description of the theory of Hodge atoms and collects the obstruction to rationality results we explore in Section~\ref{sec:Examples}. Finally, we remark that this text has two appendices. Appendix~\ref{app:toricsoftware} provides a brief account of how to use the updated software, while Appendix~\ref{ComputeFan} contains more on the relation between fans and toric varieties.

\subsection*{Acknowledgements}

The authors are deeply indebted to Leonardo F. Cavenaghi for the extensive review of this paper, whose comments substantially improved the exposition. The authors also thank him and Maxim Kontsevich for sharing their five-dimensional argument obtained jointly with the third named author. 
Milena Brockhof was supported by FAPESP grants no. 2024/22950-5, 2025/23347-3. 
G.~Galindo and L.~Katzarkov are supported by the Simons Foundation, grant SFI-MPS-T-Institutes-00007697, and the Ministry of Education and Science of the Republic of Bulgaria, grant DO1-239/10.12.2024. 
  L.~Katzarkov is supported by the Simons Investigators Award (no.003136), NSF FRG grant DMS-2245099, and NSF FRG grant DMS-2255039.
P.~A.~M.~Martins is supported by FAPESP grants no. 2024/07684-7, 2025/22312-1.

\subsection*{AI Disclosure}
The  LLM models Fable 5.1 and GPT-6 Astra were used to convert the software introduced in \cite{alexeev2026atomic}, developed for flag varieties and complete intersections in flag varieties, 
to toric varieties, by systematically changing input to fan combinatorics. Lemma~\ref{lemma: classes} and Proposition~\ref{prop:finite enum} in Appendix~\ref{app:toricsoftware} were obtained as a suggestion from GPT-6 Astra when discussing this conversion. Its proof, as well as all content of this paper, is human-written. Lastly, the figures presented here were generated using Argon, following the authors' instructions regarding the content they were intended to illustrate, and GPT-6 Astra was used to proofread the manuscript once it was written.





\bigskip

\section{Toric Varieties}\label{sec:torics}
Toric varieties are among the most studied spaces in algebraic geometry. A combinatorial structure, termed a \emph{fan} (Definition~\ref{def:fan}), governs their geometry and allows geometric properties to be read off directly from combinatorial data. In this section, we give a brief overview of the theory, focusing on the results needed later for computing Gromov--Witten invariants and quantum multiplication on toric varieties and their complete intersections. No originality is claimed.

\ 

Let $\T = (\C^{\times})^{n}$ be an $n$-dimensional algebraic torus. A \emph{character} of $\T$ is a morphism of algebraic groups $\chi \colon \T \to \C^{\times}$. A \emph{one-parameter subgroup} (1-PS), or \emph{cocharacter}, of $\T$ is a morphism of algebraic groups $\lambda \colon \C^{\times} \to \T$. Their collections form free abelian groups
\[
M := \Hom(\T, \C^{\times}) \simeq \Z^n \qquad \text{and} \qquad N := \Hom(\C^{\times}, \T) \simeq \Z^n
\]
of rank equal to the dimension $n$ of the torus $\T$; these are called the \emph{character lattice} and the \emph{cocharacter lattice} (or lattice of one-parameter subgroups) of $\T$, respectively. Following standard convention, for each $m \in M$ we denote the corresponding character by $\chi^{m}$, and for each $u \in N$ we denote the corresponding cocharacter by $\lambda^{u}$. What this notation is aimed at emphasize is the following. When $\T = (\C^{\times})^n$ with coordinates $t = (t_1, \ldots, t_n)$, a lattice element $m = (m_1, \ldots, m_n) \in \Z^n$ gives $\chi^m(t) = t_1^{m_1} \cdots t_n^{m_n}$, while $u = (u_1, \ldots, u_n) \in \Z^n$ gives $\lambda^u(t) = (t^{u_1}, \ldots, t^{u_n})$ for $t \in \C^{\times}$.

Given a character $\chi^{m}$ and a one-parameter subgroup $\lambda^{u}$, the composition $\chi^{m} \circ \lambda^{u} \colon \C^{\times} \to \C^{\times}$ is an endomorphism of the algebraic group $\C^{\times}$ and is therefore given by $t \mapsto t^l$ for a unique integer $l \in \Z$. This defines a canonical perfect $\Z$-bilinear pairing
\begin{align}
    \langle \cdot , \cdot \rangle \colon M \times N & \longrightarrow \Z \label{eq:pairing-MN} \\
    (m,u) & \longmapsto \langle m, u \rangle := l, \nonumber
\end{align}
which allows us to canonically identify $N$ with $\Hom_{\Z}(M, \Z)$ and $M$ with $\Hom_{\Z}(N, \Z)$. Moreover, there is a canonical group isomorphism
\[
\begin{aligned}
N \otimes_{\Z} \C^{\times} &\xrightarrow{\;\sim\;} \T \\
u \otimes t &\longmapsto \lambda^{u}(t).
\end{aligned}
\]
Henceforth, we write $\T_N := N \otimes_{\Z} \C^{\times}$ for the algebraic torus whose cocharacter lattice is $N$ and whose character lattice is $M$.

\medskip

\subsection{Affine Toric Varieties}
In this subsection, we give a brief introduction to affine toric varieties. They serve as the local building blocks that glue together to form general toric varieties.

\begin{definition}
    An \emph{affine toric variety} is an irreducible affine variety $V$ containing a torus $\T_N$ as a Zariski open subset, such that the group multiplication action of $\T_N$ on itself extends to an algebraic action $\T_N \times V \to V$.
\end{definition}

As an affine scheme, an affine toric variety can be viewed as the spectrum of a commutative ring whose structure arises from an affine semigroup. Recall that an \emph{affine semigroup} is a finitely generated commutative monoid $S$ (written additively, with identity element $0$) that can be embedded in a lattice $M$. Given an affine semigroup $S \subseteq M$, the associated \emph{semigroup algebra} $\C[S]$ is the $\C$-vector space with basis $\{\chi^m \mid m \in S\}$, where $\chi^m$ denotes the formal basis element (or character) corresponding to $m \in S$. Explicitly,
\[
\C[S] := \left\{ \sum_{i=1}^{k} c_i \chi^{m_i} \;\middle|\; k \in \Z_{\geq 0}, \, c_i \in \C, \text{ and } m_i \in S \text{ for all } 1 \leq i \leq k \right\},
\]
with multiplication induced by the semigroup operation in $S$:
\[
\chi^{m} \cdot \chi^{m'} = \chi^{m + m'} \quad \text{for all } m, m' \in S, \quad \text{and} \quad \chi^0 = 1.
\]

Theorem~\ref{thm:semi-group} below identifies affine toric varieties with the spectra of semigroup algebras. 

\begin{theorem}[{\cite{cox2011toric}*{Proposition 1.1.14, Theorem 1.1.17}}]\label{thm:semi-group}
    Let $S \subseteq M$ be an affine semigroup. Then $\Spec \C[S]$ is an affine toric variety whose torus has character lattice $\Z S$. Conversely, if $V$ is an affine toric variety, then $V$ is isomorphic to $\Spec \C[S]$ for some affine semigroup $S \subseteq M$.
\end{theorem}

It follows that the geometry of an affine toric variety is determined up to isomorphism by a choice of affine semigroup $S \subseteq M$. Next, we explain how convex rational polyhedral cones allow us to construct such semigroups systematically using the duality between the character lattice $M$ and the cocharacter lattice $N$. Towards this goal, we remark that tensoring the dual lattices $M$ and $N$ with $\R$ yields dual real vector spaces $M_{\R} = M \otimes_{\Z} \R \simeq \R^n$ and $N_{\R} = N \otimes_{\Z} \R \simeq \R^n$, over which the canonical pairing~\eqref{eq:pairing-MN} extends $\R$-bilinearly to $\langle \cdot, \cdot \rangle \colon M_{\R} \times N_{\R} \to \R$.

\begin{definition}
    A \emph{convex polyhedral cone} in $N_{\R}$ is a set of the form
    \[
    \sigma = \Cone(F) = \left\{ \sum_{u \in F} \lambda_u u \;\middle|\; \lambda_u \geq 0 \right\} \subseteq N_{\R},
    \]
    where $F \subseteq N_{\R}$ is a finite set (with the convention that $\Cone(\emptyset) = \{0\}$). We say that $\sigma$ is \emph{rational} if $F$ can be chosen as a finite subset of the lattice $N$. For brevity, we often refer to a convex polyhedral cone simply as a \emph{polyhedral cone}.
\end{definition}

\begin{definition}
    Given a polyhedral cone $\sigma \subseteq N_{\R}$, its \emph{dual cone} is defined by
    \[
    \sigma^{\vee} = \{ m \in M_{\R} \mid \langle m,u \rangle \geq 0 \text{ for all } u \in \sigma \}.
    \]
    The \emph{orthogonal complement} of $\sigma$ is the subspace defined by
    \[
    \sigma^{\perp} = \{ m \in M_{\R} \mid \langle m,u \rangle = 0 \text{ for all } u \in \sigma \}.
    \]
\end{definition}

By \cite{cox2011toric}*{Proposition 1.2.4}, $\sigma^{\vee}$ is a polyhedral cone in $M_{\R}$ (which is rational whenever $\sigma$ is rational), and under the canonical identification $(M_{\R})^* \simeq N_{\R}$, one has $(\sigma^{\vee})^{\vee} = \sigma$.

For $m \in M_{\R} \setminus \{0\}$, consider the hyperplane
\[
H_m = \{ u \in N_{\R} \mid \langle m,u \rangle = 0 \}
\]
and the closed half-space
\[
H^+_m = \{ u \in N_{\R} \mid \langle m,u \rangle \geq 0 \},
\]
with the convention that $H_0 = H^+_0 = N_{\R}$. If $m_1, \ldots, m_s$ generate $\sigma^{\vee}$, then the duality $(\sigma^{\vee})^{\vee} = \sigma$ implies that $\sigma = H^+_{m_1} \cap \cdots \cap H^+_{m_s}$.

\begin{figure}[htbp]
\centering
\begin{tikzpicture}[>=Stealth, scale=0.95]

    \begin{scope}[shift={(0,0)}]
        \fill[gray!12, rotate=-45] (-2.0,0) rectangle (2.0,2.6);

        \fill[blue!18] (0,0) -- (1.4,2.8) -- (2.8,1.4) -- cycle;

        \draw[->, gray!60] (-1.8,0) -- (3.0,0) node[right, black] {\small $N_{\R}$};
        \draw[->, gray!60] (0,-1.6) -- (0,3.0);

        \draw[thick, red!80!black, dashed] (-1.4,1.4) -- (1.4,-1.4)
            node[below right] {\small $H_m$};
        \node[red!80!black] at (-0.9, -0.7) {\small $H_m^+$};
        \draw[->, red!80!black, thick] (-0.8,0.8) -- (-0.35,1.25);

        \draw[->, thick, blue!70!black] (0,0) -- (1.4,2.8) node[above left] {\small $\rho_2$};
        \draw[->, thick, blue!70!black] (0,0) -- (2.8,1.4) node[above right] {\small $\rho_1$};
        \node[blue!70!black] at (1.35,1.35) {$\sigma$};

        \fill[red!80!black] (0,0) circle (2.5pt) node[below left] {\small $\tau = H_m \cap \sigma = \{0\}$};
    \end{scope}

    \begin{scope}[shift={(6.8,0)}]
        \begin{scope}
            \clip (-1.8,-1.2) rectangle (3.2,3.0);
            \fill[gray!12] (-2,-1) -- (3,1.5) -- (-1.5,3.5) -- cycle;
        \end{scope}

        \fill[blue!18] (0,0) -- (1.4,2.8) -- (2.8,1.4) -- cycle;

        \draw[->, gray!60] (-1.8,0) -- (3.0,0) node[right, black] {\small $N_{\R}$};
        \draw[->, gray!60] (0,-1.6) -- (0,3.0);

        \draw[thick, red!80!black, dashed] (-1.6,-0.8) -- (3.0,1.5)
            node[below right] {\small $H_{m'}$};
        \draw[->, red!80!black, thick] (1.0,0.5) -- (0.7,1.1);

        \draw[->, thick, blue!70!black] (0,0) -- (1.4,2.8) node[above left] {\small $\rho_2$};
        \draw[->, very thick, red!80!black] (0,0) -- (2.8,1.4)
            node[above left, xshift=-2pt, yshift=2pt] {\small $\tau = H_{m'} \cap \sigma = \rho_1$};
        \node[blue!70!black] at (1.35,1.35) {$\sigma$};

        \fill[black] (0,0) circle (2pt);
    \end{scope}

\end{tikzpicture}
\caption{Faces of a strongly convex polyhedral cone $\sigma \subseteq N_{\R}$. Left: choosing $m$ in the interior of $\sigma^{\vee}$ yields the $0$-dimensional face $\tau = H_m \cap \sigma = \{0\}$. Right: choosing $m' \in \sigma^{\vee} \cap \rho_1^{\perp}$ yields the $1$-dimensional face (both a ray and a facet) $\tau = H_{m'} \cap \sigma = \rho_1$.}
\label{fig:faces-supporting-hyperplanes}
\end{figure}
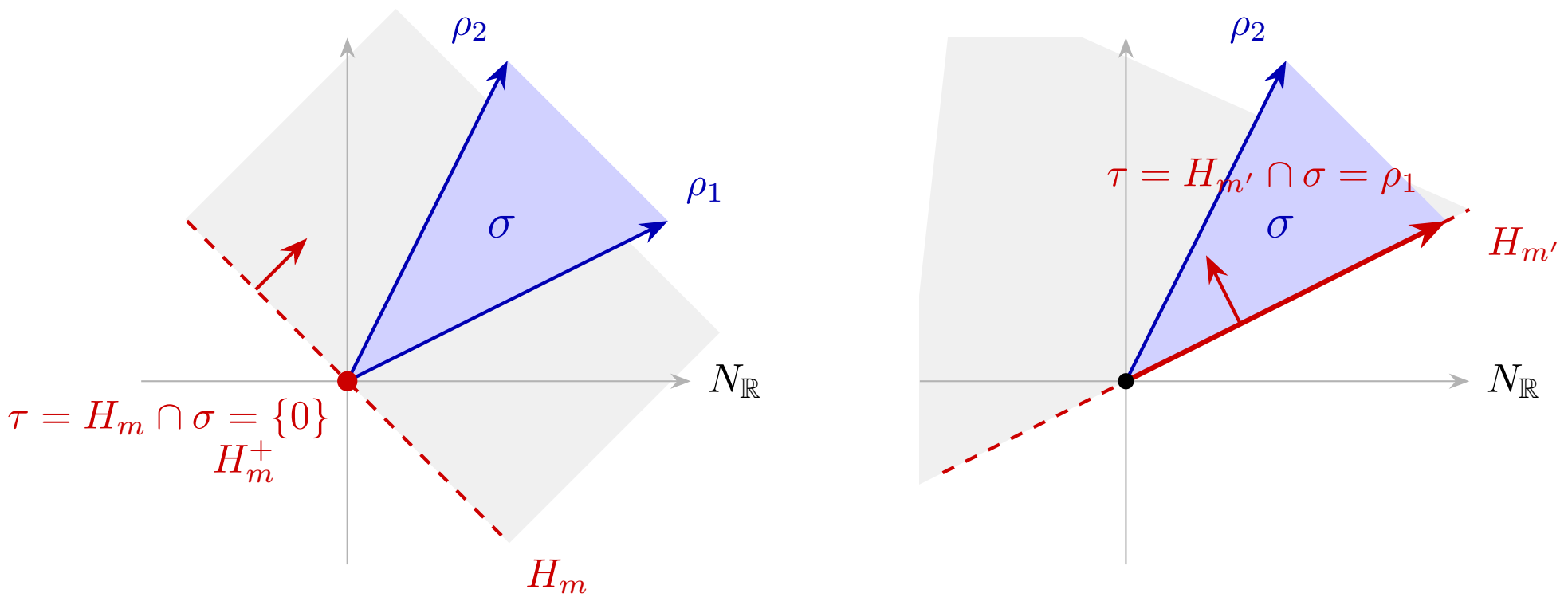

\begin{definition}
    A \emph{face} of a polyhedral cone $\sigma$ is a subset of the form $\tau = H_m \cap \sigma$ for some $m \in \sigma^{\vee}$, denoted $\tau \preceq \sigma$ (see Figure~\ref{fig:faces-supporting-hyperplanes}). Every face $\tau$ of $\sigma$ is itself a polyhedral cone. The \emph{dimension} $\dim \sigma$ is the dimension of the smallest linear subspace $\Span(\sigma) \subseteq N_{\R}$ containing $\sigma$. A face $\tau \preceq \sigma$ is called a \emph{facet} of $\sigma$ if it has codimension $1$ in $\sigma$, and a \emph{ray} if $\dim \tau = 1$.
\end{definition}

\begin{proposition}[{\cite{cox2011toric}*{Proposition 1.2.12}}]\label{prop:equiv}
    Let $\sigma \subseteq N_{\R} \simeq \R^n$ be a polyhedral cone. Then the following conditions are equivalent:
    \begin{itemize}
        \item[(a)] $\{0\}$ is a face of $\sigma$;
        \item[(b)] $\sigma$ contains no positive-dimensional subspace of $N_{\R}$;
        \item[(c)] $\sigma \cap (-\sigma) = \{0\}$;
        \item[(d)] $\dim \sigma^{\vee} = n$.
    \end{itemize}
\end{proposition}

\begin{definition}
    A polyhedral cone $\sigma \subseteq N_{\R}$ satisfying the equivalent conditions of the Proposition~\ref{prop:equiv} is called \emph{strongly convex}.
\end{definition}

Henceforth, unless stated otherwise, we assume that all cones under consideration are strongly convex rational polyhedral cones. In this setting, $\sigma$ admits a canonical minimal generating set. Specifically, every ray $\rho$ of $\sigma$ is a rational half-line, so the semigroup $\rho \cap N$ is generated by a unique primitive lattice element $u_{\rho} \in \rho \cap N$, called the \emph{ray generator} of $\rho$. By \cite{cox2011toric}*{Lemma 1.2.15}, $\sigma$ is generated by the ray generators of its rays:
\[
\sigma = \Cone(u_{\rho} \mid \rho \text{ is a ray of } \sigma).
\]

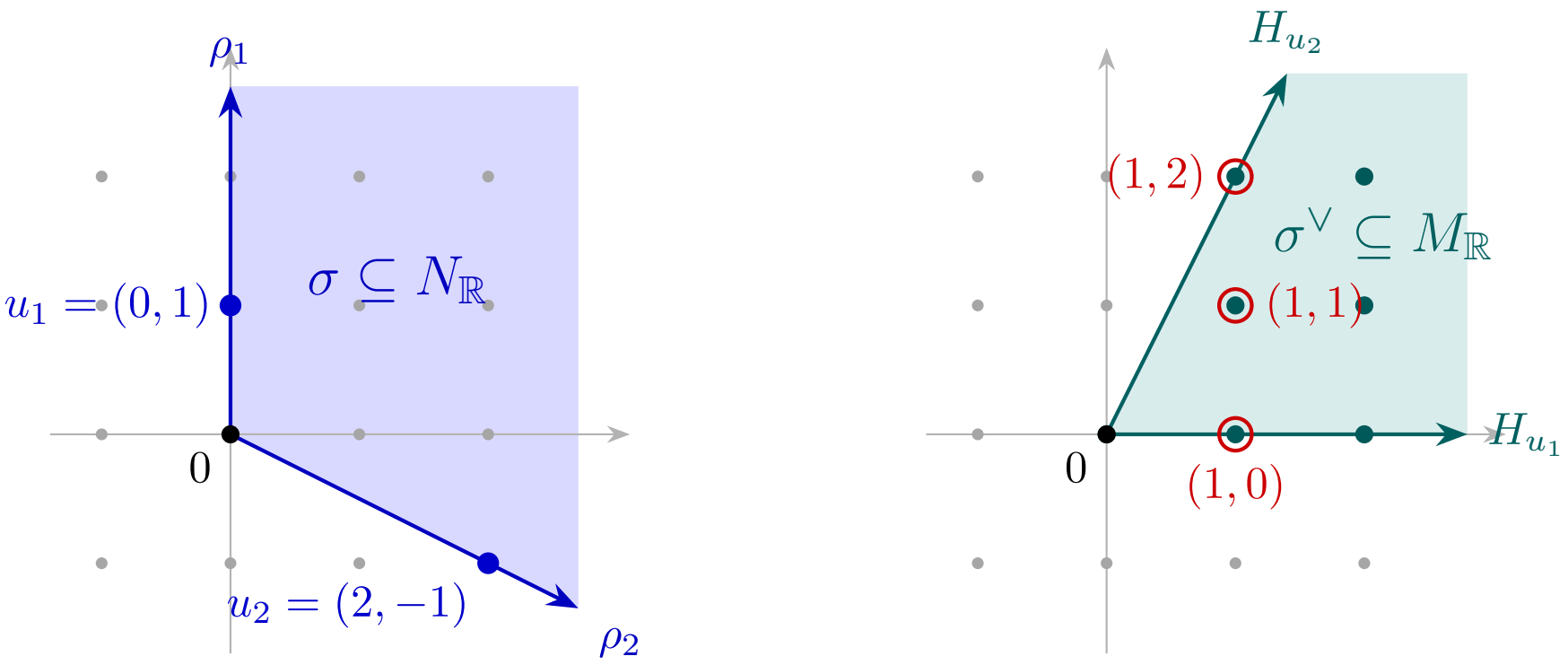
\begin{figure}[htbp]
\centering
\begin{tikzpicture}[>=Stealth, scale=0.95]

    \begin{scope}[shift={(0,0)}]
        \fill[blue!15] (0,0) -- (0,2.7) -- (2.7,2.7) -- (2.7,-1.35) -- cycle;

        \draw[->, gray!60] (-1.4,0) -- (3.1,0);
        \draw[->, gray!60] (0,-1.7) -- (0,3.0);

        \foreach \x in {-1,...,2} {
            \foreach \y in {-1,...,2} {
                \fill[gray!70] (\x,\y) circle (1.3pt);
            }
        }

        \draw[->, thick, blue!75!black] (0,0) -- (0,2.7) node[above] {\small $\rho_1$};
        \draw[->, thick, blue!75!black] (0,0) -- (2.7,-1.35) node[below right] {\small $\rho_2$};

        \fill[blue!80!black] (0,1) circle (2.4pt) node[left] {\small $u_1 = (0,1)$};
        \fill[blue!80!black] (2,-1) circle (2.4pt) node[below left] {\small $u_2 = (2,-1)$};
        \fill[black] (0,0) circle (2pt) node[below left] {\small $0$};

        \node[blue!75!black] at (1.3,1.2) {$\sigma \subseteq N_{\R}$};
    \end{scope}

    \begin{scope}[shift={(6.8,0)}]
        \fill[teal!15] (0,0) -- (2.8,0) -- (2.8,2.8) -- (1.4,2.8) -- cycle;

        \draw[->, gray!60] (-1.4,0) -- (3.1,0);
        \draw[->, gray!60] (0,-1.7) -- (0,3.0);

        \foreach \x in {-1,...,2} {
            \foreach \y in {-1,...,2} {
                \fill[gray!70] (\x,\y) circle (1.3pt);
            }
        }

        \draw[->, thick, teal!75!black] (0,0) -- (2.8,0) node[right] {\small $H_{u_1}$};
        \draw[->, thick, teal!75!black] (0,0) -- (1.4,2.8) node[above] {\small $H_{u_2}$};

        \foreach \x/\y in {0/0, 1/0, 2/0, 1/1, 2/1, 1/2, 2/2} {
            \fill[teal!70!black] (\x,\y) circle (2pt);
        }

        \draw[red!80!black, thick] (1,0) circle (3.6pt) node[below, yshift=-2pt] {\small $(1,0)$};
        \draw[red!80!black, thick] (1,1) circle (3.6pt) node[right, xshift=2pt] {\small $(1,1)$};
        \draw[red!80!black, thick] (1,2) circle (3.6pt) node[left, xshift=-2pt] {\small $(1,2)$};
        \fill[black] (0,0) circle (2pt) node[below left] {\small $0$};

        \node[teal!75!black] at (2.15,1.55) {$\sigma^{\vee} \subseteq M_{\R}$};
    \end{scope}

\end{tikzpicture}
\caption{Left: the strongly convex rational polyhedral cone $\sigma = \Cone(u_1, u_2) \subseteq N_{\R}$ with primitive ray generators $u_1 = (0,1)$ and $u_2 = (2,-1)$. Right: its dual cone $\sigma^{\vee} = \Cone((1,0), (1,2)) \subseteq M_{\R}$ and the lattice points of $S_{\sigma} = \sigma^{\vee} \cap M$, whose minimal semigroup generators $(1,0)$, $(1,1)$, and $(1,2)$ are circled.}
\label{fig:cone-and-dual-cone}
\end{figure}

Now, given a rational polyhedral cone $\sigma \subseteq N_{\R}$, the lattice points of its \emph{dual cone} $\sigma^{\vee} \subseteq M_{\R}$,
\[
S_{\sigma} := \sigma^{\vee} \cap M \subseteq M,
\]
form a finitely generated commutative monoid (i.e., an affine semigroup) by \emph{Gordan's Lemma} \cite[Proposition 1.2.17]{cox2011toric}, as illustrated in Figure~\ref{fig:cone-and-dual-cone}. Consequently, $S_{\sigma}$ defines an affine toric variety.

\begin{theorem}[{\cite{cox2011toric}*{Theorem 1.2.18}}]
    Let $\sigma \subseteq N_{\R}$ be a rational polyhedral cone with associated affine semigroup $S_{\sigma} = \sigma^{\vee} \cap M$. Then
    \[
    U_{\sigma} := \Spec \C[S_{\sigma}] = \Spec \C[\sigma^{\vee} \cap M]
    \]
    is an affine toric variety. Furthermore, $\dim U_{\sigma} = n \Longleftrightarrow$ the torus of $U_{\sigma}$ equals $\T_N \Longleftrightarrow \sigma$ is strongly convex.
\end{theorem}

\begin{remark}
    If $\tau \preceq \sigma$ is a face of a rational polyhedral cone $\sigma$, then $U_{\tau}$ embeds naturally as a principal affine open subset of $U_{\sigma}$. As we explain in the next subsection, these open immersions allow us to glue the affine charts $U_{\sigma}$ along common faces to construct general toric varieties.
\end{remark}

\medskip

\subsection{General Properties of Toric Varieties}

We now turn to general and smooth projective toric varieties, which can be constructed by gluing affine toric varieties along suitable open subsets. Since a strongly convex rational polyhedral cone describes every normal affine toric variety, the combinatorial data governing this gluing is encoded by a \emph{fan}; namely, a finite collection of cones satisfying compatibility conditions along their intersections.

\begin{definition}
    A \emph{toric variety} is an irreducible variety $X$ containing a torus $\T_N$ as a Zariski open subset such that the group multiplication action of $\T_N$ on itself extends to an algebraic action $\T_N \times X \to X$. If $X$ is a projective variety, we call $X$ a \emph{projective toric variety}.
\end{definition}

\begin{definition}\label{def:fan}
    A \emph{fan} $\Sigma$ in $N_{\R}$ is a finite collection of strongly convex rational polyhedral cones $\sigma \subseteq N_{\R}$ such that:
    \begin{enumerate}
        \item[(1)] for all $\sigma \in \Sigma$, each face of $\sigma$ is also in $\Sigma$;
        \item[(2)] for all $\sigma_1, \sigma_2 \in \Sigma$, the intersection $\sigma_1 \cap \sigma_2$ is a face of both $\sigma_1$ and $\sigma_2$.
    \end{enumerate}
    Furthermore, the \emph{support} of $\Sigma$ is the subset $|\Sigma| := \bigcup_{\sigma \in \Sigma} \sigma \subseteq N_{\R}$, and $\Sigma(r)$ denotes the set of $r$-dimensional cones of $\Sigma$.
\end{definition}

The fan $\Sigma$ provides the data needed to glue the affine pieces $\{U_{\sigma}\}_{\sigma \in \Sigma}$ together into an abstract normal separated toric variety $X_{\Sigma}$, by gluing $U_{\sigma_1}$ and $U_{\sigma_2}$ along the open subvariety $U_{\tau}$ corresponding to their common face $\tau = \sigma_1 \cap \sigma_2$. In what follows, we state the key algebraic results on semigroup algebras that justify this gluing and refer the reader to \cite{cox2011toric} for complete proofs.

\begin{figure}[htbp]
\centering
\begin{tikzpicture}[>=Stealth, scale=0.92]

    \begin{scope}[shift={(0,0)}]
        \fill[blue!15] (0,0) rectangle (2.4,2.4);
        \fill[orange!15] (-2.4,0) rectangle (0,2.4);

        \draw[->, gray!60] (-2.7,0) -- (2.7,0) node[right, black] {\small $N_{\R}$};
        \draw[->, gray!60] (0,-1.1) -- (0,2.7);

        \foreach \x in {-2,...,2} {
            \foreach \y in {0,...,2} {
                \fill[gray!70] (\x,\y) circle (1.3pt);
            }
        }

        \draw[->, thick, blue!75!black] (0,0) -- (2.4,0) node[below left] {\small $e_1$};
        \draw[->, thick, orange!80!black] (0,0) -- (-2.4,0) node[below right] {\small $-e_1$};
        \draw[dashed, thick, red!80!black] (0,-1.0) -- (0,2.6);
        \draw[->, very thick, red!80!black] (0,0) -- (0,2.4) node[above] {\small $\tau = \sigma_1 \cap \sigma_2 = H_m \cap \sigma_i$};

        \fill[black] (0,0) circle (2pt) node[below left] {\small $0$};
        \node[blue!75!black] at (1.25,1.3) {$\sigma_1 \subseteq H_m^+$};
        \node[orange!80!black] at (-1.25,1.3) {$\sigma_2 \subseteq H_{-m}^+$};
    \end{scope}

    \begin{scope}[shift={(6.8,0)}]
        \fill[blue!12] (0,0) rectangle (2.4,2.4);
        \fill[orange!12] (-2.4,0) rectangle (0,2.4);

        \draw[->, gray!60] (-2.7,0) -- (2.7,0) node[right, black] {\small $M_{\R}$};
        \draw[->, gray!60] (0,-1.1) -- (0,2.7);

        \foreach \x in {-2,...,2} {
            \foreach \y in {0,...,2} {
                \fill[teal!70!black] (\x,\y) circle (1.7pt);
            }
        }

        \draw[->, thick, blue!75!black] (0,0) -- (2.4,0);
        \draw[->, thick, orange!80!black] (0,0) -- (-2.4,0);
        \draw[->, thick, teal!80!black] (0,0) -- (0,2.4);

        \fill[red!80!black] (1,0) circle (2.4pt) node[below] {\small $m = e_1$};
        \fill[red!80!black] (-1,0) circle (2.4pt) node[below] {\small $-m = -e_1$};
        \fill[black] (0,0) circle (2pt);

        \node[blue!75!black] at (1.35,1.5) {$\sigma_1^{\vee}$};
        \node[orange!80!black] at (-1.35,1.5) {$\sigma_2^{\vee}$};
        \node[teal!80!black] at (0, -0.75) {\small $\tau^{\vee} = \sigma_1^{\vee} + \sigma_2^{\vee} \implies S_{\tau} = S_{\sigma_1} + S_{\sigma_2}$};
    \end{scope}

\end{tikzpicture}
\caption{Left: two cones $\sigma_1 = \Cone(e_1,e_2)$ and $\sigma_2 = \Cone(-e_1,e_2)$ in $N_{\R}$ separated by $H_m$ ($m = e_1$) along their common face $\tau = \sigma_1 \cap \sigma_2 = \Cone(e_2)$. Right: in $M_{\R}$, localizing $S_{\sigma_1}$ at $m$ and $S_{\sigma_2}$ at $-m$ yields $S_{\tau} = S_{\sigma_1} + S_{\sigma_2}$, gluing $U_{\sigma_1} \simeq \C^2$ and $U_{\sigma_2} \simeq \C^2$ along $U_{\tau} \simeq \C^{\times} \times \C$ to form $\mathbb{P}^1 \times \C$.}
\label{fig:gluing-separation-lemma}
\end{figure}
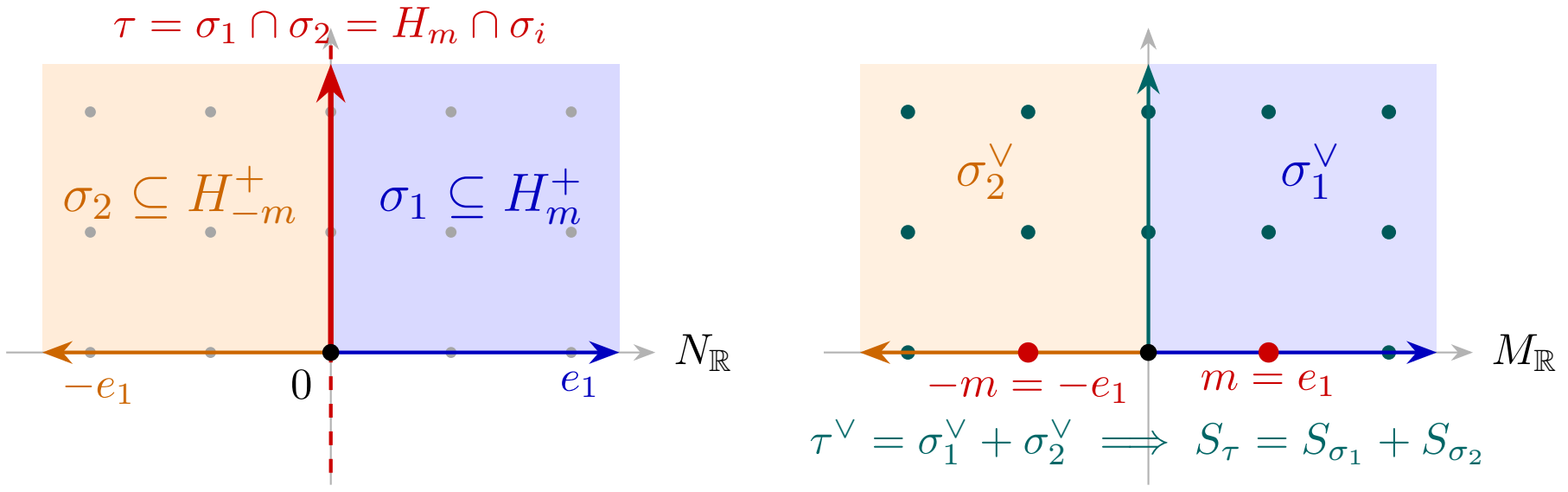

The first ingredient relates the semigroup algebras of a cone $\sigma$ and a face $\tau \preceq \sigma$. 
Recall that any face $\tau \preceq \sigma$ can be written as $\tau = \sigma \cap H_m$ for some lattice element $m \in \sigma^{\vee}$. Then \cite[Proposition 1.3.16]{cox2011toric} shows that $S_{\tau} = \tau^{\vee} \cap M = S_{\sigma} + \Z_{\geq 0}(-m)$, so the semigroup algebra $\C[S_{\tau}]$ is the localization of $\C[S_{\sigma}]$ at the character $\chi^m \in \C[S_{\sigma}]$, that is, $\C[S_{\tau}] = \C[S_{\sigma}]_{\chi^m}$. Passing to spectra yields the principal affine open subset
\[
U_{\tau} = (U_{\sigma})_{\chi^m} \subseteq U_{\sigma}.
\]
Furthermore, for any two cones $\sigma_1, \sigma_2 \in \Sigma$ intersecting along the common face $\tau = \sigma_1 \cap \sigma_2$, then there exists a lattice element $m \in \sigma_1^{\vee} \cap (-\sigma_2)^{\vee} \cap M = S_{\sigma_1} \cap (-S_{\sigma_2})$ such that
\[
\sigma_1 \cap H_m = \tau = \sigma_2 \cap H_m,
\]
with $\sigma_1 \subseteq H_m^+$ and $\sigma_2 \subseteq H_{-m}^+$ (see Figure~\ref{fig:gluing-separation-lemma}), see \cite[Lemma 1.2.13]{cox2011toric}. Combining this with the localization formula above gives
\[
U_{\sigma_1} \supseteq (U_{\sigma_1})_{\chi^m} = U_{\tau} = (U_{\sigma_2})_{\chi^{-m}} \subseteq U_{\sigma_2}.
\]
The canonical identification of both localizations with $\C[S_{\tau}]$ induces a gluing isomorphism
\[
g_{\sigma_2, \sigma_1} \colon (U_{\sigma_1})_{\chi^m} \xrightarrow{\;\sim\;} (U_{\sigma_2})_{\chi^{-m}}
\]
that equals the identity on $U_{\tau}$, ensuring that the cocycle compatibility conditions on triple overlaps are automatically satisfied. Moreover, by \cite[Proposition 3.1.3]{cox2011toric}, the dual cones satisfy $\tau^{\vee} = \sigma_1^{\vee} + \sigma_2^{\vee}$ and in fact $S_{\tau} = S_{\sigma_1} + S_{\sigma_2}$, which implies that the diagonal morphism $U_{\tau} \to U_{\sigma_1} \times U_{\sigma_2}$ is a closed immersion and hence that the resulting toric variety $X_{\Sigma}$ is separated \cite[Theorem 3.1.5]{cox2011toric}.

With this construction in hand, we can describe how the torus $\T_N$ acts on $X_{\Sigma}$ via its invariant affine charts $U_{\sigma} = \Spec \C[S_{\sigma}]$. Because the coordinate ring of $U_{\sigma}$ is $\C[U_{\sigma}] = \C[S_{\sigma}]$, specifying an algebraic action $\T_N \times U_{\sigma} \to U_{\sigma}$ is equivalent to giving a $\C$-algebra coaction homomorphism
\[
\C[S_{\sigma}] = \C[U_{\sigma}] \longrightarrow \C[\T_N \times U_{\sigma}] \simeq \C[\T_N] \otimes_{\C} \C[U_{\sigma}] = \C[M] \otimes_{\C} \C[S_{\sigma}],
\]
which on basis elements is defined by $\chi^m \longmapsto \chi^m \otimes \chi^m$ for all $m \in S_{\sigma}$.

Geometrically, by \cite[Proposition 1.3.1]{cox2011toric}, the closed points of $U_{\sigma} = \Spec \C[S_{\sigma}]$ are in natural bijection with semigroup homomorphisms $\gamma \colon (S_{\sigma}, +) \to (\C, \cdot)$\footnote{here $\C$ is viewed as a commutative monoid under multiplication with identity $1$, so $\gamma(0) = 1$ and $\gamma(m + m') = \gamma(m)\gamma(m')$}. Explicitly, a point $p \in U_{\sigma}$ corresponds to the evaluation map $m \mapsto \chi^m(p)$, and an element $t \in \T_N$ acts on $\gamma$ by
\[
(t \cdot \gamma)(m) = \chi^m(t)\gamma(m) \quad \text{for all } m \in S_{\sigma}.
\]
Because these actions are compatible with the localization maps $\C[S_{\sigma}] \to \C[S_{\tau}] = \C[S_{\sigma}]_{\chi^m}$ for $\tau \preceq \sigma$, the torus actions on the charts $U_{\sigma}$ glue to a global algebraic action of $\T_N$ on $X_{\Sigma}$.

Beyond providing a means of constructing toric varieties, fans directly encode their geometric properties. Since our primary interest lies in smooth compact (and projective) toric varieties, we highlight the corresponding combinatorial conditions on the fan $\Sigma$.

\begin{definition}
    Let $\sigma \subseteq N_{\R}$ be a strongly convex rational polyhedral cone and let $\Sigma$ be a fan in $N_{\R}$.
    \begin{enumerate}
        \item[(a)] The cone $\sigma$ is \emph{smooth} (or \emph{regular}) if its minimal ray generators form part of a $\Z$-basis of $N$. The fan $\Sigma$ is \emph{smooth} if every cone $\sigma \in \Sigma$ is smooth.
        \item[(b)] The fan $\Sigma$ is \emph{complete} if its support fills the entire ambient space, i.e., $|\Sigma| = N_{\R}$.
    \end{enumerate}
\end{definition}

\begin{remark}
    Terminology for smooth cones varies in the literature; many authors refer to them as \emph{regular} cones.
\end{remark}

\begin{theorem}[{\cite[Theorem 3.1.19]{cox2011toric}}]
    Let $X_{\Sigma}$ be the toric variety defined by a fan $\Sigma$ in $N_{\R}$. Then:
    \begin{enumerate}
        \item[(a)] $X_{\Sigma}$ is a smooth variety if and only if the fan $\Sigma$ is smooth.
        \item[(b)] $X_{\Sigma}$ is compact in the classical topology (equivalently, complete as an algebraic variety) if and only if the fan $\Sigma$ is complete.
    \end{enumerate}
\end{theorem}

\ 

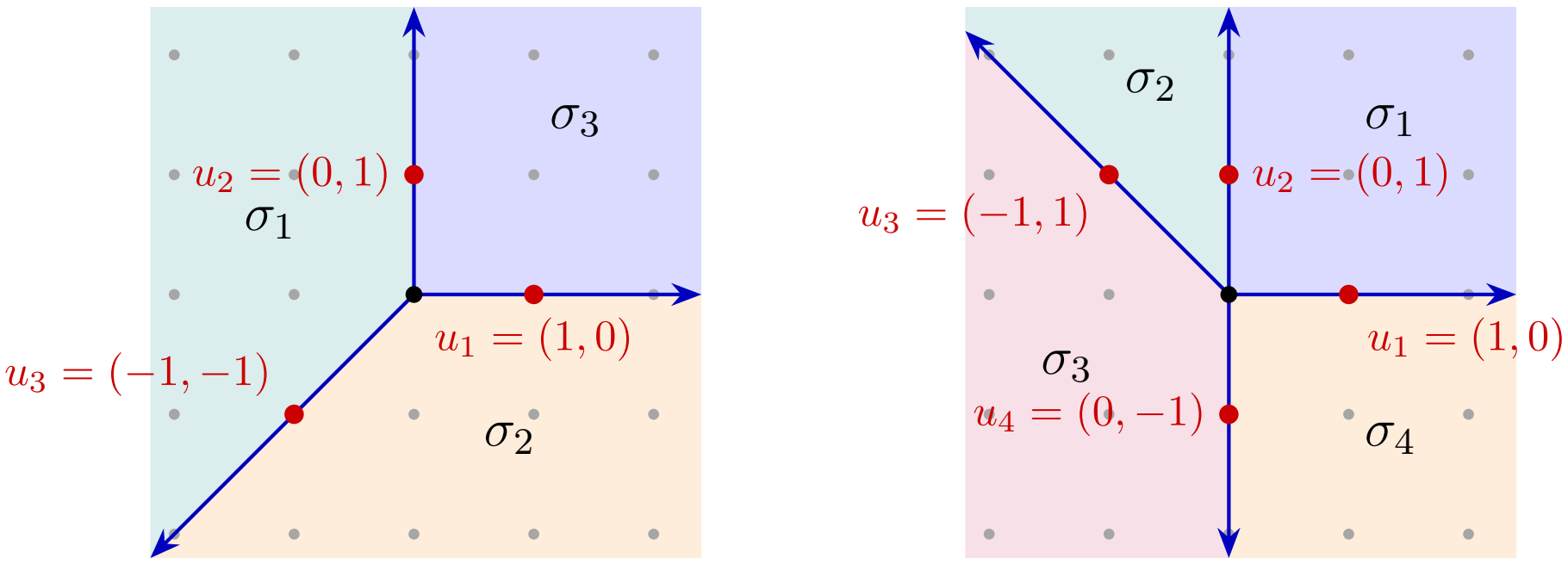
\begin{figure}[htbp]
\centering
\begin{tikzpicture}[>=Stealth, scale=0.92]

    \begin{scope}[shift={(0,0)}]
        \fill[blue!14] (0,0) -- (2.4,0) -- (2.4,2.4) -- (0,2.4) -- cycle;
        \fill[teal!14] (0,0) -- (0,2.4) -- (-2.2,2.4) -- (-2.2,-2.2) -- cycle;
        \fill[orange!14] (0,0) -- (2.4,0) -- (2.4,-2.2) -- (-2.2,-2.2) -- cycle;

        \foreach \x in {-2,...,2} {
            \foreach \y in {-2,...,2} {
                \fill[gray!70] (\x,\y) circle (1.3pt);
            }
        }

        \draw[->, thick, blue!75!black] (0,0) -- (2.4,0);
        \draw[->, thick, blue!75!black] (0,0) -- (0,2.4);
        \draw[->, thick, blue!75!black] (0,0) -- (-2.2,-2.2);

        \fill[red!80!black] (1,0) circle (2.3pt) node[below, yshift=-1pt] {\small $u_1 = (1,0)$};
        \fill[red!80!black] (0,1) circle (2.3pt) node[left, xshift=-1pt] {\small $u_2 = (0,1)$};
        \fill[red!80!black] (-1,-1) circle (2.3pt) node[above left, xshift=-1pt] {\small $u_3 = (-1,-1)$};
        \fill[black] (0,0) circle (2pt);

        \node at (1.35,1.45) {$\sigma_3$};
        \node at (-1.2,0.6) {$\sigma_1$};
        \node at (0.8,-1.2) {$\sigma_2$};
    \end{scope}

    \begin{scope}[shift={(6.8,0)}]
        \fill[blue!14] (0,0) -- (2.4,0) -- (2.4,2.4) -- (0,2.4) -- cycle;
        \fill[teal!14] (0,0) -- (0,2.4) -- (-2.2,2.4) -- (-2.2,2.2) -- cycle;
        \fill[purple!12] (0,0) -- (-2.2,2.2) -- (-2.2,-2.2) -- (0,-2.2) -- cycle;
        \fill[orange!14] (0,0) -- (2.4,0) -- (2.4,-2.2) -- (0,-2.2) -- cycle;

        \foreach \x in {-2,...,2} {
            \foreach \y in {-2,...,2} {
                \fill[gray!70] (\x,\y) circle (1.3pt);
            }
        }

        \draw[->, thick, blue!75!black] (0,0) -- (2.4,0);
        \draw[->, thick, blue!75!black] (0,0) -- (0,2.4);
        \draw[->, thick, blue!75!black] (0,0) -- (-2.2,2.2);
        \draw[->, thick, blue!75!black] (0,0) -- (0,-2.2);

        \fill[red!80!black] (1,0) circle (2.3pt) node[below right, yshift=-1pt] {\small $u_1 = (1,0)$};
        \fill[red!80!black] (0,1) circle (2.3pt) node[right, xshift=1pt] {\small $u_2 = (0,1)$};
        \fill[red!80!black] (-1,1) circle (2.3pt) node[below left] {\small $u_3 = (-1,1)$};
        \fill[red!80!black] (0,-1) circle (2.3pt) node[left, xshift=-1pt] {\small $u_4 = (0,-1)$};
        \fill[black] (0,0) circle (2pt);

        \node at (1.35,1.45) {$\sigma_1$};
        \node at (-0.65,1.75) {$\sigma_2$};
        \node at (-1.35,-0.6) {$\sigma_3$};
        \node at (1.35,-1.2) {$\sigma_4$};
    \end{scope}

\end{tikzpicture}
\caption{Examples of smooth complete fans in $N_{\R} \simeq \R^2$. Left: the fan of $\mathbb{P}^2$, with three maximal cones $\sigma_0, \sigma_1, \sigma_2$ and ray generators $u_1 = e_1$, $u_2 = e_2$, and $u_3 = -e_1-e_2$. Right: the fan of the Hirzebruch surface $\mathcal{H}_1 \simeq \operatorname{Bl}_{\mathrm{pt}}\mathbb{P}^2$, with four maximal cones $\sigma_1, \sigma_2, \sigma_3, \sigma_4$. In both fans, the ray generators of every $2$-dimensional cone form a $\Z$-basis of $N = \Z^2$ ($\det = \pm 1$), and $|\Sigma| = \R^2$.}
\label{fig:smooth-complete-fans}
\end{figure}

A natural question is whether every toric variety can be constructed from a fan. As a consequence Sumihiro's theorem \cite[Theorem 3.1.7]{cox2011toric}, every normal separated toric variety comes from a fan; in particular, since smooth projective varieties are normal and separated, the following is true:

\begin{proposition}[{\cite[Corollary 3.1.8]{cox2011toric}}]\label{Suh}
    Let $X$ be a smooth projective toric variety with torus $\T_N$. Then there exists a smooth complete fan $\Sigma$ in $N_{\R}$ such that $X \simeq X_{\Sigma}$.
\end{proposition}

Having described the torus action on the affine charts $U_{\sigma}$, we next observe that each chart contains a canonical distinguished point. The \emph{Orbit--Cone Correspondence} (Theorem~\ref{OrbitCone} below) will then establish that the $\T_N$-orbits of these distinguished points decompose $X_{\Sigma}$ and allow us to read off the orbit structure directly from the fan $\Sigma$.

\begin{proposition}[{\cite[Proposition 1.3.2 and Corollary 1.3.3]{cox2011toric}}]\label{fix}
    Let $X_{\Sigma}$ be the toric variety associated to a fan $\Sigma$ in $N_{\R}$. For each cone $\sigma \in \Sigma$, the affine chart $U_{\sigma} \subseteq X_{\Sigma}$ contains a distinguished point $\gamma_{\sigma} \in U_{\sigma}$ defined by the semigroup homomorphism $S_{\sigma} \to \C$ given by
    \[
    m \longmapsto
    \begin{cases}
        1 & \text{if } m \in S_{\sigma} \cap \sigma^{\perp} = \sigma^{\perp} \cap M, \\
        0 & \text{otherwise}.
    \end{cases}
    \]
    Furthermore, $\gamma_{\sigma}$ is a fixed point of the $\T_N$-action on $X_{\Sigma}$ if and only if $\dim \sigma = \dim N_{\R}$.
\end{proposition}

For each cone $\sigma \in \Sigma$, we denote the $\T_N$-orbit of the distinguished point $\gamma_{\sigma}$ by
\[
O(\sigma) := \T_N \cdot \gamma_{\sigma} \subseteq X_{\Sigma}.
\]

\begin{theorem}[The Orbit--Cone Correspondence {\cite[Theorem 3.2.6]{cox2011toric}}]\label{OrbitCone}
    Let $X_{\Sigma}$ be the toric variety of a fan $\Sigma$ in $N_{\R}$, where $n = \dim N_{\R}$. Then:
    \begin{itemize}
        \item[(a)] There is an inclusion-reversing bijective correspondence
        \[
        \begin{aligned}
            \{\text{cones } \sigma \in \Sigma\} &\longleftrightarrow \{\T_N\text{-orbits in } X_{\Sigma}\} \\
            \sigma &\longleftrightarrow O(\sigma) \simeq \Hom_{\Z}(\sigma^{\perp} \cap M, \C^{\times}).
        \end{aligned}
        \]
        \item[(b)] For each cone $\sigma \in \Sigma$, $\dim O(\sigma) = n - \dim \sigma$.
        \item[(c)] For each cone $\sigma \in \Sigma$, the affine open chart $U_{\sigma}$ is the disjoint union of orbits
        \[
        U_{\sigma} = \bigcup_{\tau \preceq \sigma} O(\tau).
        \]
        \item[(d)] A cone $\tau \in \Sigma$ is a face of $\sigma \in \Sigma$ (i.e., $\tau \preceq \sigma$) if and only if $O(\sigma) \subseteq \overline{O(\tau)}$, and
        \[
        \overline{O(\tau)} = \bigcup_{\substack{\sigma \in \Sigma \\ \tau \preceq \sigma}} O(\sigma),
        \]
        where $\overline{O(\tau)}$ denotes the closure of $O(\tau)$ in $X_{\Sigma}$ in both the classical and Zariski topologies.
    \end{itemize}
\end{theorem}

We conclude that the orbits of the torus action on $X_{\Sigma}$ are completely governed by the combinatorics of the fan $\Sigma$; in particular, the torus fixed points correspond bijectively to the $n$-dimensional cones $\sigma \in \Sigma(n)$ (which are the maximal cones of $\Sigma$ whenever $\Sigma$ is complete). A detailed proof of Theorem~\ref{OrbitCone} can be found in \cite[Theorem 3.2.6]{cox2011toric}. For details on how to recover the fan from the torus orbits, see Appendix~\ref{ComputeFan}.

\medskip

\subsection{Cohomology from fan combinatorics}
\label{sec:picard}

Our next goal is to describe the cohomology ring of a smooth projective toric variety in terms of its fan. According to Theorem~\ref{OrbitCone}, each ray $\rho \in \Sigma(1)$ corresponds to a $(n-1)$-dimensional torus orbit $O(\rho)$, whose closure $D_\rho = \overline{O(\rho)}$ is a $\T_N$-invariant prime divisor on $X_\Sigma$. The cohomology ring of $X_\Sigma$ is generated by the classes of these $\T_N$-invariant divisors modulo combinatorial relations read directly from $\Sigma$.

\begin{definition}
Let $\Sigma$ be a fan in $N_{\R}$ with rays $\Sigma(1) = \{\rho_1,\ldots,\rho_r\}$, and let $u_i \in N$ denote the primitive ray generator of $\rho_i$. Write $D_i = D_{\rho_i} = \overline{O(\rho_i)}$ for the corresponding $\T_N$-invariant divisor.

In the polynomial ring $\mathbb{Q}[x_1,\ldots,x_r]$, graded by $\deg(x_i) = 1$, we associate the variable $x_i$ to the ray $\rho_i$ and define the \emph{Stanley--Reisner ideal} of $\Sigma$ by
\[
\mathscr{J}
=
\left\langle
x_{i_1}\cdots x_{i_s}
\;\middle|\;
\begin{array}{l}
i_1,\ldots,i_s \text{ are distinct and}\\
\operatorname{Cone}(u_{i_1},\ldots,u_{i_s})
\notin\Sigma
\end{array}
\right\rangle.
\]
Let $\mathcal{I} \subseteq \mathbb{Q}[x_1,\ldots,x_r]$ be the homogeneous ideal generated by the linear forms
\[
\mathcal{I}
=
\left\langle
\sum_{i=1}^{r}\langle m,u_i\rangle x_i
\;\middle|\;
m\in M
\right\rangle.
\]
Since both $\mathcal{I}$ and $\mathscr{J}$ are homogeneous ideals, the quotient ring
\[
R_{\mathbb{Q}}(\Sigma)
=
\mathbb{Q}[x_1,\ldots,x_r]/(\mathcal{I}+\mathscr{J})
=
\bigoplus_{k \geq 0} R_{\mathbb{Q}}(\Sigma)_k
\]
is a graded $\mathbb{Q}$-algebra, called the \emph{rational cohomology ring of the fan} $\Sigma$.
\end{definition}

\begin{theorem}[\cite{cox2011toric}*{Theorem 12.4.1}]
\label{12.4.1}
Let $X_\Sigma$ be a smooth projective toric variety. Then the homomorphism of graded $\mathbb{Q}$-algebras
\[
R_{\mathbb{Q}}(\Sigma)
\longrightarrow
H^{\ast}(X_\Sigma,\mathbb{Q})
\]
induced by sending $x_i$ to the divisor class $[D_i] \in H^2(X_\Sigma,\mathbb{Q})$ is an isomorphism that doubles degrees. Consequently, for every $k \geq 0$,
\[
H^{2k}(X_\Sigma,\mathbb{Q}) \simeq R_{\mathbb{Q}}(\Sigma)_k
\qquad \text{and} \qquad
H^{2k+1}(X_\Sigma,\mathbb{Q}) = 0.
\]
\end{theorem}

\begin{example}[Projective spaces $\PP^n$]
Let $X_\Sigma = \PP^n$, whose fan $\Sigma$ in $N_{\R} = \R^n$ has rays $\rho_1,\ldots,\rho_{n+1}$ with primitive generators $u_i = e_i$ for $i=1,\ldots,n$ and $u_{n+1} = -\sum_{i=1}^n e_i$, where $e_1,\ldots,e_n$ denotes the standard $\Z$-basis of $N = \Z^n$. To each ray $\rho_i$ we associate the hyperplane divisor $D_i = \overline{O(\rho_i)}$, whose cohomology class $[D_i] \in H^2(\PP^n,\mathbb{Q})$ corresponds to the image of $x_i$ in $R_{\mathbb{Q}}(\Sigma)$ under Theorem~\ref{12.4.1}. The maximal cones of $\Sigma$ are
\[
\sigma_j
=
\operatorname{Cone}(u_i\mid i\neq j),
\qquad j=1,\ldots,n+1.
\]
Every proper subset of $\{u_1,\ldots,u_{n+1}\}$ is contained in $\{u_i \mid i \neq j\}$ for some $j$ and therefore spans a cone in $\Sigma$; the only subset of $\{u_1,\ldots,u_{n+1}\}$ that does not span a cone in $\Sigma$ is $\{u_1,\ldots,u_{n+1}\}$ itself. Evaluating the linear forms in $\mathcal{I}$ on the dual $\Z$-basis $e_1^{\ast},\ldots,e_n^{\ast}$ of $M$ yields $\sum_{i=1}^{n+1} \langle e_j^{\ast}, u_i \rangle x_i = x_j - x_{n+1}$ for $j=1,\ldots,n$. Hence the ideals $\mathscr{J}$ and $\mathcal{I}$ are given by
\[
\mathscr{J}
=
\langle x_1\cdots x_{n+1}\rangle
\qquad \text{and} \qquad
\mathcal{I}
=
\langle x_1-x_{n+1},\ldots,x_n-x_{n+1}\rangle.
\]
Thus, we obtain
\[
H^{\ast}(\PP^n,\mathbb{Q})
\simeq
\frac{\mathbb{Q}[x_1,\ldots,x_{n+1}]}
{\langle x_1\cdots x_{n+1},\,
x_1-x_{n+1},\ldots,x_n-x_{n+1}\rangle}
\simeq
\frac{\mathbb{Q}[x_{n+1}]}{\langle x_{n+1}^{n+1}\rangle}.
\]
\end{example}

Next, we compile several structural properties of the cohomology of smooth projective toric varieties. Throughout the remainder of this section, let $X_\Sigma$ be a smooth projective toric variety of dimension $n$. We retain the notation $\Sigma(1)=\{\rho_1,\ldots,\rho_r\}$, with primitive ray generators $u_1,\ldots,u_r\in N$ and associated $\T_N$-invariant divisors $D_1,\ldots,D_r$, and we identify $x_i\in R_{\mathbb{Q}}(\Sigma)$ with its image $[D_i] \in H^2(X_\Sigma,\mathbb{Q})$ under the isomorphism of Theorem~\ref{12.4.1}.

Let $A_d(X_\Sigma)$ denote the Chow group of $d$-dimensional algebraic cycles on $X_\Sigma$ modulo rational equivalence. Since $X_\Sigma$ is smooth of dimension $n$, setting $A^k(X_\Sigma):= A_{n-k}(X_\Sigma)$ equips $A^{\ast}(X_\Sigma) = \bigoplus_{k=0}^n A^k(X_\Sigma)$ with the intersection product, making it a graded commutative ring. Following Theorem~\ref{OrbitCone}, for each cone $\sigma\in\Sigma$, the orbit closure $\overline{O(\sigma)}$ is a subvariety of codimension $\dim\sigma$ in $X_\Sigma$.

\begin{lemma}[\cite{cox2011toric}*{Lemma 12.5.1}]
The rational equivalence classes $[\overline{O(\sigma)}]\in A^{\dim\sigma}(X_\Sigma)$ generate $A^{\ast}(X_\Sigma)$ as an abelian group as $\sigma$ runs through $\Sigma$.
\end{lemma}

Every cycle class $\alpha\in A^k(X_\Sigma)=A_{n-k}(X_\Sigma)$ determines a homology class in $H_{2n-2k}(X_\Sigma,\mathbb{Z})$. Sending $\alpha$ to the Poincar\'e dual of this homology class defines the cycle class ring homomorphism
\[
A^{\ast}(X_\Sigma)_{\mathbb{Q}} := A^{\ast}(X_\Sigma) \otimes_{\mathbb{Z}} \mathbb{Q}
\longrightarrow
H^{2\bullet}(X_\Sigma,\mathbb{Q}) = H^{\ast}(X_\Sigma,\mathbb{Q}).
\]
When $X_\Sigma$ is smooth and projective, this map is an isomorphism \cite{cox2011toric}*{Theorem 12.5.3}. Consequently, every rational cohomology class on $X_\Sigma$ is algebraic (and in particular is a Hodge class of pure type $(k,k)$), and by the lemma above, the cohomology classes $[\overline{O(\sigma)}]$ of the $\T_N$-orbit closures generate the entire cohomology ring $H^{\ast}(X_\Sigma,\mathbb{Q})$. This property will allow us to simplify the expression of the A-model F-bundle in Section~\ref{def:A-model-F-bundle}.

As our main goal is to apply these methods to study Gromov--Witten invariants, we require a toric description of the effective curves on $X_\Sigma$. Theorem~\ref{Toric Cone} below provides an explicit description of the cone of effective curves in terms of $\T_N$-invariant curves, which correspond to the $(n-1)$-dimensional cones (walls) in $\Sigma(n-1)$. In what follows, we also characterize the first Chern class of a line bundle in terms of $\T_N$-invariant divisors and compute the intersection pairings between these divisors and the invariant curves. This will be especially useful when considering complete intersections later.

Let $N_1(X_\Sigma)$ denote the finite-dimensional real vector space of numerical equivalence classes of $1$-cycles on $X_\Sigma$. The \emph{cone of effective curves} is defined by
\[
\mathsf{NE}(X_\Sigma)
=
\left\{
\sum_{i=1}^{s}a_i[C_i]
\;\middle|\;
a_i\in\mathbb{R}_{\geq0},\
C_i\subset X_\Sigma
\text{ irreducible complete curves}
\right\}
\subseteq N_1(X_\Sigma).
\]
Its closure $\overline{\mathsf{NE}}(X_\Sigma) \subseteq N_1(X_\Sigma)$ is called the \emph{Mori cone} of $X_\Sigma$.

\begin{theorem}[\cite{cox2011toric}*{Theorem 6.3.20 and Theorem 6.3.22}] \label{Toric Cone}
Let $X_\Sigma$ be a smooth projective toric variety of dimension $n$. Then $\overline{\mathsf{NE}}(X_\Sigma) = \mathsf{NE}(X_\Sigma)$ is a strongly convex rational polyhedral cone in $N_1(X_\Sigma)$. Furthermore,
\[
\overline{\mathsf{NE}}(X_\Sigma)
=
\sum_{\tau\in\Sigma(n-1)}
\mathbb{R}_{\geq0}[\overline{O(\tau)}],
\]
where $[\overline{O(\tau)}] \in N_1(X_\Sigma)$ denotes the numerical equivalence class of the $\T_N$-invariant curve $\overline{O(\tau)} \simeq \PP^1$.
\end{theorem}

\begin{remark}
The space $N_1(X_\Sigma)$ and the Néron--Severi space $N^1(X_\Sigma)$ of numerical equivalence classes of Cartier divisors on $X_\Sigma$ are dual real vector spaces under the intersection pairing $([D], [C]) \mapsto D \cdot C = \int_C [D]$. Moreover, for a smooth projective toric variety, numerical and linear equivalence of divisors coincide, so that
\[
N^1(X_\Sigma)
=
\Pic(X_\Sigma)_{\mathbb{R}}
=
\Pic(X_\Sigma)\otimes_{\mathbb{Z}}\mathbb{R}
\simeq
H^2(X_\Sigma,\mathbb{R}).
\]
\end{remark}

Next, we express the first Chern class in terms of the $\T_N$-invariant divisors $D_i$ and their associated cohomology classes $x_i = [D_i]$.

\begin{proposition}[\cite{cox2011toric}*{Theorem 8.2.3 and Example 13.1.1}]
If $L = \cO_{X_\Sigma}(D)$ is the line bundle associated with a Cartier divisor $D$ on $X_\Sigma$, then its first Chern class is $c_1(L) = [D] \in H^2(X_\Sigma,\mathbb{Q})$. In particular, since the canonical divisor of $X_\Sigma$ is given by $K_{X_\Sigma} = -\sum_{i=1}^r D_i$, the first Chern class of the tangent bundle $TX_\Sigma$ is
\[
c_1(TX_\Sigma)
=
c_1(\cO_{X_\Sigma}(-K_{X_\Sigma}))
=
\sum_{i=1}^r [D_i]
=
\sum_{i=1}^r x_i.
\]
\end{proposition}

Our final aim in this section is to describe the linear relation satisfied by the ray generators of two adjacent maximal cones sharing a wall $\tau \in \Sigma(n-1)$. This relation allows us to compute the homology classes and intersection numbers of $\T_N$-invariant curves directly from the fan.

\begin{lemma}[\cite{cox2011toric}*{Wall Relation (6.4.4)}]\label{lem:wall}
Let $X_\Sigma$ be a smooth projective toric variety of dimension $n$. Let $\tau\in\Sigma(n-1)$ be a wall, and set
\[
I_\tau
=
\{i\in\{1,\ldots,r\}\mid \rho_i\subseteq\tau\}.
\]
Let $\sigma$ and $\sigma'$ be the two maximal cones in $\Sigma(n)$ containing $\tau$ as a facet, written as
\[
\sigma=\operatorname{Cone}(u_i \mid i \in I_\tau \cup \{p\})
\qquad \text{and} \qquad
\sigma'=\operatorname{Cone}(u_i \mid i \in I_\tau \cup \{q\})
\]
for distinct indices $p, q \in \{1,\ldots,r\} \setminus I_\tau$. Then there exist unique integers $a_i \in \Z$ for $i\in I_\tau$ such that
\begin{equation}\label{eq:wall}
u_p+u_q+\sum_{i\in I_\tau}a_i u_i=0.
\end{equation}
\end{lemma}

\begin{lemma}[\cite{cox2011toric}*{Proposition 6.4.4}]\label{lem:wallintersection}
With the notation and assumptions of Lemma~\ref{lem:wall}, let $C_\tau=\overline{O(\tau)}$ be the $\T_N$-invariant curve corresponding to $\tau\in\Sigma(n-1)$, and set
\[
b_i=
\begin{cases}
1, & \text{if } i=p \text{ or } i=q,\\
a_i, & \text{if } i\in I_\tau,\\
0, & \text{otherwise}.
\end{cases}
\]
Then the intersection number of $x_i = [D_i]$ with $C_\tau$ is given by
\[
\int_{C_\tau}x_i = D_i \cdot C_\tau = b_i
\]
for every $i=1,\ldots,r$.
\end{lemma}

\medskip

 \subsection{The Picard group, ampleness, and complete intersections}
\label{sec:toricpicard}

We conclude this discussion by recalling key structural results concerning the Picard group of a smooth projective toric variety and the positivity of line bundles. These results allow us to characterize ampleness, the Fano property, and numerical invariants of complete intersections directly from the combinatorial data of the fan.

Throughout this section, let $X_\Sigma$ be a smooth projective toric variety of dimension $n$, with rays $\Sigma(1)=\{\rho_1,\ldots,\rho_r\}$ and primitive ray generators $u_1,\ldots,u_r\in N$. As before, $D_i$ denotes the $\T_N$-invariant prime divisor associated with $\rho_i$, and $x_i = [D_i] = c_1(\mathcal{O}_{X_\Sigma}(D_i)) \in H^2(X_\Sigma,\mathbb{Q})$.

We begin by recalling the relationship between the Picard group and the second cohomology group.

\begin{theorem}[\cite{cox2011toric}*{Theorem 12.3.2}]
The first Chern class induces an isomorphism
\[
\operatorname{Pic}(X_\Sigma)
\xrightarrow{\ \sim\ }
H^2(X_\Sigma,\mathbb{Z}).
\]
\end{theorem}

Furthermore, because $X_\Sigma$ is smooth and complete, every Weil divisor is Cartier and linearly equivalent to a $\T_N$-invariant divisor. This yields the fundamental short exact sequence of free abelian groups
\[
0 \longrightarrow M \xrightarrow{\quad} \mathbb{Z}^r \longrightarrow \operatorname{Pic}(X_\Sigma) \longrightarrow 0,
\]
where the first map sends $m \in M$ to the tuple $(\langle m,u_1\rangle,\ldots, \langle m,u_r\rangle)$. Consequently, the Picard group admits the explicit presentation
\[
\operatorname{Pic}(X_\Sigma)
\cong
\frac{\mathbb{Z}^r}
{\left\{
(\langle m,u_i\rangle)_{i=1}^r
\;\middle|\; m\in M
\right\}}.
\]
This description enables us to compute the second homology group $H_2(X_\Sigma,\mathbb{Z})$ explicitly. Since $H^{\ast}(X_\Sigma, \mathbb{Z})$ is free and finitely generated, the Universal Coefficient Theorem provides a perfect pairing between $H^2(X_\Sigma, \mathbb{Z})$ and $H_2(X_\Sigma, \mathbb{Z})$, giving
\[
H_2(X_\Sigma, \mathbb{Z}) \simeq \operatorname{Hom}\bigl(H^2(X_\Sigma, \mathbb{Z}),\mathbb{Z}\bigr) \simeq \operatorname{Hom}\bigl(\operatorname{Pic}(X_\Sigma),\mathbb{Z}\bigr).
\]
Applying $\operatorname{Hom}_{\mathbb{Z}}(-, \mathbb{Z})$ to the presentation exact sequence above yields the exact sequence
\[
0 \longrightarrow \operatorname{Hom}\bigl(\operatorname{Pic}(X_\Sigma),\mathbb{Z}\bigr)
\longrightarrow \mathbb{Z}^r
\xrightarrow{\quad} \operatorname{Hom}(M,\mathbb{Z}) = N,
\]
where the rightmost map sends the standard basis vector $e_i \in \mathbb{Z}^r$ to the ray generator $u_i \in N$, identifying $\operatorname{Hom}(\operatorname{Pic}(X_\Sigma), \mathbb{Z})$ with its kernel.

\begin{proposition}\label{prop:torich2}
There is a canonical isomorphism
\[
H_2(X_\Sigma,\mathbb{Z})
\cong
\left\{
(a_i)_{i=1}^r\in\mathbb{Z}^r
\;\middle|\;
\sum_{i=1}^r a_i u_i=0
\text{ in } N
\right\}.
\]
Under this identification, a homology class $\beta\in H_2(X_\Sigma,\mathbb{Z})$ corresponding to $(a_i)_{i=1}^r$ satisfies
\begin{equation}\label{Intbetatoric}
\int_\beta x_i = a_i,
\qquad
\text{and more generally}\qquad
\int_\beta\left(\sum_{i=1}^r c_i x_i\right)
=\sum_{i=1}^r c_i a_i
\end{equation}
for any $c_1,\ldots,c_r \in \mathbb{Q}$. Moreover, the class $[C_\tau] \in H_2(X_\Sigma,\mathbb{Z})$ of a torus-invariant curve $C_\tau = \overline{O(\tau)}$ for $\tau \in \Sigma(n-1)$ corresponds precisely to the coefficient tuple $(b_i)_{i=1}^r$ given in Lemma~\ref{lem:wallintersection}.
\end{proposition}

Theorem~\ref{OrbitCone} implies that each torus-invariant curve $C_\tau$ for $\tau \in \Sigma(n-1)$ can be written as the transversal intersection of the divisors corresponding to the rays containing $\tau$:
\[
C_\tau = \bigcap_{i \in I_\tau} D_i,
\qquad \text{where } I_\tau = \{i \in \{1,\ldots,r\} \mid \rho_i \subseteq \tau\}.
\]
Because $X_\Sigma$ is smooth, the normal bundle splits as a direct sum of line bundles whose degrees are given by the intersection numbers:
\[
N_{C_\tau / X_\Sigma} \simeq \bigoplus_{i \in I_\tau} \mathcal{O}_{X_\Sigma}(D_i)\big|_{C_\tau},
\qquad
\deg \left(\mathcal{O}_{X_\Sigma}(D_i)\big|_{C_\tau} \right) = D_i \cdot C_\tau = b_i.
\]
Combining this with the wall relations (Lemma~\ref{lem:wall}), the anticanonical degree of $C_\tau$ evaluates to
\[
\int_{C_\tau} c_1(TX_\Sigma)
=
\sum_{i=1}^r b_i
=
2+\sum_{i\in I_\tau}a_i.
\]
This intersection number need not be positive; when $2 + \sum_{i \in I_\tau} a_i > 0$ for all walls $\tau \in \Sigma(n-1)$, the toric variety $X_\Sigma$ is Fano.

\medskip

We now turn to the ampleness of line bundles. Since $X_\Sigma$ is smooth and projective, the notions of ampleness and very ampleness coincide (\cite{cox2011toric}*{Theorem 6.1.15}). The Toric Kleiman Criterion (Theorem~\ref{Toric Kleiman}) enables us to characterize ampleness and nefness entirely in terms of intersection pairings with the torus-invariant curves.

Recall that the Néron--Severi space $N^1(X_\Sigma) = \operatorname{Pic}(X_\Sigma) \otimes_{\mathbb{Z}} \mathbb{R}$ and the real vector space of $1$-cycles $N_1(X_\Sigma)$ are dual real vector spaces under the intersection pairing. The Mori cone $\overline{\mathsf{NE}}(X_\Sigma) \subseteq N_1(X_\Sigma)$ generated by effective curves is dual to the Nef cone $\operatorname{Nef}(X_\Sigma) \subseteq N^1(X_\Sigma)$. By Theorem~\ref{Toric Cone}, this cone is generated by finitely many invariant curves $C_\tau$, yielding
\[
\operatorname{Nef}(X_\Sigma)
=
\left\{
[D]\in N^1(X_\Sigma)
\;\middle|\;
D \cdot [\overline{O(\tau)}]\geq 0
\text{ for every }
\tau\in\Sigma(n-1)
\right\}.
\]

\begin{proposition}[\cite{cox2011toric}*{Theorem 6.3.12}] \label{Toric KleimanNef}
Let $D \in \operatorname{Pic}(X_\Sigma) \cong H^2(X_\Sigma,\mathbb{Z})$ be a Cartier divisor on $X_\Sigma$. The following conditions are equivalent:
\begin{enumerate}
    \item[(a)] $D$ is nef;
    \item[(b)] $D$ is generated by global sections;
    \item[(c)] $D \cdot [\overline{O(\tau)}] \geq 0$ for every $\tau \in \Sigma(n-1)$.
\end{enumerate}
\end{proposition}

Because the ample cone is the interior of the nef cone, ampleness corresponds to strict positivity on all nonzero classes in the Mori cone.

\begin{theorem}[Toric Kleiman Criterion, \cite{cox2011toric}*{Theorem 6.3.22}] \label{Toric Kleiman}
Let $D = \sum_{i=1}^r d_i D_i$ be a Cartier divisor on $X_\Sigma$. Then:
\begin{enumerate}
    \item[(a)] $D$ is nef if and only if
    \[
    \sum_{i=1}^r d_i b_i(\tau) \geq 0
    \quad
    \text{for every } \tau \in \Sigma(n-1).
    \]
    \item[(b)] $D$ is ample if and only if
    \[
    \sum_{i=1}^r d_i b_i(\tau) > 0
    \quad
    \text{for every } \tau \in \Sigma(n-1),
    \]
\end{enumerate}
where $b_i(\tau)$ denotes the intersection coefficient $D_i \cdot C_\tau$ computed via Lemma~\ref{lem:wallintersection}.
\end{theorem}

\medskip

We now apply these criteria to smooth complete intersections in smooth projective toric varieties.

\begin{theorem}\label{thm:ci-toric}
Let $X_\Sigma$ be a smooth projective toric variety of dimension $n$, and let
\[
\mathcal{E}
=
\bigoplus_{j=1}^e
L_j,
\qquad
L_j=\mathcal{O}_{X_\Sigma} \left(\sum_{i=1}^r d_{j,i}D_i \right),
\quad d_{j,i}\in\mathbb{Z}.
\]
Let $Y \subset X_\Sigma$ be the zero locus of a transversal global section of $\mathcal{E}$, assumed smooth of codimension $e$, and let $i \colon Y \hookrightarrow X_\Sigma$ denote the inclusion. Then:
\begin{enumerate}
\item[(a)] If every line bundle $L_j$ is ample (i.e., satisfies $\sum_{i=1}^r d_{j,i}b_i(\tau) > 0$ for all $j=1,\ldots,e$ and $\tau\in\Sigma(n-1)$) and $n-e \geq 3$, then the restriction map induces an isomorphism
\[
\operatorname{Pic}(X_\Sigma)
\xrightarrow{\ \sim\ }
\operatorname{Pic}(Y).
\]

\item[(b)] The canonical bundle and first Chern class of $Y$ are given by
\[
K_Y
=
i^*\mathcal{O}_{X_\Sigma}
\left(
\sum_{j=1}^e \sum_{i=1}^r d_{j,i}D_i-\sum_{i=1}^r D_i
\right)
\quad \text{and} \quad
c_1(TY)
=
i^*\left(
\sum_{i=1}^r x_i
-
\sum_{j=1}^e\sum_{i=1}^r d_{j,i}x_i
\right).
\]

\item[(c)] If the anticanonical bundle of $Y$ is ample, which is guaranteed when
\[
\sum_{i=1}^r
\left(1-\sum_{j=1}^e d_{j,i}\right)
b_i(\tau)>0
\]
for every $\tau\in\Sigma(n-1)$, then $Y$ is a Fano manifold.
\end{enumerate}
\end{theorem}

\begin{proof}
Item (a) follows from the Lefschetz hyperplane theorem, because $Y$ is a smooth complete intersection of ample divisors of codimension $e \le n-3$. 

For (b), the adjunction formula for vector bundle zero loci gives $K_Y = i^*(K_{X_\Sigma} \otimes \det \mathcal{E})$. Substituting the canonical divisor $K_{X_\Sigma} = \mathcal{O}_{X_\Sigma}(-\sum_{i=1}^r D_i)$ yields the desired expression.

For (c), the anticanonical bundle is $-K_Y = i^*\mathcal{O}_{X_\Sigma}(\sum_{i=1}^r D_i - \sum_{j=1}^e \sum_{i=1}^r d_{j,i}D_i)$. By Theorem~\ref{Toric Kleiman}, the assumed positivity condition ensures that this divisor is ample on $X_\Sigma$, so its restriction to $Y$ is ample, from where it follows that $Y$ is Fano.
\end{proof}



\bigskip

\section{Gromov--Witten invariants via localization}\label{sec:Loc_methods}

In~\cite{Kontsevich1995}, Kontsevich introduced a method for computing Gromov--Witten invariants of smooth projective varieties via torus localization, reducing the intersection theory on moduli spaces of stable maps to combinatorial calculations over graphs. This approach was subsequently extended to general Goresky--Kottwitz--MacPherson (GKM) spaces. Building on these ideas, several authors have developed explicit localization formulas and computational tools to perform such calculations.

Our goal in this section is to recall the genus-zero localization formulas established in~\cites{Liu_Chiu,holmes2025computationsequivariantgromovwittentheory,liu2013localization} and to explain how all of their geometric ingredients can be expressed in terms of the combinatorial data of a toric fan. We also review the corresponding localization formula for smooth complete intersections in toric varieties, following~\cite{alexeev2026atomic}, in which the contribution of the defining vector bundle is incorporated via an equivariant twisting factor. Although no new localization formulas are introduced here, we include this review to document explicitly how the theoretical results are translated into the algorithms implemented in the software we have developed.

\subsection{Fixed points, invariant curves, and tangent weights}

Let $X_\Sigma$ be a smooth projective toric variety of dimension $n$. We keep employing the notation of Section~\ref{sec:torics}: the rays of the fan are $\Sigma(1)=\{\rho_1,\ldots,\rho_r\}$, with primitive generators $u_i\in N$, and $\mathbb{T}_N$ denotes the acting algebraic torus. The main goal in this section is that of expressing the genus-zero localization formula in terms of these combinatorial data.

The first ingredients are the torus-fixed points and the one-dimensional torus orbits. The first basic ingredient comes from Theorem~\ref{OrbitCone}, which associates a $\mathbb{T}_N$-orbit $O(\sigma)$ with each cone $\sigma\in\Sigma$.

\begin{proposition}\label{prop:toricGKMgraph}
The fixed points of the $\mathbb{T}_N$-action on $X_\Sigma$ are in bijection with the maximal cones $\sigma\in\Sigma(n)$. Denote the fixed point $\gamma_\sigma$ associated with $\sigma$ (see Proposition~\ref{fix}) by $p_\sigma$.

Two distinct fixed points $p_\sigma$ and $p_{\sigma'}$ are joined by a $\mathbb{T}_N$-invariant curve if and only if $\tau=\sigma\cap\sigma'\in\Sigma(n-1)$. In this case,
\[
C_\tau=\overline{O(\tau)}
=O(\tau)\cup\{p_\sigma,p_{\sigma'}\}
\cong\mathbb{P}^1.
\]
In particular, distinct cones in $\Sigma(n-1)$ determine distinct $\mathbb{T}_N$-invariant curves.
\end{proposition}

Consequently, smooth projective toric varieties are GKM spaces, and the localization formulas of~\cites{Liu_Chiu,holmes2025computationsequivariantgromovwittentheory,liu2013localization} apply directly.

For $\sigma\in\Sigma(n)$, set
\[
I_\sigma=\{i\in\{1,\ldots,r\}\mid \rho_i\subseteq\sigma\},
\]
and for each $i\in I_\sigma$, let
\[
\tau_i=\operatorname{Cone}
(u_j\mid j\in I_\sigma\setminus\{i\}) \in \Sigma(n-1)
\]
be the facet of $\sigma$ obtained by omitting the ray $\rho_i$. Since $X_\Sigma$ is smooth, $\{u_i\}_{i\in I_\sigma}$ forms a $\mathbb{Z}$-basis of $N$.

\begin{proposition}\label{prop:toricTangentWeights}
For every $\sigma\in\Sigma(n)$, the tangent space at $p_\sigma$ decomposes as a $\mathbb{T}_N$-module into one-dimensional weight spaces
\[
T_{p_\sigma}X_\Sigma
\cong
\bigoplus_{i\in I_\sigma}
T_{p_\sigma}C_{\tau_{i}}.
\]
Moreover, if $\tau=\sigma\cap\sigma'\in\Sigma(n-1)$, then the tangent weight $w(\sigma, \tau)$ of $C_\tau$ at $p_\sigma$ is the unique irreducible element of
\[
(S_\sigma\setminus\{0\})\cap\tau^\perp,
\qquad \text{where } S_\sigma=\sigma^\vee\cap M.
\]
At the opposite fixed point $p_{\sigma'}$, the tangent weight of $C_\tau$ is $w(\sigma', \tau)=-w(\sigma, \tau)$.
\end{proposition}

\begin{proof}
By Proposition~\ref{fix}, the maximal ideal of $\C[U_\sigma] = \C[S_\sigma]$ corresponding to the fixed point $p_\sigma$ is $\fm_{p_\sigma}=\langle \chi^m \mid m\in S_\sigma\setminus\{0\}\rangle$, so that
\[
\bigoplus_{\substack{m\in S_\sigma\\
m\text{ is irreducible}}}
\C\overline{\chi^m}
\cong \fm_{p_\sigma}/\fm_{p_\sigma}^2
\simeq T_{p_\sigma}^*X_\Sigma.
\]
Since $X_\Sigma$ is smooth, the irreducible elements of $S_\sigma$ form the $\mathbb{Z}$-basis $\{m_{\sigma,i}\}_{i\in I_\sigma}$ of $M$ dual to the primitive ray generators $\{u_i\}_{i\in I_\sigma}$ of $\sigma$. The corresponding characters $\chi^{m_{\sigma, i}}$ provide $\mathbb{T}_N$-eigenfunctions that serve as local coordinates centered at $p_\sigma$, and their weights determine the one-dimensional weight spaces of $T_{p_\sigma}X_\Sigma$.

To identify the weight on the tangent line $T_{p_\sigma}C_{\tau_i}\subset T_{p_\sigma}X_\Sigma$, note that $\cO_{C_{\tau_i},p_\sigma} = \cO_{X_\Sigma,p_\sigma}/I_{C_{\tau_i},p_\sigma}$, where $I_{C_{\tau_i},p_\sigma} = I_{C_{\tau_i}}\cdot\cO_{X_\Sigma,p_\sigma}$ is the stalk at $p_\sigma$ of the vanishing ideal sheaf $I_{C_{\tau_i}}$ of $C_{\tau_i}$. This yields
\[
T_{p_\sigma}^*C_{\tau_i}
=
\fm_{C_{\tau_i},p_\sigma}/
\fm_{C_{\tau_i},p_\sigma}^2
\cong
\frac{\fm_{p_\sigma}}
{\fm_{p_\sigma}^2+I_{C_{\tau_i},p_\sigma}}.
\]
Because $C_{\tau_i}=\overline{O(\tau_i)}$, on the affine chart $U_\sigma$ the ideal of regular functions vanishing along $C_{\tau_i}$ is generated by
\[
\{\chi^m\mid
m\in S_\sigma\setminus\tau_i^\perp\}.
\]
Thus, the tangent line $T_{p_\sigma}C_{\tau_i}\subset T_{p_\sigma}X_\Sigma$ corresponds to the unique basis character $\chi^m$ with $m \in S_\sigma$ satisfying $\langle m,u\rangle=0$ for all $u\in\tau_i$. Since $\tau_i$ is obtained by omitting $\rho_i$ from $\sigma$, this character is $\chi^{m_{\sigma,i}}$, where $m_{\sigma,i} \in M$ is the dual basis element defined by
\[
\langle m_{\sigma,i},u_j\rangle
=\delta_{ij}
\qquad \text{for all } j\in I_\sigma.
\]
Consequently, $w(\sigma,\tau_i) = m_{\sigma,i}$ is the unique irreducible element of $(S_\sigma\setminus\{0\})\cap\tau_i^\perp$, and
\[
T_{p_\sigma}X_\Sigma
=
\bigoplus_{i\in I_\sigma}
T_{p_\sigma}C_{\tau_i}.
\]
Finally, if $\sigma' \in \Sigma(n)$ is the other maximal cone containing $\tau_i$, the $\mathbb{T}_N$-action on $C_{\tau_i}\cong\mathbb{P}^1$ has opposite tangent weights at its two fixed points $p_\sigma$ and $p_{\sigma'}$, so that $w(\sigma',\tau_i) = -w(\sigma,\tau_i)$.
\end{proof}

We also require the weights of the normal bundle $N_{C_\tau/X_\Sigma} \simeq \bigoplus_{i\in I_\tau} \cO_{X_\Sigma}(D_i)|_{C_\tau}$ at the two fixed points of $C_\tau$. For each $i\in I_\tau$, denote by $w_i^\perp(\sigma,\tau)$ the weight of the normal direction corresponding to $\rho_i$ at $p_\sigma$. These are precisely the weights of the fibers of $\mathcal{O}_{X_\Sigma}(D_i)|_{C_\tau}$ at $p_\sigma$, and are characterized by
\[
\langle w_i^\perp(\sigma,\tau),u_j\rangle
=\delta_{ij},
\qquad j\in I_\sigma.
\]
Similarly, we denote by $w_i^\perp(\sigma',\tau)$ the corresponding normal fiber weights at $p_{\sigma'}$. Note that when $\tau = \tau_k$ for $k \in I_\sigma$, we have $w_i^\perp(\sigma,\tau_k) = m_{\sigma,i} = w(\sigma,\tau_i)$ for every $i \in I_{\tau_k} = I_\sigma \setminus \{k\}$.

This description allows us to recover the degrees of the normal line bundle summands from the weights at the two endpoints:
\[
\deg\bigl(\mathcal{O}_{X_\Sigma}(D_i)|_{C_\tau}\bigr)
=
\frac{
w_i^\perp(\sigma,\tau)
-
w_i^\perp(\sigma',\tau)
}{
w(\sigma,\tau)
}
=a_i(\tau),
\]
where $w_i^\perp(\sigma,\tau) - w_i^\perp(\sigma',\tau) = a_i(\tau)w(\sigma,\tau)$ in $\tau^\perp \cap M$, and $a_i(\tau)$ is the coefficient of $u_i$ in the wall relation of Lemma~\ref{lem:wall}. Consequently,
\[
\deg N_{C_\tau/X_\Sigma}
=
\sum_{i\in I_\tau}a_i(\tau).
\]

\begin{figure}[htbp]
\centering
\begin{tikzpicture}[>=Stealth, scale=0.95]

    \begin{scope}[shift={(0,0)}]
        \fill[blue!10] (0,0) -- (2.3,0) -- (0,2.3) -- cycle;
        \fill[teal!10] (0,0) -- (0,2.3) -- (-1.9,-1.9) -- cycle;
        \fill[orange!10] (0,0) -- (-1.9,-1.9) -- (2.3,0) -- cycle;

        \foreach \x in {-2,...,2} {
            \foreach \y in {-2,...,2} {
                \fill[gray!55] (\x,\y) circle (1.1pt);
            }
        }

        \draw[->, very thick, blue!75!black] (0,0) -- (2.3,0)
            node[right, font=\small] {$\rho_1 = \tau_{23}$};
        \draw[->, very thick, teal!75!black] (0,0) -- (0,2.3)
            node[above, font=\small] {$\rho_2 = \tau_{13}$};
        \draw[->, very thick, red!80!black] (0,0) -- (-1.9,-1.9)
            node[below left, font=\small] {$\rho_3 = \tau_{12}$};

        \fill[black] (1,0) circle (2pt) node[below, font=\scriptsize] {$u_1$};
        \fill[black] (0,1) circle (2pt) node[left, font=\scriptsize] {$u_2$};
        \fill[black] (-1,-1) circle (2pt) node[above left, font=\scriptsize] {$u_3$};
        \fill[black] (0,0) circle (1.8pt);

        \node[font=\small, blue!75!black] at (1.1,1.1) {$\sigma_3$};
        \node[font=\small, teal!75!black] at (-1.0,0.8) {$\sigma_1$};
        \node[font=\small, orange!80!black] at (0.8,-1.0) {$\sigma_2$};
    \end{scope}

    \begin{scope}[shift={(7.2,-0.8)}]
        \fill[gray!10] (0,0) -- (3.2,0) -- (0,3.2) -- cycle;

        \draw[thick, gray!80!black] (0,0) -- (3.2,0)
            node[midway, below, yshift=-14pt, font=\footnotesize] {$C_{13}=\overline{O(\rho_2)}$};
        \draw[thick, gray!80!black] (0,0) -- (0,3.2)
            node[midway, left, xshift=-14pt, font=\footnotesize] {$C_{23}=\overline{O(\rho_1)}$};
        \draw[thick, gray!80!black] (3.2,0) -- (0,3.2)
            node[midway, above right, xshift=4pt, yshift=2pt, font=\footnotesize] {$C_{12}=\overline{O(\rho_3)}$};

        \draw[->, very thick, blue!75!black] (0,0) -- (1.15,0)
            node[below, font=\scriptsize] {$\chi_1$};
        \draw[->, very thick, blue!75!black] (0,0) -- (0,1.15)
            node[left, font=\scriptsize] {$\chi_2$};

        \draw[->, very thick, teal!75!black] (3.2,0) -- (2.05,0)
            node[below, font=\scriptsize] {$-\chi_1$};
        \draw[->, very thick, teal!75!black] (3.2,0) -- (2.35,0.85)
            node[above right, xshift=-2pt, font=\scriptsize] {$\chi_2-\chi_1$};

        \draw[->, very thick, orange!85!black] (0,3.2) -- (0,2.05)
            node[left, font=\scriptsize] {$-\chi_2$};
        \draw[->, very thick, orange!85!black] (0,3.2) -- (0.85,2.35)
            node[above right, yshift=-2pt, font=\scriptsize] {$\chi_1-\chi_2$};

        \fill[blue!75!black] (0,0) circle (2.6pt) node[below left, font=\small] {$p_{\sigma_3}$};
        \fill[teal!75!black] (3.2,0) circle (2.6pt) node[below right, font=\small] {$p_{\sigma_1}$};
        \fill[orange!85!black] (0,3.2) circle (2.6pt) node[above left, font=\small] {$p_{\sigma_2}$};
    \end{scope}

\end{tikzpicture}
\caption{Left: the fan of $\mathbb{P}^2$ with maximal cones $\sigma_1,\sigma_2,\sigma_3$ and walls $\tau_{ab}=\sigma_a\cap\sigma_b$. Right: the corresponding GKM graph of $\mathbb{P}^2$, showing the fixed points $p_{\sigma_j}$, the $\mathbb{T}_N$-invariant lines $C_{ab}$, and the tangent weights $w(\sigma_j,\tau_{ab})$ along the two incident edges at each vertex.}
\label{fig:gkm-weights-p2}
\end{figure}
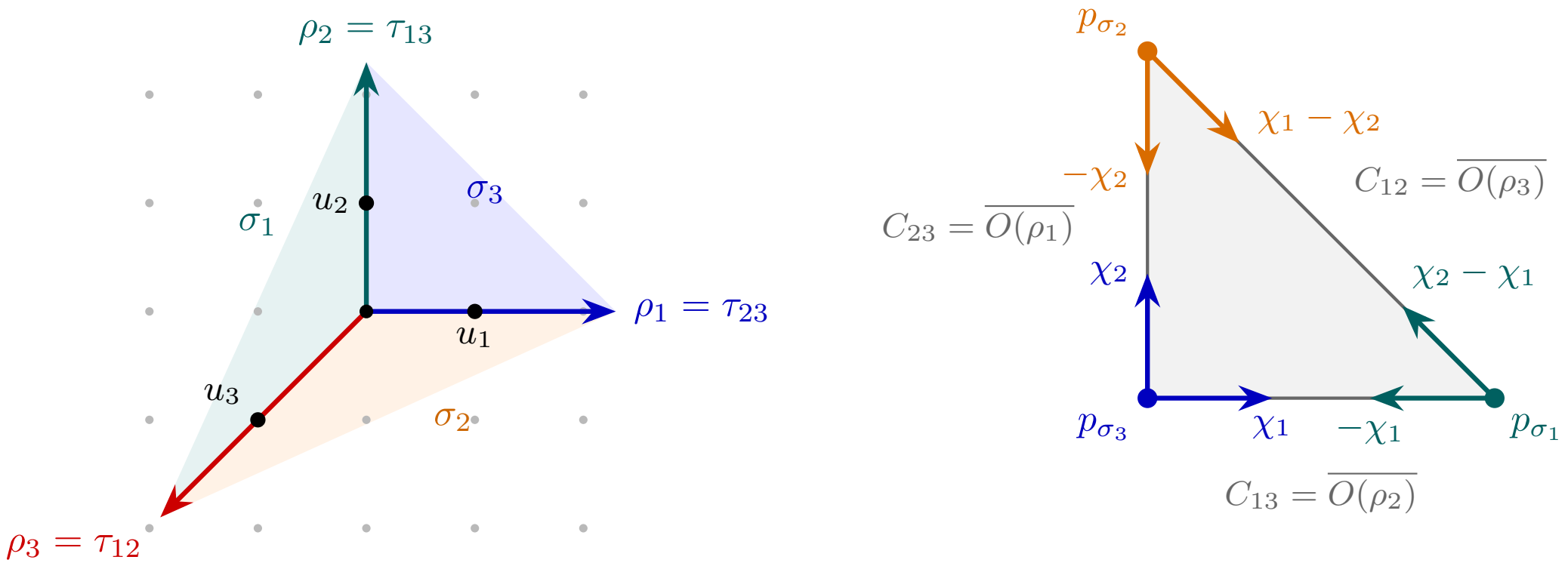

\begin{example}[Tangent and normal weights of $\mathbb{P}^2$]
Let $N=\mathbb{Z}^2$ with standard basis $e_1,e_2$, and consider the fan of $\mathbb{P}^2$ with primitive ray generators
\[
u_1=e_1=(1,0),\qquad
u_2=e_2=(0,1),\qquad
u_3=-e_1-e_2=(-1,-1).
\]
The maximal cones are
\[
\sigma_j=\operatorname{Cone}(u_i\mid i\neq j),
\qquad j=1,2,3.
\]
Write $C_{ab}$ for the $\mathbb{T}_N$-invariant curve joining $p_{\sigma_a}$ and $p_{\sigma_b}$ for $1 \le a < b \le 3$, and let $\chi_1,\chi_2\in M$ denote the dual basis to $e_1,e_2$.

The corresponding affine semigroups $S_{\sigma_j} = \sigma_j^\vee \cap M$ are generated over $\mathbb{Z}_{\geq 0}$ by
\[
S_{\sigma_1}
=
\mathbb{Z}_{\geq 0}(\chi_2-\chi_1) + \mathbb{Z}_{\geq 0}(-\chi_1),
\qquad
S_{\sigma_2}
=
\mathbb{Z}_{\geq 0}(\chi_1-\chi_2) + \mathbb{Z}_{\geq 0}(-\chi_2),
\qquad
S_{\sigma_3}
=
\mathbb{Z}_{\geq 0}\chi_1 + \mathbb{Z}_{\geq 0}\chi_2.
\]
For $a<b$, let $\tau_{ab}=\sigma_a\cap\sigma_b$. Then
\[
C_{12}=\overline{O(\rho_3)},\qquad
C_{13}=\overline{O(\rho_2)},\qquad
C_{23}=\overline{O(\rho_1)}.
\]
The tangent weights $w(\sigma_a,\tau_{ab})$ at $p_{\sigma_a}$, the normal fiber weights at the two endpoints $p_{\sigma_a}$ and $p_{\sigma_b}$, and the degrees of the normal bundles are summarized in the table below (and illustrated in Figure~\ref{fig:gkm-weights-p2}):
\[
\begin{array}{c|c|c|c}
& C_{12} & C_{13} & C_{23}\\
\hline
w(\sigma_a,\tau_{ab})
& \chi_2-\chi_1 & -\chi_1 & -\chi_2\\
w^\perp(\sigma_a,\tau_{ab})
& -\chi_1 & \chi_2-\chi_1 & \chi_1-\chi_2\\
w^\perp(\sigma_b,\tau_{ab})
& -\chi_2 & \chi_2 & \chi_1\\
\deg N_{C_{ab}/\mathbb{P}^2}
& 1 & 1 & 1
\end{array}
\]
\end{example}

\subsection{Equivariant cohomology and equivariant lifts}

Throughout the remainder of this section, we abbreviate the acting torus $\mathbb{T}_N$ to $T$. A key ingredient in the localization formula for genus-zero Gromov--Witten invariants is a choice of $\mathbb{Q}$-linear equivariant lift (i.e. a section)
\begin{center}
\begin{tikzcd}[cramped, column sep=large, baseline=-0.5ex]
H_{T}^*(X_\Sigma, \mathbb{Q}) \arrow[r, two heads] & H^*(X_\Sigma, \mathbb{Q}). \arrow[l, bend left=15, dashed, "s"]
\end{tikzcd}
\end{center}
In general, there is no canonical section $s$, but the resulting non-equivariant Gromov--Witten invariants are independent of this choice. We begin by recalling several structural properties of equivariant cohomology that underlie the construction of our equivariant lifts.

Let $X_{\Sigma}^T$ denote the fixed-point locus of the $T$-action on $X_{\Sigma}$. The inclusion $X_{\Sigma}^T \hookrightarrow X_{\Sigma}$ induces a restriction ring homomorphism $H_T^{\ast}(X_{\Sigma}, \mathbb{Z}) \rightarrow H_T^{\ast}(X_{\Sigma}^T, \mathbb{Z})$. By Proposition~\ref{prop:toricGKMgraph}, $X_\Sigma^T = \{p_\sigma \mid \sigma \in \Sigma(n)\}$ is finite. Denote the equivariant cohomology ring of a point by
\[
\Lambda_T = H_T^{\ast}(\mathrm{pt},\mathbb{Z})
\cong \operatorname{Sym}_{\mathbb{Z}}(M),
\qquad
(\Lambda_T)_{\mathbb{Q}} = \Lambda_T \otimes_{\mathbb{Z}} \mathbb{Q} \cong \operatorname{Sym}_{\mathbb{Q}}(M_{\mathbb{Q}}).
\]
Since the fixed-point locus is finite, we have
\[
H_T^{\ast}(X^T_{\Sigma}, \mathbb{Q}) = \bigoplus_{\sigma \in \Sigma(n)} H_T^{\ast}(p_{\sigma}, \mathbb{Q}) = \bigoplus_{\sigma \in \Sigma(n)} (\Lambda_T)_{\mathbb{Q}},
\]
giving the localization restriction homomorphism
\[
\Delta \colon H_T^{\ast}(X_{\Sigma}, \mathbb{Q}) \longrightarrow \bigoplus_{\sigma \in \Sigma(n)} (\Lambda_T)_{\mathbb{Q}}.
\]

The following results summarize the relationship between the equivariant and ordinary cohomology of $X_\Sigma$.

\begin{proposition}[\cite{cox2011toric}*{Proposition 12.4.7}]
    Let $X_{\Sigma}$ be a smooth projective toric variety. Then there is an isomorphism of $(\Lambda_T)_{\mathbb{Q}}$-modules
    \[
    H_T^{\ast}(X_{\Sigma}, \mathbb{Q}) \simeq (\Lambda_T)_{\mathbb{Q}} \otimes_{\mathbb{Q}} H^{\ast}(X_\Sigma, \mathbb{Q}).
    \]
\end{proposition}

Consequently, the ordinary cohomology ring $H^{\ast}(X_\Sigma,\mathbb{Q})$ is recovered from the $(\Lambda_T)_{\mathbb{Q}}$-algebra structure of $H_T^{\ast}(X_{\Sigma}, \mathbb{Q})$:

\begin{corollary}[\cite{cox2011toric}*{Corollary 12.4.8}]
    There is a canonical ring isomorphism
    \[
    H^{\ast}(X_{\Sigma}, \mathbb{Q}) \simeq H_T^{\ast}(X_{\Sigma}, \mathbb{Q})/I_T H_T^{\ast}(X_{\Sigma}, \mathbb{Q}) \simeq H_T^{\ast}(X_{\Sigma}, \mathbb{Q}) \otimes_{(\Lambda_T)_{\mathbb{Q}}} \mathbb{Q},
    \]
    where $I_T \subseteq (\Lambda_T)_{\mathbb{Q}}$ denotes the augmentation ideal generated by elements of positive degree.
\end{corollary}

Our goal is to describe a combinatorial procedure for computing the restrictions of equivariant lifts to the torus-fixed points using the Cartier data of $T$-invariant divisors. Let $D=\sum_{i=1}^r a_i D_i$ be a $T$-invariant Cartier divisor on $X_\Sigma$. For each cone $\sigma \in \Sigma$, there exists an element $m_\sigma\in M$ satisfying
\[
\langle m_\sigma,u_i\rangle = -a_i
\qquad \text{for all } i\in I_\sigma,
\]
and $m_{\sigma}$ is unique modulo $M(\sigma) = \sigma^{\perp} \cap M$ \cite{cox2011toric}*{Theorem 4.2.8}. Since for each $m \in M$ the character $\chi^m$ is a rational function on $X_{\Sigma}$ with principal divisor
\[
\operatorname{div}(\chi^m) = \sum_{i=1}^r \langle m, u_i \rangle D_i,
\]
we have $D|_{U_{\sigma}} = \operatorname{div}(\chi^{-m_{\sigma}})|_{U_{\sigma}}$. Thus $\{(U_{\sigma}, \chi^{-m_{\sigma}})\}_{\sigma \in \Sigma}$ is a system of local defining data for $D$, and we call $\{m_{\sigma}\}_{\sigma \in \Sigma}$ the \emph{Cartier data} of $D$.

To see how the Cartier data define an equivariant cohomology class, recall the Chang--Skjelbred exact sequence for the smooth projective toric variety $X_\Sigma$:
\begin{equation}\label{CS}
0 \longrightarrow H_T^{\ast}(X_\Sigma, \mathbb{Q}) \xrightarrow{\ \Delta\ } \bigoplus_{\sigma \in \Sigma(n)} H_T^{\ast}(O(\sigma), \mathbb{Q}) \longrightarrow \bigoplus_{\tau \in \Sigma(n-1)} H_T^{\ast}(O(\tau), \mathbb{Q}).
\end{equation}
For each $\sigma\in\Sigma(n)$, $O(\sigma)=\{p_\sigma\}$ is a fixed point, while for each wall $\tau\in\Sigma(n-1)$, the stabilizer of a point in $O(\tau)$ is the $(n-1)$-dimensional subtorus $T_\tau = (\operatorname{Span}(\tau)\cap N)\otimes_{\mathbb{Z}}\C^\times$ with character lattice $M_\tau = M/(\tau^{\perp} \cap M)$. Hence
\[
H_T^{\ast}(O(\sigma),\mathbb{Q})
\cong \operatorname{Sym}_{\mathbb{Q}}(M_{\mathbb{Q}})
\qquad \text{and} \qquad
H_T^{\ast}(O(\tau),\mathbb{Q})
\cong \operatorname{Sym}_{\mathbb{Q}}((M_\tau)_{\mathbb{Q}}),
\]
so sequence~\eqref{CS} takes the algebraic form
\[
0 \longrightarrow H_T^{\ast}(X_\Sigma, \mathbb{Q}) \longrightarrow \bigoplus_{\sigma \in \Sigma(n)} \operatorname{Sym}_{\mathbb{Q}}(M_{\mathbb{Q}}) \longrightarrow \bigoplus_{\tau \in \Sigma(n-1)} \operatorname{Sym}_{\mathbb{Q}}((M_\tau)_{\mathbb{Q}}).
\]
Viewing $\operatorname{Sym}_{\mathbb{Q}}(M_{\mathbb{Q}})$ and $\operatorname{Sym}_{\mathbb{Q}}((M_\tau)_{\mathbb{Q}})$ as the rings of polynomial functions with rational coefficients on $\sigma$ and on $\operatorname{Span}(\tau)$, respectively, exactness implies that an equivariant cohomology class in $H_T^{\ast}(X_\Sigma,\mathbb{Q})$ is uniquely determined by a tuple of polynomials $\{f_\sigma\}_{\sigma\in\Sigma(n)}$ satisfying the compatibility condition
\[
f_\sigma|_\tau = f_{\sigma'}|_\tau
\qquad \text{whenever } \sigma\cap\sigma'=\tau\in\Sigma(n-1).
\]
Because the Cartier data $\{m_\sigma\}_{\sigma \in \Sigma(n)}$ of $D$ satisfy $m_\sigma - m_{\sigma'} \in \tau^\perp \cap M$ whenever $\tau = \sigma \cap \sigma'$, the linear functions $\{-m_\sigma\}_{\sigma \in \Sigma(n)}$ agree along all shared walls $\tau$ and therefore define an equivariant lift of $[D] \in H^2(X_\Sigma,\mathbb{Q})$.

Explicitly, for each maximal cone $\sigma \in \Sigma(n)$, let $\{\eta_{\sigma,i}\}_{i\in I_\sigma} \subseteq M$ be the $\mathbb{Z}$-basis dual to the primitive ray generators $\{u_i\}_{i\in I_\sigma}$ of $\sigma$ (Note that $\eta_{\sigma,i} = m_{\sigma,i}$ in the notation of Proposition~\ref{prop:toricTangentWeights}):
\[
\langle \eta_{\sigma,j}, u_i \rangle = \delta_{ij}
\qquad \text{for all } i,j \in I_\sigma.
\]
For a $T$-invariant divisor $D = \sum_{i=1}^r a_i D_i$ with Cartier data $\{m_\sigma\}_{\sigma\in\Sigma(n)}$ defined by $\langle m_\sigma, u_i \rangle = -a_i$ ($i \in I_\sigma$), we have
\[
m_\sigma = -\sum_{i\in I_\sigma} a_i \eta_{\sigma,i},
\]
and we define the restriction of the equivariant lift of $D$ at $p_\sigma$ by
\[
s_\sigma(D) := -m_\sigma = \sum_{i\in I_\sigma} a_i \eta_{\sigma,i}.
\]
In particular, for each generator $x_i = [D_i]$ ($1 \le i \le r$), its Cartier data $\{m_\sigma^i\}_{\sigma\in\Sigma(n)}$ yield
\[
s_\sigma(x_i) = -m_\sigma^i =
\begin{cases}
    \eta_{\sigma,i}, & \text{if } i \in I_\sigma,\\
    0, & \text{if } i \notin I_\sigma.
\end{cases}
\]
For a monomial $x_{i_1}\cdots x_{i_k}$ in the generators, we set
\begin{equation}\label{LiftExtended}
s_\sigma(x_{i_1}\cdots x_{i_k}) := s_\sigma(x_{i_1})\cdots s_\sigma(x_{i_k}).
\end{equation}

Note that Equation~\eqref{LiftExtended} does not descend to a ring homomorphism on $H^{\ast}(X_\Sigma,\mathbb{Q})$, because it does not annihilate the linear ideal $\mathcal{I}$. Indeed, for any $m\in M$, the linear relation $\sum_{i=1}^r \langle m,u_i\rangle x_i = 0$ holds in $H^2(X_\Sigma,\mathbb{Q})$, whereas applying $s_\sigma$ yields
\[
\sum_{i=1}^r \langle m,u_i\rangle s_\sigma(x_i)
=
\sum_{i\in I_\sigma} \langle m,u_i\rangle \eta_{\sigma,i}
= m,
\]
which is nonzero for $m \neq 0$. To obtain a well-defined $\mathbb{Q}$-linear section $s \colon H^{\ast}(X_\Sigma,\mathbb{Q}) \to H_T^{\ast}(X_\Sigma,\mathbb{Q})$, we fix a basis of each graded piece $H^{2k}(X_\Sigma,\mathbb{Q})$ consisting of monomials in $x_1,\ldots,x_r$, evaluate $s_\sigma$ on these chosen basis monomials via~\eqref{LiftExtended}, and extend $\mathbb{Q}$-linearly to all of $H^{\ast}(X_\Sigma,\mathbb{Q})$.

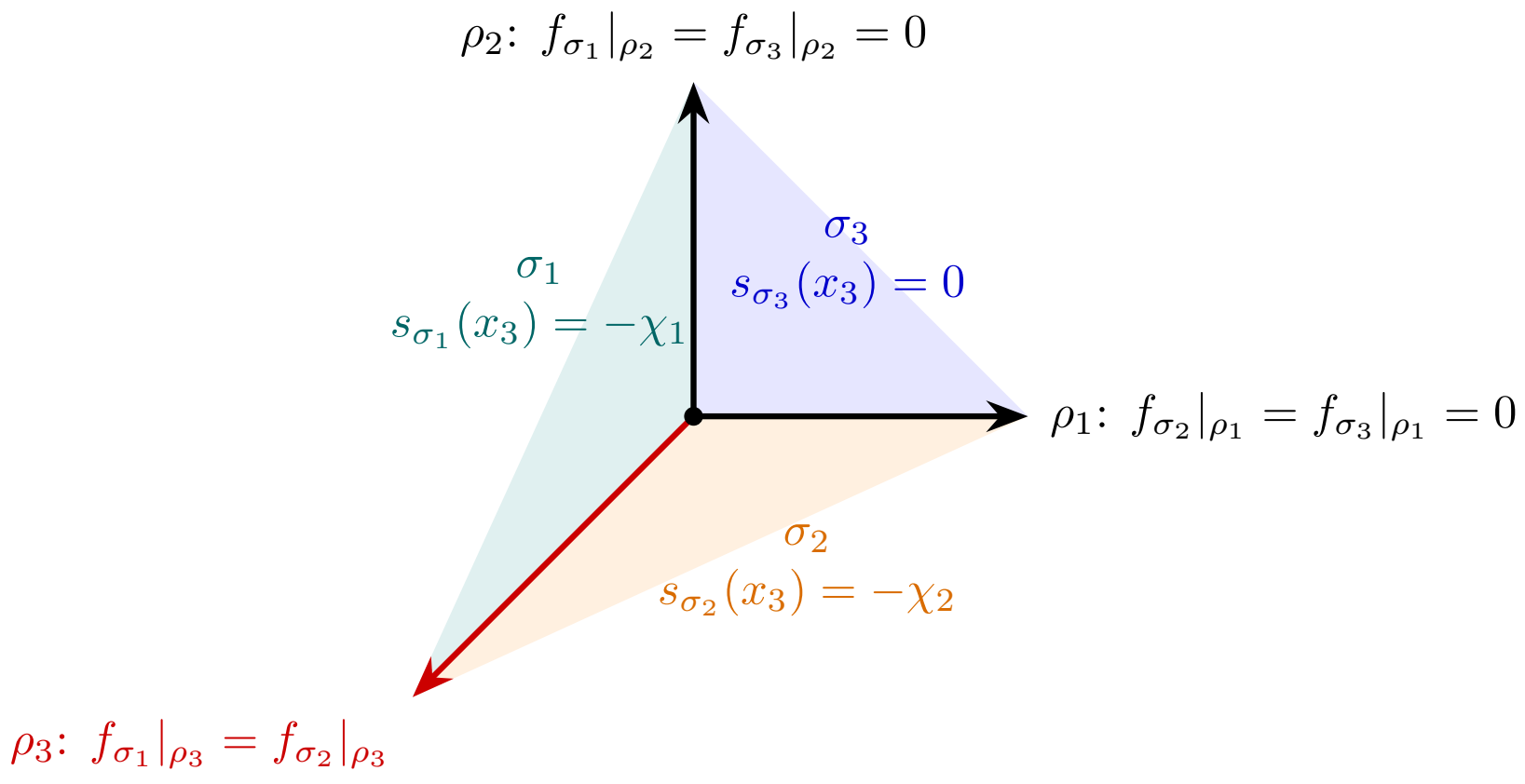
\begin{figure}[htbp]
\centering
\begin{tikzpicture}[>=Stealth, scale=0.95]

    \fill[blue!10] (0,0) -- (2.5,0) -- (0,2.5) -- cycle;
    \fill[teal!12] (0,0) -- (0,2.5) -- (-2.1,-2.1) -- cycle;
    \fill[orange!12] (0,0) -- (-2.1,-2.1) -- (2.5,0) -- cycle;

    \draw[->, very thick, black] (0,0) -- (2.5,0)
        node[right, font=\small] {$\rho_1$: $f_{\sigma_2}|_{\rho_1} = f_{\sigma_3}|_{\rho_1} = 0$};
    \draw[->, very thick, black] (0,0) -- (0,2.5)
        node[above, font=\small] {$\rho_2$: $f_{\sigma_1}|_{\rho_2} = f_{\sigma_3}|_{\rho_2} = 0$};
    \draw[->, very thick, red!80!black] (0,0) -- (-2.1,-2.1)
        node[below left, font=\small] {$\rho_3$: $f_{\sigma_1}|_{\rho_3} = f_{\sigma_2}|_{\rho_3}$};

    \fill[black] (0,0) circle (2pt);

    \node[align=center, font=\small, blue!80!black] at (1.15,1.15)
        {$\sigma_3$\\$s_{\sigma_3}(x_3) = 0$};
    \node[align=center, font=\small, teal!80!black] at (-1.15,0.85)
        {$\sigma_1$\\$s_{\sigma_1}(x_3) = -\chi_1$};
    \node[align=center, font=\small, orange!85!black] at (0.85,-1.15)
        {$\sigma_2$\\$s_{\sigma_2}(x_3) = -\chi_2$};

\end{tikzpicture}
\caption{The equivariant lift $s(x_3) \in H_T^2(\mathbb{P}^2,\mathbb{Q})$ viewed as a piecewise linear function $\{f_{\sigma_j} = s_{\sigma_j}(x_3)\}_{j=1}^3$ on the fan $\Sigma$: the linear forms $-\chi_1$ on $\sigma_1$, $-\chi_2$ on $\sigma_2$, and $0$ on $\sigma_3$ agree along the three walls $\rho_1, \rho_2, \rho_3$.}
\label{fig:piecewise-linear-lift}
\end{figure}

\begin{example}[Equivariant lifts for $\mathbb{P}^2$]
Let $e_1, e_2$ be the standard basis of $N = \mathbb{Z}^2$ with dual basis $\chi_1,\chi_2\in M$. As before, the fan of $\mathbb{P}^2$ has ray generators $u_1=e_1$, $u_2=e_2$, $u_3=-e_1-e_2$ and maximal cones $\sigma_j=\operatorname{Cone}(u_i\mid i\neq j)$ for $j=1,2,3$.

Choosing the monomial basis $\{1, x_3, x_3^2\}$ for $H^\ast(\mathbb{P}^2,\mathbb{Q}) \cong \mathbb{Q}[x_3]/\langle x_3^3\rangle$, it suffices to compute the Cartier data $\{m_{\sigma_j}\}_{j=1}^3$ of $D_3$:
\begin{itemize}
    \item On $\sigma_1=\operatorname{Cone}(u_2,u_3)$, solving $\langle m_{\sigma_1},u_2\rangle=0$ and $\langle m_{\sigma_1},u_3\rangle=-1$ gives $m_{\sigma_1}=\chi_1$, so $s_{\sigma_1}(x_3)=-\chi_1$.
    \item On $\sigma_2=\operatorname{Cone}(u_1,u_3)$, solving $\langle m_{\sigma_2},u_1\rangle=0$ and $\langle m_{\sigma_2},u_3\rangle=-1$ gives $m_{\sigma_2}=\chi_2$, so $s_{\sigma_2}(x_3)=-\chi_2$.
    \item On $\sigma_3=\operatorname{Cone}(u_1,u_2)$, solving $\langle m_{\sigma_3},u_1\rangle=0$ and $\langle m_{\sigma_3},u_2\rangle=0$ gives $m_{\sigma_3}=0$, so $s_{\sigma_3}(x_3)=0$.
\end{itemize}
Notice that along the wall $\rho_3 = \sigma_1 \cap \sigma_2 = \mathbb{R}_{\geq 0}(-1,-1)$, the linear functions $s_{\sigma_1}(x_3) = -\chi_1$ and $s_{\sigma_2}(x_3) = -\chi_2$ both restrict to $t \mapsto t$ (Figure~\ref{fig:piecewise-linear-lift}). Extending multiplicatively to $x_3^2$ via~\eqref{LiftExtended}, the restrictions of our equivariant lifts at the three fixed points are given by:
\[
\begin{array}{c|c|c|c}
& p_{\sigma_1} & p_{\sigma_2} & p_{\sigma_3}\\
\hline
1 & 1 & 1 & 1\\
x_3 & -\chi_1 & -\chi_2 & 0\\
x_3^2 & \chi_1^2 & \chi_2^2 & 0
\end{array}
\]
\end{example}

\subsection{Gromov--Witten invariants for smooth projective toric varieties}

We now assemble the geometric and combinatorial ingredients above to compute genus-zero Gromov--Witten invariants of smooth projective toric varieties via virtual localization. Following~\cites{alexeev2026atomic, Liu_Chiu,liu2013localization,holmes2025computationsequivariantgromovwittentheory}, they can be expressed as finite sums over decorated trees whose weights are determined entirely by the fan $\Sigma$.

\begin{figure}[htbp]
\centering
\begin{tikzpicture}[>=Stealth, scale=0.92]

    \begin{scope}[shift={(0,0)}]
        \draw[very thick, blue!75!black] (-2.2,0.5) .. controls (-1.1,1.0) and (-0.2,0.6) .. (0.4,0)
            node[midway, above, font=\scriptsize] {$C_{e_1} \xrightarrow{d_{e_1}:1} C_{\tau_{e_1}}$};
        \draw[very thick, teal!75!black] (0.4,0) .. controls (1.0,-0.6) and (1.9,-1.0) .. (3.0,-0.5)
            node[midway, below, font=\scriptsize] {$C_{e_2} \xrightarrow{d_{e_2}:1} C_{\tau_{e_2}}$};

        \draw[very thick, red!80!black] (0.0,-0.9) .. controls (0.3,-0.2) and (0.5,0.4) .. (0.8,1.1)
            node[above, font=\scriptsize] {$C_{v_2} \to \{p_{\ell(v_2)}\}$};

        \fill[blue!80!black] (-2.2,0.5) circle (2.3pt)
            node[below, font=\scriptsize] {$p_{\ell(v_1)}$};
        \fill[red!80!black] (0.4,0) circle (2.3pt);
        \fill[teal!80!black] (3.0,-0.5) circle (2.3pt)
            node[above, font=\scriptsize] {$p_{\ell(v_3)}$};

        \draw[thick] (-2.2,0.5) -- (-2.6,1.1) node[above, font=\scriptsize] {$\mathrm{mk}=1$};
        \draw[thick] (0.6,0.65) -- (1.2,0.8) node[right, font=\scriptsize] {$\mathrm{mk}=2$};
        \draw[thick] (0.15,-0.55) -- (-0.45,-0.75) node[left, font=\scriptsize] {$\mathrm{mk}=3$};
    \end{scope}

    \draw[<->, thick, gray!75!black] (3.9,0) -- (5.1,0)
        node[midway, above, font=\scriptsize] {corresponds to};

    \begin{scope}[shift={(8.2,0)}]
        \draw[very thick, blue!75!black] (-2.2,0.4) -- (0,0)
            node[midway, above, font=\footnotesize] {$e_1,\, d_{e_1}$};
        \draw[very thick, teal!75!black] (0,0) -- (2.2,-0.4)
            node[midway, above, font=\footnotesize] {$e_2,\, d_{e_2}$};

        \draw[thick, red!75!black] (-2.2,0.4) -- (-2.6,1.2)
            node[above, font=\scriptsize] {$1$};
        \draw[thick, red!75!black] (0,0) -- (-0.3,1.0)
            node[above, font=\scriptsize] {$2$};
        \draw[thick, red!75!black] (0,0) -- (0.3,1.0)
            node[above, font=\scriptsize] {$3$};

        \fill[black] (-2.2,0.4) circle (2.6pt)
            node[below, align=center, font=\scriptsize] {$v_1,\, \ell(v_1)$\\($n_{v_1}=2$)};
        \fill[black] (0,0) circle (2.6pt)
            node[below, align=center, font=\scriptsize] {$v_2,\, \ell(v_2)$\\($n_{v_2}=4$)};
        \fill[black] (2.2,-0.4) circle (2.6pt)
            node[below, align=center, font=\scriptsize] {$v_3,\, \ell(v_3)$\\($n_{v_3}=1$)};
    \end{scope}

\end{tikzpicture}
\caption{Left: a $T$-fixed genus-zero stable map $f \colon C \to X_\Sigma$ with two non-contracted components $C_{e_1}, C_{e_2}$ of degrees $d_{e_1}, d_{e_2}$, one contracted component $C_{v_2}$ mapped to $p_{\ell(v_2)}$, and $m=3$ marked points. Right: the associated decorated tree $\Gamma = (G,\ell,d,\mathrm{mk})$.}
\label{fig:decorated-tree-localization}
\end{figure}
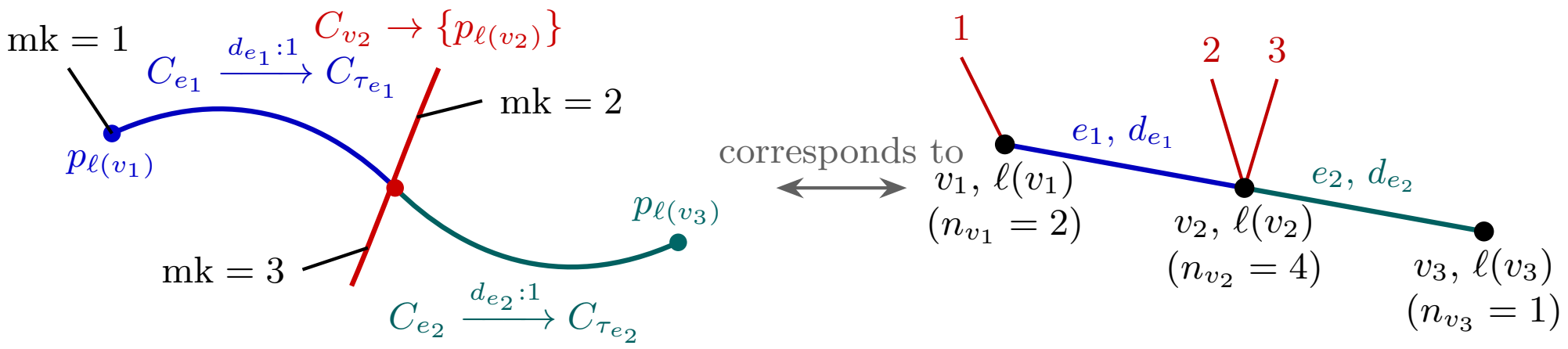

Let $X_\Sigma$ be a smooth projective toric variety of dimension $n$, and let $\beta\in H_2(X_\Sigma,\mathbb{Z})$ be a nonzero effective curve class. The virtual dimension of the moduli stack $\overline{\mathcal{M}}_{0,m}(X_\Sigma,\beta)$ of genus-zero stable maps to $X_\Sigma$ of class $\beta$ with $m$ marked points is
\[
\boldsymbol{\delta}(m,\beta)
=
n-3+m+\int_\beta c_1(TX_\Sigma).
\]
Consequently, for homogeneous cohomology classes $\gamma_j\in H^{2k_j}(X_\Sigma,\mathbb{Q})$ ($1 \le j \le m$), the genus-zero Gromov--Witten invariant
\[
\left\langle
\gamma_1,\ldots,\gamma_m
\right\rangle^{X_\Sigma}_{0,\beta}
\]
can be nonzero only if the dimension constraint below is verified
\begin{equation}\label{eq:GWdimensiontoric}
\sum_{j=1}^m k_j
=
n-3+m+\int_\beta c_1(TX_\Sigma)
\end{equation}
Next, we explain how one can use virtual localization to reduce the evaluation of this integral to the connected components of the $T$-fixed locus $\overline{\mathcal{M}}_{0,m}(X_\Sigma,\beta)^T$; see also Figure~\ref{fig:decorated-tree-localization}.

\begin{definition}\label{def:decoratedtree}
Let $m\geq 0$ and let $\beta\in H_2(X_\Sigma,\mathbb{Z})$ be an effective nonzero curve class. A \emph{decorated tree of type $(m,\beta)$} is a tuple $\Gamma=(G,\ell,d,\mathrm{mk})$ consisting of:
\begin{enumerate}
\item[(a)] a finite connected tree $G$ with vertex set $V$ and edge set $E$;
\item[(b)] a vertex labeling $\ell \colon V \to \Sigma(n)$ such that, for every edge $e=\{v,v'\} \in E$,
\[
\tau_e:=\ell(v)\cap\ell(v')
\in\Sigma(n-1);
\]
\item[(c)] an edge degree map $d \colon E \to \mathbb{Z}_{>0}$ assigning a positive integer $d_e$ to each edge $e \in E$;
\item[(d)] a marking map $\mathrm{mk} \colon \{1,\ldots,m\} \to V$ distributing the $m$ marked points among the vertices;
\end{enumerate}
subject to the total homology degree constraint
\[
\sum_{e\in E} d_e[C_{\tau_e}] = \beta
\qquad
\text{in } H_2(X_\Sigma,\mathbb{Z}).
\]
An isomorphism of decorated trees is a graph isomorphism between the underlying trees that preserves the maps $\ell$, $d$, and $\mathrm{mk}$. We denote by $\operatorname{Aut}\Gamma$ the automorphism group of $\Gamma$. For each vertex $v\in V$, let $\operatorname{val}(v)$ denote its valence (number of incident edges) and set
\[
n_v := \operatorname{val}(v) + |\mathrm{mk}^{-1}(v)|
\]
for the total number of special points on the component corresponding to $v$.
\end{definition}

Geometrically, a vertex $v \in V$ with $n_v \geq 3$ corresponds to a maximal genus-zero contracted component of the domain curve mapped constantly to the fixed point $p_{\ell(v)}$. In contrast, a vertex with $n_v \leq 2$ corresponds to an unmarked or singly marked ramification point over $p_{\ell(v)}$. Each edge $e=\{v,v'\} \in E$ corresponds to a non-contracted rational component $C_e \cong \mathbb{P}^1$ mapped onto the invariant curve $C_{\tau_e}$ as a $d_e$-fold cover totally ramified over $p_{\ell(v)}$ and $p_{\ell(v')}$.

\begin{definition}\label{def:omega}
Let $\Gamma=(G,\ell,d,\mathrm{mk})$ be a decorated tree and $e=\{v,v'\} \in E$ an edge. The \emph{normalized tangent weight} of $e$ at $v$ is
\[
\omega_{v,e}
:=
\frac{w(\ell(v),\tau_e)}{d_e}
\in M_{\mathbb{Q}},
\]
where $w(\ell(v),\tau_e)$ is the tangent weight of the $T$-action on $C_{\tau_e}$ at the fixed point $p_{\ell(v)}$.
\end{definition}

We now define the vertex and edge contributions appearing in the localization formula. First, by Proposition~\ref{prop:toricTangentWeights}, the equivariant Euler class of the tangent space $T_{p_\sigma}X_\Sigma$ at a fixed point $p_\sigma$ ($\sigma\in\Sigma(n)$) is the product of the tangent weights of the incident invariant curves:
\[
e_\sigma
:=
e_T(T_{p_\sigma}X_\Sigma)
=
\prod_{\substack{\tau\in\Sigma(n-1)\\
\tau\preceq\sigma}}w(\sigma,\tau).
\]

Second, the edge contribution $\Omega_e$ accounts for the equivariant Euler class of the deformation and obstruction spaces of the $d_e$-fold cover of $C_{\tau_e}$, determined by the splitting of the normal bundle $N_{C_{\tau_e}/X_\Sigma} \simeq \bigoplus_{i\in I_{\tau_e}} \cO_{X_\Sigma}(D_i)|_{C_{\tau_e}}$.

\begin{definition}\label{def:OmegaTheta}
Let $\Gamma=(G,\ell,d,\mathrm{mk})$ be a decorated tree. For an edge $e=\{v,v'\}\in E$ read at one of its endpoints $v$, set $\sigma=\ell(v)$, $\tau=\tau_e=\ell(v)\cap\ell(v')$, $d=d_e$, and $\omega=\omega_{v,e}$. For each $i \in I_\tau = \{i \in \{1,\ldots,r\} \mid \rho_i \subseteq \tau\}$, let $m_i := d \cdot a_i(\tau)$ be the degree of the pullback of $\cO_{X_\Sigma}(D_i)|_{C_\tau}$ along the $d$-fold cover, where $a_i(\tau)$ is the wall coefficient from Lemma~\ref{lem:wall}, and let $\tau_i = \operatorname{Cone}(u_j \mid j \in I_\sigma \setminus \{i\})$ so that $w(\sigma,\tau_i) = w_i^\perp(\sigma,\tau)$. We define
\[
  \Omega_{e} :=
  \frac{1}{(-1)^{d}(d!)^{2}\,\omega^{2d}}
  \prod_{i \in I_\tau}
  \begin{cases}
    \displaystyle\prod_{j=1}^{-m_i-1}
      \bigl(w(\sigma, \tau_i)+j\,\omega\bigr),
      & \text{if } m_i\le -2,\\[1.2em]
    1, & \text{if } m_i = -1,\\[0.4em]
    \displaystyle\prod_{k=0}^{m_i}
      \bigl(w(\sigma, \tau_i)-k\,\omega\bigr)^{-1},
      & \text{if } m_i\ge 0.
  \end{cases}
\]
\end{definition}
\begin{remark}
    Although Definition~\ref{def:OmegaTheta} uses an endpoint $v$, its value is independent of this choice. In fact, let $\sigma'=\ell(v')$ and write
$\omega'=\omega_{v',e}=-\omega$. For $i\in I_\tau$, denote by $m_{\sigma,i}$ and $m_{\sigma',i}$ the dual-basis elements associated with the ray $\rho_i$ in $\sigma$ and $\sigma'$, respectively. Then
\[
m_{\sigma',i}=m_{\sigma,i}-a_i(\tau)w(\sigma,\tau)
=m_{\sigma,i}-m_i\omega.
\]
It follows that the normal weights at the opposite endpoint are  $w_i^\perp(\sigma',\tau)=m_{\sigma', i}$, which, upon direct substitution into the products in Definition~\ref{def:OmegaTheta}, gives the same $\Omega_e$ at either endpoint.
\end{remark}

\begin{definition}
Let $\Gamma=(G,\ell,d,\mathrm{mk})$ be a decorated tree. For a vertex $v \in V$ with $n_v \le 2$, we define
\[
    \Theta_{v} :=
    \begin{cases}
        \omega_{v,e}, & \text{if } n_{v}=1 \text{ (one incident edge } e \text{ and no marks)},\\[0.3em]
        1, & \text{if } n_{v}=2 \text{ with one incident edge and one mark},\\[0.3em]
        \bigl(\omega_{v,e_{1}}+\omega_{v,e_{2}}\bigr)^{-1},
        & \text{if } n_{v}=2 \text{ with two incident edges } e_1, e_2 \text{ and no marks}.
  \end{cases}
\]
\end{definition}

We can now state the genus-zero localization formula for smooth projective toric varieties.

\begin{theorem}\label{thm:localization}
Let $X_\Sigma$ be a smooth projective toric variety of dimension $n$. Let $\beta\in H_2(X_\Sigma,\mathbb{Z})$ be an effective nonzero curve class, and let $\gamma_1,\ldots,\gamma_m\in H^\ast(X_\Sigma,\mathbb{Q})$, with $\gamma_j\in H^{2k_j}(X_\Sigma,\mathbb{Q})$, satisfy $\sum_{j=1}^m k_j = \boldsymbol{\delta}(m,\beta)$. Then
\[
\big\langle
\gamma_1,\ldots,\gamma_m
\big\rangle^{X_\Sigma}_{0,\beta}
=
\sum_\Gamma
\frac{1}{|\operatorname{Aut}\Gamma|
\prod_{e\in E}d_e}
\ \prod_{v\in V}
e_{\ell(v)}^{\operatorname{val}(v)-1}
\ \prod_{e\in E}\Omega_e
\ \prod_{v\in V}I_v
\ \prod_{j=1}^m
s_{\ell(\mathrm{mk}(j))}(\gamma_j),
\]
where the sum runs over the finitely many isomorphism classes of decorated trees $\Gamma=(G,\ell,d,\mathrm{mk})$ of type $(m,\beta)$, with $I_v = \Theta_v$ for $n_v \leq 2$ and
\[
I_v =
\left(\prod_{e\ni v}\omega_{v,e}^{-1}\right)
\left(\sum_{e\ni v}\omega_{v,e}^{-1}\right)^{n_v-3}
\qquad \text{for } n_v \geq 3.
\]
\end{theorem}

\subsection{Gromov--Witten invariants for complete intersections in toric varieties}

Following~\cite{alexeev2026atomic}, we now extend the genus-zero localization formula to smooth complete intersections $Y \subset X_\Sigma$ by incorporating the equivariant Euler class of the obstruction bundle induced by the defining vector bundle $\cE \to X_\Sigma$. We begin by computing the fiber weights of equivariant line bundles on $X_\Sigma$ at the torus-fixed points.

A $T$-equivariant line bundle on $X_\Sigma$ is a line bundle $L \to X_\Sigma$ equipped with a linearization, that is, a lift of the $T$-action to the total space of $L$ that is linear on fibers. Over each fixed point $p_{\sigma}$, the torus $T$ acts on the one-dimensional fiber $L|_{p_{\sigma}}$ by a character in $M$, which we now determine.

Since every line bundle $L \in \Pic(X_\Sigma)$ is isomorphic to $\cO_{X_\Sigma}(D)$ for a $T$-invariant Cartier divisor $D = \sum_{i=1}^r a_i D_i$, $L$ inherits a natural $T$-linearization from the left action of $T$ on rational functions given by $(t \cdot f)(x) = f(t^{-1} \cdot x)$. If $\{(U_\sigma, \chi^{-m_\sigma})\}_{\sigma\in\Sigma(n)}$ is the local Cartier data for $D$, where $\langle m_\sigma, u_i \rangle = -a_i$ for $i \in I_\sigma$, then over the affine chart $U_{\sigma}$ we have
\[
\begin{aligned}
\mathcal O_{X_\Sigma}(D)(U_\sigma)
&=\{0\}\cup
\left\{f\in\mathbb C(X_\Sigma)^\times\ \middle|\
\left.(\operatorname{div}(f)+D)\right|_{U_\sigma}\geq0\right\}\\
&=\mathbb C[U_\sigma]\,\chi^{m_\sigma}.
\end{aligned}
\]
Thus, $\chi^{m_\sigma}$ is a local frame for $\cO_{X_\Sigma}(D)$ over $U_\sigma$. Under the action $(t \cdot f)(x) = f(t^{-1} \cdot x)$, an element $t \in T$ acts on this generator by
\[
(t \cdot \chi^{m_\sigma})(x) = \chi^{m_\sigma}(t^{-1} \cdot x) = \chi^{-m_\sigma}(t)\chi^{m_\sigma}(x),
\]
so the fiber of $\cO_{X_\Sigma}(D)$ at $p_\sigma$ is a $T$-weight space of weight $-m_\sigma = s_\sigma(D)$.

\begin{proposition}\label{weightsbundles}
    With respect to the natural $T$-linearization on $L \simeq \cO_{X_\Sigma}(D)$ induced by the $T$-invariant divisor $D = \sum_{i=1}^r a_i D_i$ with Cartier data $\{m_\sigma\}_{\sigma\in\Sigma(n)}$, the weight of the fiber $L|_{p_\sigma}$ is
    \[
    \eta(\sigma) = -m_\sigma = s_\sigma(D) = \sum_{i \in I_\sigma} a_i \eta_{\sigma,i}.
    \]
\end{proposition}

\begin{remark}
    Every line bundle on a normal toric variety admits a $T$-linearization, and two choices of linearization on the same underlying line bundle differ by twisting by a global character $m \in M$, which shifts all fiber weights $\eta(\sigma)$ uniformly by $m$. This freedom corresponds precisely to changing the $T$-invariant divisor representative $D$ within its linear equivalence class by $\operatorname{div}(\chi^{-m})$.
\end{remark}

Now let $X_\Sigma$ be a smooth projective toric variety of dimension $n$, and consider a split vector bundle of rank $s$:
\[
\cE=\bigoplus_{j=1}^{s}L_j,
\qquad
L_j=\cO_{X_\Sigma}\left(\sum_{i=1}^{r}d_{j,i}D_i \right),
\quad d_{j,i}\in\mathbb{Z}.
\]
We assume that $\cE$ is globally generated (equivalently, each $L_j$ is nef by Proposition~\ref{Toric KleimanNef}), and let $Y=\mathscr{Z}(X_\Sigma,\cE) \subset X_\Sigma$ denote the zero locus of a general global section of $\cE$, assumed nonempty and smooth of codimension $s=\rank\cE$. Equipping each summand $L_j$ with its natural $T$-linearization makes $\cE$ a $T$-equivariant vector bundle on $X_\Sigma$.

Let $\overline{\mathcal{M}}_{0,m}(X_\Sigma,\beta)$ be the moduli stack of genus-zero stable maps to $X_\Sigma$ of class $\beta$, and consider the universal curve and universal evaluation map
\[
\pi \colon \mathcal{C}_{0,m,\beta}
\longrightarrow
\overline{\mathcal{M}}_{0,m}(X_\Sigma,\beta),
\qquad
f \colon \mathcal{C}_{0,m,\beta}\longrightarrow X_\Sigma.
\]
Since $\cE$ is globally generated, it is convex, meaning that $H^1(C,f^*\cE)=0$ for every genus-zero stable map $[f \colon C \to X_\Sigma]$. Consequently, the direct image
\[
\cE_{0,m,\beta} := \pi_*f^*\cE
\]
is a $T$-equivariant vector bundle on $\overline{\mathcal{M}}_{0,m}(X_\Sigma,\beta)$ of rank
\[
\rank\cE_{0,m,\beta}
=
s+\int_\beta c_1(\cE),
\]
whose fiber over $[f \colon C \to X_\Sigma]$ is $H^0(C,f^*\cE)$.

Let $\iota_\beta$ denote the proper morphism from the disjoint union of the moduli stacks of stable maps to $Y$ with ambient pushforward class $\beta$ to $\overline{\mathcal M}_{0,m}(X_\Sigma,\beta)$. Since $\mathcal E$ is convex and $Y$ is the regular zero locus of a section, we get
\[
(\iota_\beta)_*
\sum_{\iota_*\beta'=\beta}
[\overline{\mathcal M}_{0,m}(Y,\beta')]^{\mathrm{vir}}
=
c_{\mathrm{top}}(\mathcal E_{0,m,\beta})\cap
[\overline{\mathcal M}_{0,m}(X_\Sigma,\beta)]^{\mathrm{vir}}.
\]


For each edge $e=\{v,v'\}\in E$ with wall $\tau_e = \ell(v)\cap\ell(v')\in\Sigma(n-1)$ and each summand $L_j$, the restriction $L_j|_{C_{\tau_e}}$ has non-negative degree
\[
c_j(\tau_e) := \deg\bigl(L_j|_{C_{\tau_e}}\bigr) = \sum_{i=1}^r d_{j,i} b_i(\tau_e) \geq 0,
\]
so its pullback to the edge component $C_e \cong \mathbb{P}^1$ has degree $b_j := d_e c_j(\tau_e) \geq 0$. The space of sections $H^0(C_e, f_e^* L_j)$ therefore decomposes into $b_j + 1$ weight spaces with weights $\eta_j(\ell(v)) - k\omega_{v,e}$ for $k = 0, \ldots, b_j$, where $\eta_j(\ell(v)) = s_{\ell(v)}(L_j)$ is the fiber weight of $L_j$ at $p_{\ell(v)}$ from Proposition~\ref{weightsbundles}. With independent auxiliary scalar weights
$\lambda_1,\ldots,\lambda_s$, define the edge twisting factor
by
\[
\mathrm{Tw}_e(\boldsymbol\lambda)
:= \prod_{j=1}^{s} \prod_{k=0}^{d_ec_j(\tau_e)} \bigl(\lambda_j+\eta_j(\ell(v))-k\,\omega_{v,e}\bigr).
\]
\begin{remark}
The edge factor is independent of the chosen end point, as the weights at the other end point satisfy \(\eta_j(\ell(v'))=\eta_j(\ell(v)) -d_ec_j(\tau_e)\omega_{v,e}\) 
\end{remark}

\begin{proposition}\label{prop:toricTwist}
Let $\Gamma=(G,\ell,d,\mathrm{mk})$ be a decorated tree of type $(m,\beta)$. In the fraction field of the equivariant coefficient ring enlarged by the auxiliary parameters $\lambda_1,\ldots,\lambda_s$, the restriction of the auxiliary-equivariant Euler class of $\mathcal E_{0,m,\beta}$ to the fixed locus $F_\Gamma$ equals
\[
\mathrm{Tw}(\Gamma;\boldsymbol\lambda)
=
\prod_{e\in E}\mathrm{Tw}_e(\boldsymbol\lambda)
\prod_{v\in V}
\prod_{j=1}^{s}
\bigl(\lambda_j+\eta_j(\ell(v))\bigr)^{1-\operatorname{val}(v)},
\]
where
\[
\mathrm{Tw}_e(\boldsymbol\lambda)
=
\prod_{j=1}^{s}
\prod_{k=0}^{d_ec_j(\tau_e)}
\bigl(\lambda_j+\eta_j(\ell(v))-k\omega_{v,e}\bigr).
\]
\end{proposition}

We can now state the genus-zero localization formula for smooth complete intersections in smooth projective toric varieties.

\begin{theorem}\label{thm:localizationCI}
Let $X_\Sigma$ be a smooth projective toric variety of dimension $n$, and let $\mathcal E=\bigoplus_{j=1}^{s}L_j$ be a globally generated split vector bundle. Suppose that a regular section has smooth zero locus $Y$ of codimension $s$, with inclusion $\iota: Y\hookrightarrow X_\Sigma$. Let $\beta\in H_2(X_\Sigma,\mathbb Z)$ be a nonzero effective class, and let $\gamma_k=\iota^*\alpha_k$, where $\alpha_k\in H^{2q_k}(X_\Sigma,\mathbb Q)$. Assume
\[
\sum_{k=1}^{m}q_k
=n-s+m-3+\int_\beta(c_1(TX_\Sigma)-c_1(\mathcal E)).
\]
Then
\[
\begin{aligned}
&\sum_{\substack{\beta'\in H_2(Y,\mathbb Z)\\
                   \iota_*\beta'=\beta}}
\left\langle\gamma_1,\ldots,\gamma_m\right\rangle^Y_{0,\beta'}
\\
&\quad=
\left.\left[
\sum_{\Gamma}
\frac{1}{|\operatorname{Aut}\Gamma|\prod_{e\in E}d_e}
\prod_{v\in V}e_{\ell(v)}^{\operatorname{val}(v)-1}
\prod_{e\in E}\Omega_e
\prod_{v\in V}I_v
\,\mathrm{Tw}(\Gamma;\boldsymbol\lambda)
\prod_{k=1}^{m}s_{\ell(\mathrm{mk}(k))}(\alpha_k)
\right]\right|_{\boldsymbol\lambda=0,\,H_T^{>0}(\mathrm{pt})=0},
\end{aligned}
\]
where the sum runs over the decorated trees of type $(m,\beta)$. 
\end{theorem}
\begin{remark}
The factors $\Omega_e$ and $I_v$ are rational functions of the torus weights so that Theorem~\ref{thm:localization} is evaluated in the fraction field of $H_T^*(\mathrm{pt},\mathbb Q)$, with weights kept generic until the decorated-tree sum is computed. The dimensional constraint in the statement implies that the sum is constant and equal to the ordinary Gromov--Witten invariant. 
\end{remark}

\section{Obstructions to rationality}\label{sec:atoms}

The theory of Hodge atoms, introduced by Katzarkov--Kontsevich--Pantev--Yu~\cite{KKPY}, explores the bridge between Gromov--Witten theory and birational geometry. This section is a brief recap of the theory of Hodge atoms and additional invariants that can be computed from the A-model F-bundle, explaining how to obstruct rationality using these tools. For a complete discussion on this topic, with all the proofs, the reader is referred to~\cites{KKPY, Leo_Luiz_Pedro_Lud_Max, CKK-lectures}; those who are familiar with Hodge atoms and F-bundles can skip to Section~\ref{sec:Examples} for the applications.

\subsection{A-model F-bundle}

Let $X$ be a smooth projective variety of dimension $n$. Fix a homogeneous $\rH^*(X,\Q)$-basis $\{\Delta_i\}$ with $\Delta_0=1$ such that $\Delta_1,\ldots,\Delta_\rho$ form an integral basis of $\mathsf{NS} (X)/\mathrm{torsion}$, where $\rho$ is the Picard number. Let $\delta_i$ be the coordinates dual to $\Delta_i$, and let $q_1,\ldots,q_\rho$ be the Novikov coordinates. We work over the non-Archimedean field
\[
\Bbbk=\widehat{\bigcup_{j\geq1}\overline{\Q}((T^{1/j}))},
\]
with its $T$-adic valuation, normalized by $|T|_\Bbbk=e^{-1}$. The restriction of this norm to $\overline\Q$ is trivial.

\begin{definition}[\cite{KKPY}]
The \emph{A-model base $B_X$} is the $\Bbbk$-analytic superspace
\[
B_X=U_{\mathrm{Nov}}\times U_{\mathrm{even}}\times U_{\mathrm{odd}},
\]
where
\begin{enumerate}
\item $U_{\mathrm{odd}}$ is the purely odd superspace whose structure sheaf is the exterior algebra on the coordinates $\delta_i$ of odd degree;
\item $U_{\mathrm{even}}$ is the product of the analytic affine line in the unit direction $\delta_0$ with the open unit polydisc in the remaining even coordinates $\delta_i$, $i\notin\{0,1,\ldots,\rho\}$;
\item $U_{\mathrm{Nov}}$ is the torus domain
\[
U_{\mathrm{Nov}}=\left\{q\in(\mathbb G_{m,\Bbbk}^{\rho})^{\mathrm{an}}:
-\sum_{m=1}^{\rho}\log|q_m|\,\Delta_m
\text{ lies in the ample cone of }\mathsf{NS}(X)_{\R}\right\}.
\]
\end{enumerate}
\end{definition}

The super tangent sheaf of $B_X$ is trivialized by $\partial_{\delta_i}\mapsto\Delta_i$ in the additive directions and $q_m\partial_{q_m}\mapsto\Delta_m$ in the divisor directions. Thus, the full tangent superspace is identified with $\rH^*(X,\Bbbk)$.

Denote by $\mathbb{D}$ the germ at $0$ of the $\Bbbk$-analytic disc with coordinate $u$.

\begin{definition}\label{def:A-model-F-bundle}
The \emph{A-model $F$-bundle} is the pair $(\mathcal H,\nabla)/B_X$, where
\[
\mathcal H=\rH^*(X,\Q)\otimes_{\Q}\mathcal O_{B_X\times\mathbb D}
\]
and the quantum connection is
\begin{align}
\nabla_{\partial_u}&=\partial_u-\frac{1}{u^2}(\mathsf{Eu}\star(-)) +\frac{1}{u}\mathsf{Gr},\\ 
\nabla_{q_m\partial_{q_m}}&=q_m\partial_{q_m}+\frac{1}{u}(\Delta_m\star(-)),\\
\nabla_{\partial_{\delta_i}}&=\partial_{\delta_i}+\frac{1}{u}(\Delta_i\star(-)).
\end{align}
Here $\star$ is the big quantum product and
\begin{equation}\label{eq:euler-grading}
\mathsf{Gr}|_{\rH^a(X)}=\frac{a-n}{2}\,\mathrm{id}, \qquad \mathsf{Eu}=c_1(TX)+\sum_{i\notin\{1,\ldots,\rho\}} \left(1-\frac{\deg\Delta_i}{2}\right)\delta_i\Delta_i.
\end{equation}
\end{definition}

\begin{remark}
For a smooth projective toric variety $X_\Sigma$ one has $\rH^{\mathrm{odd}}(X_\Sigma,\Q)=0$, and all its rational cohomology is of type $(p,p)$, as recalled in Section~\ref{sec:picard}. Hence $B_{X_\Sigma}=U_{\mathrm{Nov}}\times U_{\mathrm{even}}$ and $B_{X_\Sigma}^{\mathsf{Hod}}=B_{X_\Sigma}$.
\end{remark}

The \emph{small locus} is the subspace $U_{\mathrm{Nov}}\subset B_X$ obtained by setting all additive coordinates $\delta_i$ equal to zero. On this locus, the $u$-direction of the quantum connection becomes
\begin{equation}\label{eq:sign}
\nabla_{\partial_u}=\partial_u-\frac{1}{u^2}\bigl(c_1(TX)\star_{\mathrm{small}}(-)\bigr)
+\frac{1}{u}\mathsf{Gr}.
\end{equation}
For Fano varieties, each coefficient of the small quantum product involves only finitely many Gromov--Witten invariants. Thus the software can compute the matrix of $c_1(TX)\star_{\mathrm{small}}(-)$. Whether its eigenvalue factors represent the generic atomic splitting requires a further argument, such as Lemma~\ref{lem:cyclic}.

 
 

Write
\[
\rH^*_{\mathrm{Hod}}(X)=\bigoplus_p\bigl(\rH^{2p}(X,\Q)\cap\rH^{p,p}(X)\bigr).
\]
Over the small locus $B_X^{\mathrm{small}}=U_{\mathrm{Nov}}\subset B_X^{\mathsf{Hod}}$, the Hodge equivariance of the quantum product
implies that the quantum connection preserves
\[
\mathcal H^{\mathrm{Hod}}=\rH^*_{\mathrm{Hod}}(X)\otimes_{\Q}
\mathcal O_{B_X^{\mathrm{small}}\times\mathbb D}.
\]
We shall use this restricted connection in Section~\ref{subsec:rat-obstructions}.


\subsection{Hodge Atoms}
Let $\mathsf{MT}$ be the Mumford--Tate group of $\rH^*(X,\Q)$, and let $\mathsf{Hod}$ be the kernel of its character on the Tate line $\Q(1)$, as in \cite{moonen1999notes}. Its invariant subspace is
\[
\rH^*(X,\Q)^{\mathsf{Hod}}=\bigoplus_p\bigl(\rH^{2p}(X,\Q)\cap\rH^{p,p}(X)\bigr)=\rH^*_{\mathrm{Hod}}(X).
\]
We denote by $B_X^{\mathsf{Hod}}$ the Hodge-invariant analytic subspace of $B_X$ obtained by retaining the invariant additive directions and all Novikov coordinates.

Consider $\boldsymbol\kappa_b=\mathsf{Eu}_b\star_b(-)$ on $\mathcal H_{(b,0)}$. Its characteristic polynomial is
\[
\chi_b(\lambda)=\det\bigl(\lambda\,\mathrm{id}-\boldsymbol\kappa_b\bigr).
\]
Let $U_X^{\mathrm{gen}}\subset B_X^{\mathsf{Hod}}$ be the open locus where $\chi_b$ has the maximal number of distinct roots. Over this locus, let
\[
\pi\colon\widetilde U_X^{\mathrm{gen}}
\longrightarrow U_X^{\mathrm{gen}}
\]
be the \emph{reduced} spectral cover, whose fiber over $b$ consists of the distinct eigenvalues of $\boldsymbol\kappa_b$.
\begin{definition}
The \emph{local Hodge atoms} of $X$ are the connected components $\alpha\in\pi_0(\widetilde U_X^{\mathrm{gen}})$.
\end{definition}

At $b\in U_X^{\mathrm{gen}}(\Bbbk)$ the Euler operator has the generalized eigenspace decomposition
\begin{equation}\label{eq:split_F_bundle_fiber}
\mathcal H_{(b,0)}=\bigoplus_{\lambda\in\operatorname{spec}(\boldsymbol\kappa_b)}
V_{b,\lambda}.
\end{equation}
Each eigenvalue determines a point of the reduced spectral cover above $b$; the corresponding connected component is a local Hodge atom. By the local spectral decomposition theorem
\cite{KKPY}*{Theorem 4.1}, the splitting lifts to the $F$-bundle in a neighborhood of $b$ (using the graded version in the presence of odd cohomology):
\begin{equation}\label{eq:split_F_bundle}
(\mathcal H,\nabla)|_{\mathcal U} \cong \bigoplus_{\lambda\in\operatorname{spec}(\boldsymbol\kappa_b)} (\mathcal H^{\lambda},\nabla^{\lambda}).
\end{equation}
Since $\boldsymbol\kappa_b$ is $\mathsf{Hod}$-equivariant, each $V_{b,\lambda}$ is a Hodge-group representation. Its rational representation class is locally constant along the corresponding component of the generic spectral cover; see \cite{KKPY}*{Lemma 5.25}.

\begin{definition}
The \emph{coarse Hodge atom} associated with a local spectral component $\alpha$ is the isomorphism class of the rational $\mathsf{Hod}$-representation $E^\alpha$ represented by $V_{b,\lambda}$ at any point $(b,\lambda)$ of that component. The \emph{coarse atomic decomposition} is the corresponding multiset of representation classes, with multiplicities counted by the distinct generic eigenvalues. Keeping track of the local $F$-bundle factors (secured by Equation~\eqref{eq:split_F_bundle}) and their connection data gives the finer notion of Hodge atoms described in \cites{KKPY,CKK-lectures}.
\end{definition}

To each Hodge atom $\alpha$, we can associate a representation of the group $\mathsf{Hod}_{\overline{\Q}}$ in the following way: Let $\lambda$ be the eigenvalue associated with the atom $\alpha$ and $E^{\alpha}$ be the $\overline{\Q}$-linear representation of $\mathsf{Hod}_{\overline{\Q}}$ corresponding to $V_{b,\lambda}$. 

Since the isomorphism class, as a representation, of $V_{b,\lambda}$ does not depend on the point $b$, the construction is well-defined and just depends on the Hodge atom $\alpha$.

Then, we can associate the following invariants to atoms.

\begin{definition}\label{def:atominvariants}
For a coarse Hodge atom $\alpha$, write $E^\alpha_{\C}=E^\alpha\otimes_{\Q}\C$. Define
\begin{defenum}
\item Its \emph{number of Hodge cycles} by
\[
\rho_\alpha=\dim_{\C}(E^\alpha_{\C})^{\mathsf{Hod}_{\C}}.
\]
\label{item:Hodge-class}
\item Its \emph{Hodge polynomial} by
\[
P_\alpha(t)=\sum_{k\in\Z}\dim_{\C}(E^\alpha_{\C})^k\,t^k \in\Z[t,t^{-1}],
\]
where $(E^\alpha_{\C})^k$ is the weight-$k$ subspace for the Hodge cocharacter $\mathbb G_{m,\C}\to\mathsf{Hod}_{\C}$.
\end{defenum}
\end{definition}

The coarse atomic decomposition is not itself invariant under birational transformations. However, the blowup formula governs how atoms and their multiplicities change under blowups, making it useful for rationality obstructions, namely:
\begin{theorem}[\cite{KKPY}*{Proposition 5.30}]
\label{theo:atom_criteria}
Let $Z$ be a smooth projective variety of dimension $n\geq2$. If $Z$ is rational, then every coarse Hodge atom of $Z$ occurs as a coarse Hodge atom of some smooth projective variety of dimension at most $n-2$.
\end{theorem}

\subsection{Quantum differential equation}\label{subsec:rat-obstructions}
The formal quantum connection on the Hodge-invariant part provides numerical rationality obstructions when a suitable cyclic generic factor can be identified. 

Fix a local factor $(\mathcal H^\alpha,\nabla^\alpha)$ above a point $b\in U_X^{\mathrm{gen}}(\Bbbk)$. Completing along $u=0$ and then inverting $u$ produces a meromorphic differential module $\mathcal M_{b,\alpha}$ over $\Bbbk((u))$, called a \emph{generic Hodge factor}.
We say such a factor has \emph{unramified exponential type with rational residues}
if
\[
\mathcal M_{b,\alpha}\cong\mathcal E^{\lambda(b)/u}\otimes\mathcal R_{b,\lambda},
\]
where $\mathcal R_{b,\lambda}$ is regular singular with rational residue eigenvalues and $\mathcal E^{\lambda/u}$ denotes the rank-one connection $\nabla_{\partial_u}=\partial_u-\lambda/u^2$. With our convention, $\lambda(b)$ is the eigenvalue of $\boldsymbol\kappa_b$ associated with this factor. Its \emph{regularized exponent multiset} is
\[
F_{b,\lambda}:=\ms{\,-\theta+\Z\mid
\theta\text{ is a residue eigenvalue of }\mathcal R_{b,\lambda}\,},
\]
with each class counted according to its algebraic multiplicity.

A generic Hodge factor corresponding to $\lambda$ is \emph{cyclic} at $b$ if the restriction of $\boldsymbol\kappa_b$ to $V_{b,\lambda}$ has a cyclic vector. Equivalently, writing $r=\dim_\Bbbk V_{b,\lambda}$, there is $v\in V_{b,\lambda}$ such that
\[
v,\boldsymbol\kappa_bv,\ldots,
\boldsymbol\kappa_b^{r-1}v
\]
form a basis of $V_{b,\lambda}$. This is equivalent to having a single Jordan block on that generalized eigenspace.

Now we can state a criterion for proving non-rationality. Theorem~\ref{thm:criterion} below uses the formal quantum connection and its Hodge-invariant factors. 

\begin{theorem}[\cite{Leo_Luiz_Pedro_Lud_Max}*{Theorem 4.5}]\label{thm:criterion}
  Let $X$ be a connected smooth complex projective variety of
  dimension $2\leq n\leq4$. Suppose that a generic Hodge factor has a cyclic representative of unramified exponential type
  with rational residues. Let $m$ be its rank and $F$ its
  regularized exponent multiset. If $X$ is rational, then:
  \begin{enumerate}
   \item if $n=2$, then
   \[
   (m,F)=(1,\ms{\cl0});
   \]
   \item if $n=3$, then
   \[
   (m, F)\in
   \left\{
   (1,\ms{\cl0}),
   (2,\ms{\cl{1/2},\cl{1/2}})
   \right\};
   \]
   \item if $n=4$, then
   \[
   (m, F)\in
   \left\{
   (1,\ms{\cl0}),
   (2,\ms{\cl{1/2},\cl{1/2}}),
   (3,\ms{\cl0,\cl0,\cl0})
   \right\}.
   \]
  \end{enumerate}
 \end{theorem}
 \begin{remark}
For a Fano variety, Chen's exponential-type theorem \cite{Chen}*{Theorem 1.3} applies to the small quantum connection at the monotone specialization (after the usual normalization of the K\"ahler class). The extension to a chosen generic Hodge factor uses the persistence results of \cite{Leo_Luiz_Pedro_Lud_Max}. In computations, the software tests whether the unit is a cyclic vector under multiplication by the Euler vector field at the selected small-locus point, and then applies the algorithm of \cite{Leo_Luiz_Pedro_Lud_Max}*{Section 5} to obtain the regularized exponent multiset.
\end{remark}

\subsection{A new criterion in dimension five}

For fivefolds, Cavenaghi--Katzarkov--Kontsevich obtain a further obstruction from coarse Hodge atoms; we recall their statement and the main steps of their proof without claiming originality.

\begin{theorem}[Cavenaghi--Katzarkov--Kontsevich,
\cite{CKK-lectures}*{Theorem 7.7}]
\label{thm:non_rat_fivfolds} Let $X$ be a connected smooth projective complex fivefold.
If its coarse atomic decomposition contains an atom $\alpha$ with
\[
\rho_\alpha\leq3, \qquad [t^3]P_\alpha(t)>0,
\]
then $X$ is nonrational.
\end{theorem}
\begin{proof}[Sketch, following \cite{CKK-lectures}*{Chapter 7}]
Suppose that $X$ is rational. By weak factorization and the coarse atomic blowup formula, $\alpha$ occurs on a smooth projective threefold $Z$. Indeed, $[t^3]P_\alpha(t)>0$ excludes centers of dimension at most two. The threefold obstruction of \cite{CKK-lectures}*{Theorem 7.6} asserts that every coarse atom of a smooth projective threefold with $[t^3]P_\alpha(t)>0$ satisfies $\rho_\alpha\geq4$, thus contradicting $\rho_\alpha\leq3$.

For completeness, the key point is the following. After passing to a resolution $\pi\colon W\to Z_{\min}$ of a terminal minimal model, the spectral idempotent $e$ corresponding to the holomorphic three-forms satisfies $e\star P=P$, where $P$ is the point class. For an ample Cartier divisor $A$ on $Z_{\min}$, set $D=\pi^*c_1(A)$ and $C_i=D^{\cup i}$. The vanishing theorem
\cite{CKK-lectures}*{Lemma 7.5} yields, on the Hodge-fixed locus,
\[
C_i\star C_j=0\quad(i+j>3),\qquad C_i\star C_j=\delta P\quad(i+j=3), \qquad \delta=A^3>0.
\]
Consequently, $e$, $e\star C_1$, $e\star C_2$, and $e\star C_3=\delta P$ are four linearly independent Hodge-invariant vectors in the factor carrying $\alpha$. We refer to \cite{CKK-lectures}*{Lemmas 7.4--7.5 and Theorem 7.6} for the proof of these identities and the argument on terminal resolutions.
\end{proof}


\subsection{Non-further splitting}

Computing the generic atomic decomposition usually requires the big quantum product, which is difficult to determine. The small quantum operator $c_1(TX)\star_{\mathrm{small}}(-)$ is more accessible, but its spectrum need not describe the generic splitting without further hypotheses. The following criterion gives a sufficient test for the absence of additional splitting on the Hodge fixed locus.

\begin{lemma}[\cite{Leo_Luiz_Pedro_Lud_Max}*{Lemma 3.9}]
\label{lem:cyclic}
Put $N=\dim_{\Q}\rH^*_{\mathrm{Hod}}(X)$ and let $b\in B_X^{\mathsf{Hod}}(\Bbbk)$. Suppose that
\[
\mathbf 1,\mathsf{Eu}_b,\mathsf{Eu}_b^{\star2},\ldots,
\mathsf{Eu}_b^{\star(N-1)}
\]
form a basis of $T_bB_X^{\mathsf{Hod}}$. Equivalently, the minimal polynomial of $\boldsymbol\kappa_b$ on the Hodge-invariant fiber has degree $N$. Then there exists an admissible neighborhood $U\subset B_X^{\mathsf{Hod}}$ of $b$ such that:
\begin{enumerate}
\item The number of distinct eigenvalues of $\boldsymbol\kappa$ on the Hodge-invariant part, their algebraic multiplicities, and their Jordan block sizes are constant on $U$. The number of distinct eigenvalues on full cohomology is also constant on $U$.
\item If the formal meromorphic connection at $b$, obtained by completion at $u=0$ followed by inversion of $u$, is of unramified exponential type, then the same is true at every $b'\in U(\Bbbk)$. After matching eigenvalues, the corresponding regular-singular factors are isomorphic over $\Bbbk((u))$.
\end{enumerate}
\end{lemma}

\ 

We are finally ready to use the software computation jointly with the theory explained here to obtain new applications to rationality questions.


\section{Examples}\label{sec:Examples} 

Using the geometric and combinatorial tools developed in Section~\ref{sec:torics}, together with the localization formulas established in Section~\ref{sec:Loc_methods}, we have extended the software of~\cite{alexeev2026atomic} to compute genus-zero Gromov--Witten invariants of smooth projective toric varieties and their smooth complete intersections.

In particular, the implementation computes the small quantum multiplication operator by the first Chern class, as well as the cyclic determinant and the regularized exponent multisets discussed in Lemma~\ref{lem:cyclic} and Theorem~\ref{thm:criterion}, respectively. Throughout the examples, $A(q)$ is the matrix of $c_1(TY)\star_{\mathrm{small}}(-)$ on the Hodge-invariant (here ambient) cohomology, and $G$ is the matrix of the grading operator $\mathsf{Gr}$ on that space. The displayed exponent multisets refer to the corresponding Hodge-invariant differential modules.

We now present several examples illustrating the method and its applications to rationality questions. The complete collection of computations, as well as the source code, can be found at \href{https://leonardofcavenaghi.github.io/atomicdecompositions/}{GEOATOMS}; details on the software appear in Appendix~\ref{app:toricsoftware}. The smooth members of each family considered below form a connected parameter space, so the displayed Gromov--Witten invariants and small quantum multiplication matrices are deformation-invariant. Numerical specializations such as $q=2$ and $q=12$ are used only to display the eigenvalues. For the non-Archimedean $F$-bundle and Lemma~\ref{lem:cyclic}, we instead take a point in the admissible small quantum domain, for example with positive Novikov valuation. The cyclic determinants computed below are nonzero independently of the Novikov variables.

\begin{example}[Intersection of quadric and cubic in $\mathbb{P}^5$]
Let 
\[
Y = \{ f = 0\} \cap \{ g = 0\} \subset \PP^5
\]
be a smooth intersection of $f \in H^0(\cO_{\PP^5}(2))$ and $g \in H^0(\cO_{\PP^5}(3))$.

We use the realization of $Y$ as a complete intersection,
\[
Y=\mathscr Z\big(\PP^{5},\cO_{\PP^{5}}(3)\oplus \cO_{\PP^{5}}(2)\big).
\]
As $Y$ is a $(2,3)$ complete intersection, adjunction gives $\omega_Y\simeq\cO_Y(-5-1+3+2)=\cO_Y(-1)$, so $Y$ is a Fano threefold of index $1$.

Sommese's theorem tells us that
\[
H^k(\PP^5, \Q) \xrightarrow{\sim} H^k(Y, \Q) \, \, \text{for all} \, \, k<3,
\]
Thus, $H^*_{\text{amb}}(Y) = H^*_{\PP^5}(Y)= H^*_{\mathrm{Hod}}(Y)$, that is, every Hodge class comes from the ambient cohomology.

Since $H^{\ast}(\PP^5, \Q) \simeq \frac{\mathbb Q[x_{0}]}{\langle x_{0}^{6}\rangle}$, and denoting $h = c_1(\cO_Y(1))$, the ambient cohomology is generated by
\[
\{1,h,h^2,h^3\}
\]
To compute the small quantum multiplication matrix by the first Chern class, we input the $6$ rays of $\PP^5$, together with the 6 maximal cones they generate, and the line bundles that define the complete intersection as indicated in the following box.

\begin{inputbox}[Input]
rays = -1,-1,-1,-1,-1;0,0,0,0,1;0,0,0,1,0;0,0,1,0,0;0,1,0,0,0;1,0,0,0,0
cones = 0,1,2,3,5;0,1,2,3,4;0,2,3,4,5;0,1,3,4,5;0,1,2,4,5;1,2,3,4,5
K = 0,0,0,0,0,2;0,0,0,0,0,3
\end{inputbox}

\begin{outputbox}[Output]
    Matrix of $c_1(TY) \star$ in the basis $1,h,h^2,h^3$ with Novikov variable $y_1$. 
    \[
    A(y_1) = 
    \begin{bmatrix}
    0 & 792 y_{1}^{2} & 43632 y_{1}^{3} & 793152 y_{1}^{4} \\
    1 & 30 y_{1} & 2340 y_{1}^{2} & 43632 y_{1}^{3} \\
    0 & 1 & 30 y_{1} & 792 y_{1}^{2} \\
    0 & 0 & 1 & 0
    \end{bmatrix}
    \]
    Grading operator: 
    \[
    G=\operatorname{diag}\left(- \frac{3}{2}, - \frac{1}{2}, \frac{1}{2}, \frac{3}{2}\right).
    \]
    Eigenvalues at the specialization $y_1 =2$:
    \begin{itemize}
        \item $\lambda=-24$ with multiplicity 3 and exponent multiset $\{0,\frac{1}{6}, \frac{5}{6}\}$. 
    \item $\lambda=192$ with multiplicity 1 and exponent multiset $\{0\}$. 
    \end{itemize}

    Cyclic determinant of $A$: $\Delta_A = 1$.
\end{outputbox}  

The Hodge diamond of $Y$ splits on the ambient projective part, in \textcolor{blue}{blue}, and the vanishing part in \textcolor{red}{red}:
\[
\begin{array}[c]{ccccccc}
& & & \textcolor{blue}{1} & & & \\
& & 0 & & 0 & & \\
& 0 & &\textcolor{blue}{1} & & 0 & \\
0 & & \textcolor{red}{20} & & \textcolor{red}{20} & & 0 \\
& 0 & & \textcolor{blue}{1} & & 0 & \\
& & 0 & & 0 & & \\
& & & \textcolor{blue}{1} & & &
\end{array}
=
\begin{array}[c]{ccccccc}
& & & \textcolor{blue}{1} & & & \\
& & 0 & & 0 & & \\
& 0 & & \textcolor{blue}{1} & & 0 & \\
0 & & 0 & & 0 & & 0 \\
& 0 & & \textcolor{blue}{1} & & 0 & \\
& & 0 & & 0 & & \\
& & &\textcolor{blue}{1} & & &
\end{array}
+
\begin{array}[c]{ccccccc}
& & & 0 & & & \\
& & 0 & & 0 & & \\
& 0 & & 0 & & 0 & \\
0 & & \textcolor{red}{20} & & \textcolor{red}{20} & & 0 \\
& 0 & & 0 & & 0 & \\
& & 0 & & 0 & & \\
& & & 0 & & &
\end{array}
\]
Thus $h^{p,p} = 1$ for $0 \leq p \leq 3$ and $h^{1,2}(Y) = h^{2,1}(Y) =20$.

The characteristic polynomial of the small quantum multiplication matrix is
\[
\det(tI-A(y_1))
=(t-96y_1)(t+12y_1)^3.
\]
Thus, for $y_1\neq 0$, we obtain two distinct eigenvalues, one of multiplicity $1$ and the other of multiplicity $3$.

Moreover, the cyclic determinant satisfies
\[
\Delta_A(y_1)=\det\bigl(1,A(y_1)1,A(y_1)^2 1,A(y_1)^3 1\bigr)=1.
\]
Since $\Delta_A(y_1)=1$, Lemma~\ref{lem:cyclic} shows that the two distinct eigenvalues on the Hodge-invariant small quantum sector do not split further under Hodge-fixed deformations. They therefore determine two generic spectral Hodge factors, with atoms $\alpha_0$ and $\alpha_1$ satisfying
\[
\rho_{\alpha_0}=3,\qquad \rho_{\alpha_1}=1.
\]

We now apply Theorem~\ref{theo:atom_criteria}. Since every Hodge atom appearing in the atomic decomposition of a point or a curve has at most two Hodge cycles and $\rho_{\alpha_0}=3$, the atom $\alpha_0$ cannot appear in the atomic decomposition of
any smooth projective variety of dimension at most $1$. Since $\dim Y=3$, Theorem~\ref{theo:atom_criteria} implies that $Y$ is not rational.

\begin{theorem}\label{thm:quadric_cubic_P5}
Every smooth complete intersection of a quadric and a cubic in $\mathbb P^5$ is nonrational.
\end{theorem}

\end{example}
\begin{example}[Intersection of quadric and cubic in $\mathbb{P}^7$]
Now we have the smooth intersection 
\[
Y = \{f = 0\} \cap \{g = 0\} \subset \PP^7
\]
of $f \in H^0(\cO_{\PP^7}(2))$ and $g \in H^0(\cO_{\PP^7}(3))$, which can be realized as 
\[
Y=\mathscr Z\big(\PP^{7},\cO_{\PP^{7}}(3)\oplus \cO_{\PP^{7}}(2)\big).
\]
a $(2,3)$ complete intersection in $\PP^7$. By adjunction, $\omega_Y\simeq\cO_Y(-3)$, so $Y$ is a Fano fivefold of index $3$. The Lefschetz theorem for smooth complete intersections gives 
\[
H^k(\PP^7, \Q) \xrightarrow{\sim} H^k(Y, \Q) \, \, \text{for all} \, \, k<5
\]
And thus, once again, every Hodge class comes from the ambient cohomology, $H^*_{\text{amb}}(Y) = H^*_{\PP^7}(Y)= H^*_{\mathrm{Hod}}$. Denoting $h = c_1(\cO_Y(1))$, the ambient cohomology is generated by
\[
\{1,h,h^2,h^3,h^4,h^5\}
\]
The inputs to compute the small quantum multiplication matrix via the first Chern class are the 8 rays that generate the fan of $\PP^7$ together with the cones they generate and the associated line bundles that define the complete intersection
\begin{inputbox}
  rays = 1,0,0,0,0,0,0;0,1,0,0,0,0,0;0,0,1,0,0,0,0;0,0,0,1,0,0,0;0,0,0,0,1,0,0
  ;0,0,0,0,0,1,0;0,0,0,0,0,0,1;-1,-1,-1,-1,-1,-1,-1
cones = 0,1,2,3,4,5,6;0,1,2,3,4,5,7;0,1,2,3,4,6,7;0,1,2,3,5,6,7;
0,1,2,4,5,6,7;0,1,3,4,5,6,7;0,2,3,4,5,6,7;1,2,3,4,5,6,7
K = 2,0,0,0,0,0,0,0;3,0,0,0,0,0,0,0
\end{inputbox}

\begin{outputbox}[Output]
    Matrix of $c_1(TY) \star$ in the basis $1,h,h^2,h^3,h^4,h^5$ with Novikov variable $y_1$. 
    \[
    A(y_1)=\left[\begin{matrix}0 & 0 & 36 y_1 & 0 & 0 & 432 y_1^{2}\\3 & 0 & 0 & 126 y_1 & 0 & 0\\0 & 3 & 0 & 0 & 126 y_1 & 0\\0 & 0 & 3 & 0 & 0 & 36 y_1\\0 & 0 & 0 & 3 & 0 & 0\\0 & 0 & 0 & 0 & 3 & 0\end{matrix}\right].
    \]
    Grading operator:
    \[
    G=\operatorname{diag}\left( -\frac{5}{2},-\frac{3}{2},-\frac{1}{2},\frac{1}{2},\frac{3}{2},\frac{5}{2}\right)
    \]
    Eigenvalues at the specialization $y_1=2$
    \begin{itemize}
\item $\lambda=0$, with multiplicity $3$ on the ambient part and regularized exponent multiset $\left\{-\frac32,-1,-\frac12\right\}$;
        \item $\lambda=18$, with multiplicity $1$ and exponent multiset $\{0\}$;
        \item $\lambda=-9\pm9\sqrt3\,i$, each with multiplicity $1$ and exponent multiset $\{0\}$.
    \end{itemize}
    Cyclic determinant of $A$: $\Delta_A=14348907=3^{15}$.
        

\end{outputbox}
The Hodge diamond of $Y$ splits on the ambient projective part, in \textcolor{blue}{blue}, and the vanishing part in \textcolor{red}{red}:
\[
\begin{array}[c]{ccccccccccc}
 & & & & & \textcolor{blue}{1} & & & & & \\
 & & & & 0 & & 0 & & & & \\
 & & & 0 & & \textcolor{blue}{1} & & 0 & & & \\
 & & 0 & & 0 & & 0 & & 0 & & \\
 & 0 & & 0 & & \textcolor{blue}{1} & & 0 & & 0 & \\
0 & & \textcolor{red}{1} & & \textcolor{red}{83} & & \textcolor{red}{83} & & \textcolor{red}{1} & & 0 \\
 & 0 & & 0 & & \textcolor{blue}{1} & & 0 & & 0 & \\
 & & 0 & & 0 & & 0 & & 0 & & \\
 & & & 0 & & \textcolor{blue}{1} & & 0 & & & \\
 & & & & 0 & & 0 & & & & \\
 & & & & & \textcolor{blue}{1} & & & & &
\end{array}
\]
Thus $h^{p,p}(Y)=1$ for $0\leq p\leq5$, while $h^{4,1}(Y)=h^{1,4}(Y)=1$ and $h^{2,3}(Y)=h^{3,2}(Y)=83$.

The characteristic polynomial of the small quantum multiplication matrix is 
\[
\det(tI - A(y_1)) = t^3(t^3-2916 y_1).
\]
For $y_1 = 2$, we have zero as an eigenvalue of multiplicity 3. By the nonvanishing cyclic determinant and Lemma~\ref{lem:cyclic}, the three simple nonzero eigenvalues determine three point atoms, while the generalized zero eigenspace determines a factor $\alpha_0$ with three Hodge-invariant directions. Its full rank also includes the primitive middle cohomology. In this example, the odd cohomology is concentrated in degree $5$. Since $c_1(TY)=3h$ and the only nonzero odd cohomology of $Y$ is $H^5(Y)$, small quantum multiplication by $c_1(TY)$ vanishes on $H^5(Y)$. Indeed, the contribution of a curve class $\beta$ maps $H^5(Y)$ to cohomological degree \[ 7-2c_1(TY)\cdot\beta. \] For this degree to equal $5$, one would need $c_1(TY)\cdot\beta=1$, which is impossible because the Fano index is $3$. All other odd cohomological degrees vanish, and the classical product lies in $H^7(Y)=0$. In other words, we have shown that the primitive middle cohomology belongs to the zero generalized eigenspace. Because $h^{4,1}(Y)=1$, this atom satisfies $[t^3]P_{\alpha_0}(t)=1$.

\begin{theorem}[\cite{CKK-lectures}*{Theorem 7.10}]\label{theo:int_quadric_cubic}
    Every smooth intersection of a quadric and a cubic in $\mathbb{P}^{7}$ is nonrational.
\end{theorem}
\begin{proof}
The generalized zero eigenspace of the ambient multiplication matrix has dimension $3$, hence $\rho_{\alpha_0}=3$. As explained in \cite{CKK-lectures}*{Proposition 7.9}, the primitive middle cohomology belongs to the same spectral factor and has dimension \[ \dim H^5(Y,\C)=2h^{4,1}(Y)+2h^{3,2}(Y)=168. \] Thus the full factor has rank $171$, with Hodge polynomial \[ P_{\alpha_0}(t)=t^3+83t+3+83t^{-1}+t^{-3}, \] so $[t^3]P_{\alpha_0}(t)=1$ and $\rho_{\alpha_0}=3$. Theorem~\ref{thm:non_rat_fivfolds} proves the claim.
\end{proof}
\end{example}

\begin{example}[Quartic double fivefold]
Consider $f_4(x_0, \ldots,x_5)$ a homogeneous quartic polynomial with smooth zero set. $Y$ is the hypersurface given by
\[
z^2 = f_4(x_0, \ldots, x_5)
\]
Here $z$ must be a variable of weight 2, so that $Y$ sits inside a weighted projective space $Y \subset \PP(1^6,2)$.  The projection $\pi\colon Y\to\PP^5$ forgetting $z$ is a finite double cover branched along the smooth quartic hypersurface $\{f_4=0\}\subset\PP^5$. 

The point $(0:0:0:0:0:0:1)$ is a singular point of the weighted projective space, but $Y$ sits inside the smooth locus. To be more careful with our previous hypothesis, we can blow up this point and consider $Y$ inside the smooth projective toric variety $X_\Sigma = \PP_{\PP^5}(\cO \oplus \cO(2))$ with coordinates $(x_i,z,w)$. In the quotient convention, the two fiber coordinates are sections of line bundles with classes $\xi$ and $\xi-2H$, respectively: $z\in H^0(X_\Sigma,\cO_{X_\Sigma}(\xi))$ and $w\in H^0(X_\Sigma,\cO_{X_\Sigma}(\xi-2H))$. Hence both $z^2$ and $f_4(x_0,\ldots,x_5)w^2$ are sections of $\cO_{X_\Sigma}(2\xi)$. The equation $z^2=f_4w^2$ therefore defines a global divisor of class $2\xi$, disjoint from the exceptional divisor $\{w=0\}$.

Let $\rho: X_\Sigma\to\PP^5$ be the natural projection of the projective bundle
\[
X_\Sigma=\PP_{\PP^5}\bigl(\cO_{\PP^5}\oplus\cO_{\PP^5}(2)\bigr).
\]
We denote by $H=\rho^*c_1(\cO_{\PP^5}(1))$ the pullback of the hyperplane class of $\PP^5$, and by $\xi=c_1(\cO_{X_\Sigma}(1))$ the relative hyperplane class of the projective bundle, using the quotient convention. These two classes generate the Picard group
\[
\operatorname{Pic}(X_\Sigma)=\Z H\oplus\Z\xi.
\]
Geometrically, $H$ corresponds to the inverse image of a hyperplane in $\PP^5$, while $\xi$ restricts to the hyperplane class on each fiber $\rho^{-1}(p)\simeq\PP^1$. The exceptional divisor of the blowup $X_\Sigma\to\PP(1^6,2)$ has class $\xi-2H$.

And we realize $Y$ as a general locus 
\[
Y=\mathscr Z\big(X_\Sigma,\cO_{X_\Sigma}( 2 \xi) \big).
\]
The defining line bundle is nef but not ample; therefore, Sommese's theorem does not apply. However, using the geometry of the variety, we can still understand the structure of $Y$ in the ambient toric space.

Since $Y$ is a smooth quartic hypersurface in $\PP(1^6,2)$ avoiding the singular point, the Lefschetz theorem for smooth weighted-projective hypersurfaces gives $H^{2p}(Y,\Q)=\Q h^p$ for $0\leq p\leq5$, with the remaining cohomology concentrated in $H^5(Y,\Q)$. In particular, the map $i^*\colon H^2(X_\Sigma,\Q)\to H^2(Y,\Q)$ is not an isomorphism: its source has dimension $2$ and its target has dimension $1$. More precisely, $i^*(\xi)=2i^*(H)$; we write $h=i^*(H)=c_1(\cO_Y(1))$. The ambient cohomology is generated by
\[
\{1,h,h^2,h^3,h^4,h^5\}.
\]
Hence, the ambient cohomology is the span of these six even-degree classes. As the remaining cohomology is concentrated in the odd middle degree, there are no further Hodge-invariant classes; thus $H^*_{\mathrm{amb}}(Y)=H^*_{\mathrm{Hod}}(Y)$.

Since the first Chern class of $X_\Sigma$ is
\[
c_1(TX_\Sigma)=2\xi+4H,
\]
the adjunction formula gives
\[
K_Y=(K_{X_\Sigma}+[Y])|_Y
=(-2\xi-4H+2\xi)|_Y=-4h.
\]
Therefore,
\[
c_1(TY)=4h,
\]
and $Y$ is a Fano fivefold of index $4$.

Although $i_*\colon H_2(Y,\Q)\to H_2(X_\Sigma,\Q)$ is not surjective, it is injective. Indeed, $H_2(Y,\Q)$ has dimension $1$, and $i^*H=h\neq0$ implies that $i_*$ is nonzero by the projection formula. Equivalently, this map is injective on numerical curve classes. Consequently, each ambient curve class has at most one preimage among effective classes of $Y$, and Theorem~\ref{thm:localizationCI} computes individual Gromov--Witten invariants rather than only sums over distinct preimages.

Thus, we can use the software to compute the small quantum multiplication by the first Chern class, using the following inputs: the 8 rays that generate the fan of $X_\Sigma$, the cones they generate, and the associated line bundles that define the complete intersection.
\begin{inputbox}
rays = 1,0,0,0,0,0;0,1,0,0,0,0;0,0,1,0,0,0;0,0,0,1,0,0;0,0,0,0,1,0;
-1,-1,-1,-1,-1,2;0,0,0,0,0,1;0,0,0,0,0,-1
cones = 0,1,2,3,4,6;0,1,2,3,4,7;0,1,2,3,5,6;0,1,2,3,5,7;0,1,2,4,5,6;0,1,2,4,5,7;
0,1,3,4,5,6;0,1,3,4,5,7;0,2,3,4,5,6;0,2,3,4,5,7;1,2,3,4,5,6;1,2,3,4,5,7
K = 0,0,0,0,0,0,0,2
\end{inputbox}

\begin{outputbox}[Output]
     Matrix of $c_1(TY) \star$ in the basis $1,h,h^2,h^3,h^4,h^5$ with Novikov variable $q=y_1^2y_2$.
    \[A(q)=\left[\begin{matrix}0 & 0 & 0 & 48 q & 0 & 0\\4 & 0 & 0 & 0 & 160 q & 0\\0 & 4 & 0 & 0 & 0 & 48 q\\0 & 0 & 4 & 0 & 0 & 0\\0 & 0 & 0 & 4 & 0 & 0\\0 & 0 & 0 & 0 & 4 & 0\end{matrix}\right]\]
    The grading operator is given by 
    \[
     G=\operatorname{diag}\!\left(-\tfrac52,-\tfrac32,-\tfrac12,\tfrac12,\tfrac32,\tfrac52\right).
    \]
    Eigenvalues at specialization $q=12$
    \begin{itemize}
        \item $\lambda=0$ has multiplicity $2$ on ambient cohomology, with exponent multiset $\left\{-\frac32,\frac12\right\}$ on the Hodge-invariant factor.
        \item The other four eigenvalues are simple with exponent multiset $\{0\}$ at each root. 
    \end{itemize}
    Cyclic determinant of $A$: $\Delta_A=1073741824=4^{15}$.
\end{outputbox}

An important remark here is that the software works with the Novikov variables of the toric variety $X_\Sigma$, denoted by $y_1$ and $y_2$. The relation between the Novikov variables of $Y$ and
those of $X_\Sigma$ follows from the restriction map
$i^*:H^2(X_\Sigma,\Z)\to H^2(Y,\Z)$ and its dual
$i_*:H_2(Y,\Z)\to H_2(X_\Sigma,\Z)$, via the natural
cohomology-homology pairing. Since
\[
i^*(H)=h,\qquad i^*(\xi)=2h,
\]
we obtain
\[
i_*(\gamma)=2\beta_f+\beta_e,
\]
where $\gamma$ generates $H_2(Y,\Z)$ and
$\beta_f,\beta_e$ are the curve classes dual
to $\xi,H$, respectively. Consequently, if $y_1$ and $y_2$ denote the ambient Novikov coordinates dual to $\xi$ and $H$, respectively, the Novikov variable associated with $\gamma$ is
\[
q=y_1^2y_2.
\]

The Hodge diamond splits on the ambient projective part, in \textcolor{blue}{blue}, and the vanishing part, in \textcolor{red}{red}
\[
\begin{array}{ccccccccccc}
&&&&&\textcolor{blue}{1}&&&&&\\
&&&&0&&0&&&&\\
&&&0&&\textcolor{blue}{1}&&0&&&\\
&&0&&0&&0&&0&&\\&0&&0&&\textcolor{blue}{1}&&0&&0&\\0&&\textcolor{red}{1}&&\textcolor{red}{90}&&\textcolor{red}{90}&&\textcolor{red}{1}&&0\\&0&&0&&\textcolor{blue}{1}&&0&&0&\\&&0&&0&&0&&0&&\\&&&0&&\textcolor{blue}{1}&&0&&&\\&&&&0&&0&&&&\\
&&&&&\textcolor{blue}{1}&&&&&
\end{array}
\]
Hence, $h^{4,1}(Y)=h^{1,4}(Y)=1$, $h^{3,2}(Y)=h^{2,3}(Y)=90$.
The characteristic polynomial of the small quantum multiplication is
\[
\det(tI-A(q))=t^2(t^4-16384q).
\]
Thus, for $q\neq0$, the ambient operator has four simple nonzero eigenvalues and a generalized zero eigenspace of dimension $2$. Since $\Delta_A\neq0$, Lemma~\ref{lem:cyclic} rules out further splitting at a generic Hodge-fixed point. The only nonzero odd cohomology is $H^5(Y)$, and quantum multiplication by $c_1(TY)=4h$ vanishes on it: a curve class of $h$-degree $m\geq0$ would send it to degree $7-8m$, where the odd cohomology vanishes. Thus $H^5(Y)$ belongs to the zero factor $\alpha_0$, which has $\rho_{\alpha_0}=2$ Hodge-invariant directions. The four simple factors are point atoms. 

\begin{theorem}\label{thm:non_rat_quartic_double_fivefold}
    Every smooth quartic double fivefold is nonrational.
\end{theorem}
\begin{proof}
    The zero-eigenvalue factor has $\rho_{\alpha_0}=2$ Hodge-invariant directions. Since $\dim H^5(Y,\C)=2(1+90)=182$, its full rank is $184$ and its Hodge polynomial is \[ P_{\alpha_0}(t)=t^3+90t+2+90t^{-1}+t^{-3}. \] In particular, $[t^3]P_{\alpha_0}(t)=1>0$ and $\rho_{\alpha_0}=2\leq3$. Theorem~\ref{thm:non_rat_fivfolds} therefore implies the nonrationality of $Y$.
\end{proof}
\end{example}

\begin{example}[Quartic fivefold in $\mathbb{P}^{6}$]
\label{ex:quartic_fivefold}
Let
\[
Y=\{f_4=0\}\subset\PP^{6}
\]
be a smooth hypersurface defined by a homogeneous polynomial
$f_4\in H^0(\PP^{6},\cO_{\PP^{6}}(4))$. We use the realization
\[
Y=\mathscr Z\bigl(\PP^{6},\cO_{\PP^{6}}(4)\bigr).
\]
Writing $h=c_1(\cO_Y(1))$, adjunction gives
\[
\omega_Y\simeq\cO_Y(-7+4)=\cO_Y(-3),\qquad c_1(TY)=3h.
\]
Thus $Y$ is a Fano fivefold of index $3$.

By the Lefschetz hyperplane theorem,
\[
H^k(\PP^{6},\Q)\xrightarrow{\sim}H^k(Y,\Q)
\qquad\text{for }k<5.
\]
Together with Poincar\'e duality, this identifies the even cohomology
with the ambient cohomology. Since the vanishing cohomology lies in
the odd degree $5$, every rational Hodge class is ambient:
\[
H^*_{\mathrm{amb}}(Y)=H^*_{\PP^{6}}(Y)=H^*_{\mathrm{Hod}}(Y).
\]
We use the ordered basis
\[
\{1,h,h^2,h^3,h^4,h^5\}.
\]
The input for the software consists of the seven rays and seven maximal cones of the
fan of $\PP^{6}$, together with the divisor defining the quartic.

\begin{inputbox}[Input]
rays = 1,0,0,0,0,0;0,1,0,0,0,0;0,0,1,0,0,0;0,0,0,1,0,0;0,0,0,0,1,0;0,0,0,0,0,1;
-1,-1,-1,-1,-1,-1
cones = 0,1,2,3,4,5;0,1,2,3,4,6;0,1,2,3,5,6;0,1,2,4,5,6;0,1,3,4,5,6;0,2,3,4,5,6;
1,2,3,4,5,6
K = 4,0,0,0,0,0,0
\end{inputbox}

\begin{outputbox}[Output]
Matrix of $c_1(TY)\star$ in the basis $1,h,h^2,h^3,h^4,h^5$,
with Novikov variable $y_1$:
\[
A(y_1)=
\begin{bmatrix}
0&0&72y_1&0&0&1728y_1^2\\
3&0&0&312y_1&0&0\\
0&3&0&0&312y_1&0\\
0&0&3&0&0&72y_1\\
0&0&0&3&0&0\\
0&0&0&0&3&0
\end{bmatrix}.
\]
Grading operator:
\[
G=\operatorname{diag}\!\left(
-\frac52,-\frac32,-\frac12,\frac12,\frac32,\frac52
\right).
\]
 At $y_1=2$, the eigenvalues of the ambient matrix $A(y_1)$ and the corresponding regularized exponent multisets are:
\begin{itemize}
\item $\lambda=0$, with multiplicity $3$ and exponent multiset
$\left\{-\frac74,-1,-\frac14\right\}$;
\item $\lambda=24$, with multiplicity $1$ and exponent multiset $\{0\}$;
\item $\lambda=-12\pm12\sqrt3\,i$, each with multiplicity $1$ and
exponent multiset $\{0\}$.
\end{itemize}

The cyclic determinant is
\[
\Delta_A=3^{15}=14348907.
\]
\end{outputbox}

The Hodge diamond separates the ambient projective part, in
\textcolor{blue}{blue}, from the primitive part, in
\textcolor{red}{red}:
\[
\begin{array}{ccccccccccc}
&&&&&\textcolor{blue}{1}&&&&&\\
&&&&0&&0&&&&\\
&&&0&&\textcolor{blue}{1}&&0&&&\\
&&0&&0&&0&&0&&\\
&0&&0&&\textcolor{blue}{1}&&0&&0&\\
0&&\textcolor{red}{7}&&\textcolor{red}{266}&&\textcolor{red}{266}&&\textcolor{red}{7}&&0\\
&0&&0&&\textcolor{blue}{1}&&0&&0&\\
&&0&&0&&0&&0&&\\
&&&0&&\textcolor{blue}{1}&&0&&&\\
&&&&0&&0&&&&\\
&&&&&\textcolor{blue}{1}&&&&&
\end{array}
\]

The characteristic polynomial of the small quantum multiplication
matrix is
\[
\det(tI-A(y_1))=t^3\bigl(t^3-6912y_1\bigr).
\]
For $y_1\neq0$, there are three simple nonzero eigenvalues and a zero
eigenvalue of multiplicity $3$. Moreover,
\[
\Delta_A(y_1)
=\det\bigl(1,A(y_1)1,A(y_1)^2 1,\ldots,A(y_1)^5 1\bigr)
=3^{15}=14348907.
\]

In the framework of Lemma~\ref{lem:cyclic}, the nonvanishing cyclic
determinant gives three simple Hodge factors and a zero-eigenvalue
Hodge atom $\alpha_0$ with
\[
\rho_{\alpha_0}=3.
\]
Since $c_1(TY)=3h$, a curve class $\beta$ of $h$-degree $m\geq0$ sends $H^5(Y)$ under quantum multiplication by $c_1(TY)$ to degree $7-6m$. This cohomology group vanishes for every $m\geq0$: it is $H^7(Y)=0$ when $m=0$, $H^1(Y)=0$ when $m=1$, and has negative degree when $m\geq2$. Hence $H^5(Y)$ lies in the generalized zero-eigenvalue factor $\alpha_0$. Its Hodge polynomial is \[ P_{\alpha_0}(t)=7t^3+266t+3+266t^{-1}+7t^{-3}, \] so $\rho_{\alpha_0}=3$, $[t^3]P_{\alpha_0}(t)=7>0$, and the full factor has rank $3+\dim H^5(Y,\C)=549$. Theorem~\ref{thm:non_rat_fivfolds} gives the claimed nonrationality.

\begin{theorem}
\label{thm:quartic_fivefold}
Every smooth quartic fivefold in $\PP^{6}$ is nonrational.
\end{theorem}

\end{example}

    \bibliographystyle{alpha}
\bibliography{min.bib}

\bigskip

\appendix

\section{Computing the fan of a toric variety}\label{ComputeFan}

In this appendix, we outline how to recover the fan $\Sigma$ of a toric variety $X_\Sigma$ from the geometry of its torus action.

By Theorem~\ref{OrbitCone}, there is a bijective correspondence between the cones $\sigma \in \Sigma$ and the $\T_N$-orbits $O(\sigma) \subseteq X_\Sigma$, under which $\tau \preceq \sigma$ if and only if $O(\sigma) \subseteq \overline{O(\tau)}$. Since $O(\{0\}) = \T_N$ is the unique open dense orbit of dimension $n$, its complement is the union of the closures of the $(n-1)$-dimensional orbits corresponding to the rays $\rho \in \Sigma(1)$:
\[
X_\Sigma \setminus \T_N = \bigcup_{\rho \in \Sigma(1)} D_\rho, \qquad \text{where } D_\rho = \overline{O(\rho)}.
\]
Thus, the rays of $\Sigma$ correspond bijectively to the codimension-$1$ irreducible components of $X_\Sigma \setminus \T_N$, which are precisely the $\T_N$-invariant prime divisors $D_\rho$.

To reconstruct the higher-dimensional cones of $\Sigma$, we must characterize when a collection of rays $\rho_{i_1}, \ldots, \rho_{i_k} \in \Sigma(1)$ spans a cone in $\Sigma$. Recall from Theorem~\ref{OrbitCone}(d) that for each $\tau \in \Sigma$,
\[
\overline{O(\tau)} = \bigcup_{\substack{\sigma \in \Sigma \\ \tau \preceq \sigma}} O(\sigma),
\]
and that distinct $\T_N$-orbits in $X_\Sigma$ are disjoint. Consequently, a point $x \in X_\Sigma$ lies in the intersection $D_{i_1} \cap \cdots \cap D_{i_k}$ if and only if the unique orbit $O(\sigma)$ containing $x$ satisfies $O(\sigma) \subseteq D_{i_j} = \overline{O(\rho_{i_j})}$ for every $j = 1, \ldots, k$. By the Orbit--Cone Correspondence, this holds if and only if each $\rho_{i_j}$ is a face of $\sigma$. Hence, $D_{i_1} \cap \cdots \cap D_{i_k} \neq \varnothing$ if and only if the rays $\rho_{i_1}, \ldots, \rho_{i_k}$ are contained in a common cone $\sigma \in \Sigma$. In general, each cone $\sigma \in \Sigma$ is recovered from its orbit closure via $\sigma = \operatorname{Cone}(u_i \mid \overline{O(\sigma)} \subseteq D_i)$. When the fan $\Sigma$ is simplicial (in particular, when $X_\Sigma$ is smooth), every subset of the ray generators of a cone $\sigma \in \Sigma$ spans a face of $\sigma$, so that
\[
\operatorname{Cone}(u_{i_1}, \ldots, u_{i_k}) \in \Sigma
\quad \Longleftrightarrow \quad
D_{i_1} \cap \cdots \cap D_{i_k} \neq \varnothing,
\]
where $u_{i_1}, \ldots, u_{i_k} \in N$ denote the primitive ray generators of $\rho_{i_1}, \ldots, \rho_{i_k}$, respectively. It therefore suffices to determine the ray generators $u_\rho$.

\begin{proposition}[\cite{cox2011toric}*{Proposition 4.1.1}]
    Let $X_\Sigma$ be the toric variety associated to a fan $\Sigma$. If $\rho \in \Sigma(1)$ is a ray with primitive generator $u_\rho \in N$ and $\chi^m$ is the character corresponding to $m \in M$, then
    \[
    \nu_\rho(\chi^m) = \langle m, u_\rho \rangle,
    \]
    where $\nu_\rho = \nu_{D_\rho} \colon \C(X_\Sigma)^{\times} \to \Z$ denotes the discrete valuation of the function field $\C(X_\Sigma) = \C(\T_N)$ associated to the prime divisor $D_\rho$.
\end{proposition}

Choosing a $\Z$-basis $m_1, \ldots, m_n$ of $M$ and identifying $N = \operatorname{Hom}_{\Z}(M, \Z) \simeq \Z^n$ via the dual $\Z$-basis of $N$, the coordinates of the primitive ray generator $u_\rho$ are given explicitly by
\[
u_\rho = \bigl(\nu_\rho(\chi^{m_1}), \ldots, \nu_\rho(\chi^{m_n})\bigr) \in \Z^n.
\]

\begin{remark}
    The correspondence between $\T_N$-invariant prime divisors and rays via discrete valuations relies only on the assumption that the toric variety $X_\Sigma$ is normal (see \cite{cox2011toric}*{Theorem 3.1.5 and Proposition 4.1.1}), whereas the criterion $\operatorname{Cone}(u_{i_1}, \ldots, u_{i_k}) \in \Sigma \Longleftrightarrow D_{i_1} \cap \cdots \cap D_{i_k} \neq \varnothing$ uses that $\Sigma$ is simplicial (for instance, when $X_\Sigma$ is smooth).
\end{remark}

For a smooth (or simplicial) toric variety $X$, the procedure for recovering the fan $\Sigma$ can be summarized in four steps:
\begin{enumerate}
    \item Determine the codimension-$1$ irreducible components of $X \setminus \T_N$, which are the $\T_N$-invariant prime divisors $D_\rho$.
    \item Fix a $\Z$-basis $m_1, \ldots, m_n$ of the character lattice $M$ and compute the order of vanishing $\nu_\rho(\chi^{m_j})$ of each basis character $\chi^{m_j}$ along each divisor $D_\rho$.
    \item Form the primitive ray generators $u_\rho = (\nu_\rho(\chi^{m_1}), \ldots, \nu_\rho(\chi^{m_n})) \in N \simeq \Z^n$ and the rays $\rho = \operatorname{Cone}(u_\rho)$.
    \item Determine which collections of divisors $D_{i_1}, \ldots, D_{i_k}$ have non-empty intersection $D_{i_1} \cap \cdots \cap D_{i_k} \neq \varnothing$; each such collection spans a cone $\operatorname{Cone}(u_{i_1}, \ldots, u_{i_k}) \in \Sigma$.
\end{enumerate}

\begin{example}
    Let us apply this algorithm to $\PP^2$ with homogeneous coordinates $[x_0 : x_1 : x_2]$. Here the torus is $\T_N = \{[x_0 : x_1 : x_2] \in \PP^2 \mid x_0 x_1 x_2 \neq 0\}$, so its complement decomposes into three irreducible components of codimension $1$:
    \[
    \PP^2 \setminus \T_N = \{x_0 x_1 x_2 = 0\} = D_0 \cup D_1 \cup D_2,
    \]
    where the $\T_N$-invariant prime divisors are the coordinate lines
    \[
    D_0 = \{x_0 = 0\}, \qquad D_1 = \{x_1 = 0\}, \qquad D_2 = \{x_2 = 0\}.
    \]
    On the affine chart $U_0 = \{x_0 \neq 0\} \simeq \C^2$, the torus $\T_N$ is identified with $(\C^{\times})^2$ via the coordinates $t_1 = x_1/x_0$ and $t_2 = x_2/x_0$. The characters of $\T_N$ are given by $\chi^m = t_1^a t_2^b$ for $m = (a,b) \in M = \Z^2$. Taking the standard basis $m_1 = (1,0)$ and $m_2 = (0,1)$ of $M$, the corresponding basis characters are the rational functions
    \[
    \chi^{m_1} = t_1 = \frac{x_1}{x_0}
    \qquad \text{and} \qquad
    \chi^{m_2} = t_2 = \frac{x_2}{x_0}.
    \]
    Computing the order of vanishing or pole of $\chi^{m_1}$ and $\chi^{m_2}$ along each divisor $D_i = \{x_i = 0\}$ yields
    \[
    \begin{aligned}
        \nu_{0}(\chi^{m_1}) &= -1, &\qquad \nu_{1}(\chi^{m_1}) &= 1, &\qquad \nu_{2}(\chi^{m_1}) &= 0, \\
        \nu_{0}(\chi^{m_2}) &= -1, &\qquad \nu_{1}(\chi^{m_2}) &= 0, &\qquad \nu_{2}(\chi^{m_2}) &= 1.
    \end{aligned}
    \]
    Consequently, the primitive generators of the rays $\rho_0, \rho_1, \rho_2$ associated to $D_0, D_1, D_2$ are
    \[
    u_0 = (-1,-1), \qquad u_1 = (1,0), \qquad u_2 = (0,1).
    \]
    Finally, we examine the intersections of $D_0, D_1, D_2$. Every pair of distinct divisors intersects at a torus fixed point, whereas the triple intersection is empty ($D_0 \cap D_1 \cap D_2 = \varnothing$):
    \[
    \begin{aligned}
        D_1 \cap D_2 &= \{[1:0:0]\} \neq \varnothing &&\longleftrightarrow \sigma_0 = \operatorname{Cone}(u_1,u_2), \\
        D_0 \cap D_2 &= \{[0:1:0]\} \neq \varnothing &&\longleftrightarrow \sigma_1 = \operatorname{Cone}(u_0,u_2), \\
        D_0 \cap D_1 &= \{[0:0:1]\} \neq \varnothing &&\longleftrightarrow \sigma_2 = \operatorname{Cone}(u_0,u_1).
    \end{aligned}
    \]
    Together with their faces $\rho_0, \rho_1, \rho_2$ and $\{0\}$, the maximal cones $\sigma_0, \sigma_1, \sigma_2$ form the fan $\Sigma$ of $\PP^2$.

    Although this valuation procedure applies to any normal toric variety, certain families admit a more direct construction. For projective space $\PP^r$, let $e_0, \ldots, e_r$ denote the standard basis of $\R^{r+1}$. The positive orthant $\operatorname{Cone}(e_0, \ldots, e_r) \subseteq \R^{r+1}$ has $r+1$ facets
    \[
    F_i = \operatorname{Cone}(e_0, \ldots, \widehat{e}_i, \ldots, e_r), \qquad i = 0, \ldots, r,
    \]
    where $\widehat{e}_i$ indicates that $e_i$ is omitted. Define the quotient lattice
    \[
    \overline{N} = \Z^{r+1} / \Z(e_0 + \cdots + e_r).
    \]
    In $\overline{N}$, the images $\overline{e}_0, \ldots, \overline{e}_r$ of the basis vectors satisfy $\overline{e}_0 + \cdots + \overline{e}_r = 0$ (and any $r$ of them form a $\Z$-basis of $\overline{N} \simeq \Z^r$). The images $\overline{F}_i \subseteq \overline{N}_{\R}$ of the facets $F_i$, together with all of their faces, form the fan of $\PP^r$.

Moreover, this quotient perspective serves as the prototype for constructing the fans of projectivizations of split vector bundles
\[
\PP\left( \mathcal{O}_{X_\Sigma}(D_0) \oplus \cdots \oplus \mathcal{O}_{X_\Sigma}(D_r) \right)
\]
over a toric variety $X_\Sigma$, where $D_0,\ldots,D_r$ are torus-invariant Cartier divisors. We refer to \cite{cox2011toric}*{Proposition 7.3.3} for the corresponding fan construction.
\end{example}

\section{On the software}
\label{app:toricsoftware}

In this appendix, we detail what the software computes and explain why the assumptions of Theorem \ref{thm:ci-toric} provide an ideal setting for its application. We then discuss how the software can be used when these assumptions are not satisfied. Finally, we explain how to provide the necessary input data to the code.

\medskip

To compute Gromov--Witten invariants for a smooth projective toric variety $X_{\Sigma}$, the implemented software uses the localization
formula proved in Theorem~\ref{thm:localization}. The inputs for the calculation are the cohomology classes, identified with polynomials in the divisor classes $x_i=[D_i]$ indexed by the rays of $\Sigma$, as in Section~\ref{sec:picard}, and an effective class $\beta\in\rH_2(X_\Sigma,\mathbb{Z})$, given as an integral relation among the primitive ray generators, as in Proposition~\ref{prop:torich2}.

For a complete intersection
$Y:=\mathscr{Z}(X_\Sigma,\mathcal{E})$, the formula implemented is the one proved
in Theorem~\ref{thm:localizationCI}. The bundles supported by the
software are direct sums of nontrivial nef line bundles,
\[
\mathcal{E}
=
\bigoplus_{\ell=1}^{e}\mathcal{O}_{X_\Sigma}(D_{\ell}).
\]
Their weights are determined by the divisor coefficients and the fan,
as explained in Proposition~\ref{weightsbundles}.

The inputs required to compute the sum of Gromov--Witten invariants
are cohomology classes coming from the ambient part
\[
\rH_{X_\Sigma}(Y):=\operatorname{im}
\bigl(\iota^*:\rH^*(X_\Sigma,\mathbb{Q})
\longrightarrow\rH^*(Y,\mathbb{Q})\bigr),
\]
represented by restrictions of the divisor-polynomial classes described in Theorem~\ref{12.4.1}, together with a homology class $\beta\in\rH_2(X_\Sigma,\mathbb{Z})$, given as an integral relation among the ray generators.

\ 

Assume that the Poincaré pairing restricts nondegenerately to $\rH_{X_\Sigma}(Y)$, and let $\{\sigma_i\}$ be a basis of this space. The software computes the matrix of the ambient projection of small quantum multiplication by $c_1(TY)$, with Novikov variables identified through $\iota_*$:
\[
A(q)\sigma_i
=
\sum_{\beta\in\mathsf{NE}(X_\Sigma)}
\sum_{j,k}
\left(
\sum_{\substack{
\beta'\in\mathsf{NE}(Y)\\
\iota_*\beta'=\beta}}
\big\langle
c_1(TY),\sigma_i,\sigma_j
\big\rangle^Y_{0,\beta'}
\right)
q^\beta g^{jk}\sigma_k,
\]
where
\[
g_{jk}=\int_Y\sigma_j\cup\sigma_k,
\]
and $(g^{jk})$ is the inverse of $(g_{jk})$.

Notice that the Novikov variables are indexed by the second homology
of the toric variety, not of the complete intersection. More explicitly,
for the divisor basis $A_1,\ldots,A_r$ selected by the software,
\[
q^\beta=\prod_{\nu=1}^{r}y_\nu^{A_\nu\cdot\beta}.
\]
The software uses a nef integral basis when it finds one; otherwise, it uses its chosen divisor basis. In the latter case, some exponents $A_\nu\cdot\beta$ may be negative, and the corresponding expressions are Laurent monomials.

Under the ampleness and dimension assumptions of
Theorem~\ref{thm:ci-toric}, namely that every defining divisor $D_\ell$
is ample and $\dim Y\geq3$, there are isomorphisms
\[
\rH^2(X_\Sigma,\mathbb{Q})\cong\rH^2(Y,\mathbb{Q}),
\qquad
\iota_*:\rH_2(Y,\mathbb{Z})
\xrightarrow{\sim}\rH_2(X_\Sigma,\mathbb{Z}).
\]
Thus, Theorem~\ref{thm:localizationCI} actually computes Gromov--Witten invariants. Under these assumptions, the ambient cohomology is preserved by small quantum multiplication by $c_1(TY)$, so $A(q)$ describes the restriction of this operator to the ambient part. Classes outside the ambient part are not included in this matrix.

When the map
$\iota_*:\rH_2(Y,\mathbb{Z})\rightarrow\rH_2(X_\Sigma,\mathbb{Z})$
fails to be an isomorphism, one should keep track of its behavior
to understand the relation between $A(q)$ and
$c_1(TY)\star_{\mathrm{small}}$. In particular, effective classes of
$Y$ with the same image in $X_\Sigma$ contribute to the same Novikov monomial.

There is also a distinction concerning the enumeration of curve classes.
If
\[
L=-K_{X_\Sigma}-\sum_{\ell=1}^{e}D_\ell
\]
is ample on $X_\Sigma$, the dimension constraint and the positivity of $L$ on the ambient Mori cone give a finite list of classes that can contribute to $A(q)$. When $Y$ is Fano but $L$ is not ample on $X_\Sigma$, Proposition~\ref{prop:finite enum} provides a ``finite search region'' under its additional positivity hypothesis.

The complete-intersection calculation requires every defining line
bundle to be nef, equivalently, globally generated on $X_\Sigma$.
The software checks this condition by testing its degree on every
torus-invariant curve.

\medskip

Lastly, we briefly describe how to use the software.

\medskip

The inputs required for the computation are the following: the
primitive ray generators and maximal cones, which characterize the toric variety; and the defining divisors, which specify the vector bundle.

First, one must specify the primitive generators of
the rays of $\Sigma$. We number these starting from zero, in accordance
with the convention used by the code:
\[
v_0,\ldots,v_{k-1}\in N\cong\mathbb{Z}^n.
\]
Each vector is entered as a comma-separated list of coordinates.
Different vectors are separated by a semicolon ``;'' or entered on
separate lines in the ``Rays'' field. Thus, for example,
\[
\texttt{1,0;0,1;-1,-1}
\]
specifies the three primitive ray generators
\[
v_0=(1,0),\qquad v_1=(0,1),\qquad v_2=(-1,-1).
\]

Second, we must specify the maximal cones of the fan. A maximal cone
\[
\sigma=\operatorname{Cone}(v_{i_1},\ldots,v_{i_n})
\]
is recorded by the corresponding indices
\[
i_1,\ldots,i_n.
\]
Different maximal cones are separated by a semicolon ``;'' or entered
on separate lines in the ``Maximal cones'' field. For the rays above,
the input
\[
\texttt{0,1;1,2;0,2}
\]
gives the fan of $\mathbb{P}^2$.
The input fan must define a smooth projective toric variety.

Lastly, we need to specify the vector bundle under consideration.
We do so in the ``Complete intersection'' field, corresponding to
the parameter \texttt{K}. Let $D_i$ denote the invariant prime divisor
corresponding to the ray generated by $v_i$. The line bundle
\[
\mathcal{O}_{X_\Sigma}(D),
\qquad
D=\sum_{i=0}^{k-1}a_iD_i,
\]
is encoded by
\[
a_0,\ldots,a_{k-1}.
\]
Different summands are separated by a semicolon ``;''. Thus, the
vector bundle
\[
\mathcal{O}_{X_\Sigma}\left(\sum_{i=0}^{k-1}a_iD_i\right)
\oplus
\mathcal{O}_{X_\Sigma}\left(\sum_{i=0}^{k-1}b_iD_i\right)
\]
is entered as
\[
a_0,\ldots,a_{k-1};b_0,\ldots,b_{k-1}.
\]

Alternatively, a reduced expression can be used. If
$x_{j_1},\ldots,x_{j_r}$ is the basis of $\rH^2(X_\Sigma)$ displayed by
the software, where $r=k-n$, a divisor class
\[
[D]=d_1x_{j_1}+\cdots+d_rx_{j_r}
\]
can be entered as
\[
d_1,\ldots,d_r.
\]
For example, on $\mathbb{P}^n$, the line bundle
$\mathcal{O}(d)$ is entered as \texttt{d}, and
$\mathcal{O}(a)\oplus\mathcal{O}(b)$ is entered as
\texttt{a;b}. Leaving this field empty gives computations for
the ambient toric variety itself.

For an individual Gromov--Witten calculation, the cohomology insertions
are entered as polynomial expressions in the divisor classes
\[
x_i=[D_i],
\]
written as \texttt{x0}, \texttt{x1}, and so on.
The notation \texttt{one} denotes the unit class. The symbol \texttt{pt} denotes the point class of the ambient variety. Different insertions are separated by \texttt{|}; for example, \texttt{pt|pt} specifies two point insertions. The curve class is entered as a vector
\[
\beta=(b_0,\ldots,b_{k-1})
\]
satisfying
\[
\sum_{i=0}^{k-1}b_iv_i=0,
\qquad
\int_\beta x_i=b_i.
\]
Alternatively, one may enter its coordinates in the integral basis
of $\rH_2(X_\Sigma,\mathbb{Z})$ displayed by the software.
The coefficients $b_i$ need not all be nonnegative for an effective
curve class.

Furthermore, there are shortcuts implemented for any projective space, products of such, and projectivization of vector bundles over projective spaces of the form
\[
\mathbb P_{\mathbb P^n}\!\left( \bigoplus_{i=1}^{k}\mathcal O_{\mathbb P^n}(a_i) \right),\qquad a_i\ge0\quad(1\le i\le k).
\]

For implementation purposes, we have to establish a criterion for when the software should stop summing over effective curve classes. Let us first recall why this problem is simpler for the complete intersections in flag varieties considered in \cite{alexeev2026atomic}. 

Let $F=G/P$, and let $Y\subset F$ be a smooth Fano complete intersection of dimension $d$, defined as the zero locus of a regular section of a globally generated vector bundle $\mathcal E$. Denote by $\iota\colon Y\hookrightarrow F$ the inclusion and set
\[
L=c_1(TF)-c_1(\mathcal E),
\qquad
\iota^*L=c_1(TY).
\]
The cone of effective curves of $F$ is generated by Schubert curve classes $[C_1],\ldots,[C_r]$, which form an integral basis of $H_2(F,\Z)$. Thus, the ambient image of every effective curve class $\beta'\in H_2(Y,\Z)$ can be written uniquely as
\[
\beta=\iota_*\beta'
=\sum_{i=1}^r b_i[C_i],
\qquad b_i\in\Z_{\geq0}.
\]
For three-point invariants contributing to small quantum multiplication
by $c_1(TY)$, the dimension constraint gives
\[
0\leq c_1(TY)\cdot\beta'
=L\cdot\beta
=\sum_{i=1}^r b_i(L\cdot C_i)
\leq d+1.
\]
If $L\cdot C_i>0$ for every $i$, it follows that
\[
0\leq b_i\leq
\left\lfloor\frac{d+1}{L\cdot C_i}\right\rfloor,
\qquad i=1,\ldots,r,
\]
providing an explicit finite search region.

For a smooth projective toric variety $X_\Sigma$, the ambient Mori cone is likewise generated by finitely many torus-invariant curve classes $[\Gamma_1],\ldots,[\Gamma_s]$. These generators, however, need not form a basis. More importantly, the bundle
\[
L=-K_{X_\Sigma}-c_1(\mathcal E)
\]
need not be positive on all of them, even if $Y$ is Fano.  Hence, even though an effective ambient class can be written as
\[
\beta=\sum_{j=1}^s a_j[\Gamma_j],
\qquad a_j\in\R_{\geq0},
\]
which dimension constraint gives
\[
0\leq L\cdot\beta
=\sum_{j=1}^s a_j(L\cdot\Gamma_j)
\leq d+1,
\]
if some $L\cdot\Gamma_j$ are nonpositive, this inequality alone need not bound the coefficients $a_j$. We therefore seek additional geometric restrictions on the ambient images of effective curves in $Y$, allowing us to obtain a bounded search region without requiring $L$ to be ample on $X_\Sigma$. The idea is to restrict the search region for the contributing effective classes to a bounded rational polytope.

Consider $X_\Sigma$ to be a smooth projective toric variety and 
\[
Y=\mathscr{Z}\left(X_\Sigma,\bigoplus_{\ell=1}^{e}\mathcal{O}_{X_\Sigma}(D_\ell) \right)
\]
a smooth Fano complete intersection. Denote $d=\dim Y$ and $L=-K_{X_\Sigma}-\sum_{\ell=1}^{e}D_\ell$, with $L|_Y=-K_Y$.

Because of the dimension axiom, the only nonzero invariants $\langle c_1(TY),\alpha,\gamma\rangle^Y_{0,\beta'}$ are those which satisfy 
\[
1+\deg_{\mathbb C}\alpha+
\deg_{\mathbb C}\gamma
=d+c_1(TY)\cdot\beta'.
\]
Since $Y$ is Fano, every effective class $\beta'\in H_2(Y,\Z)$ satisfies $c_1(TY)\cdot\beta'\geq0$. Thus, by adjunction, we may restrict our attention to the ambient classes $\beta=\iota_*\beta'\in H_2(X_\Sigma,\Z)$ satisfying \[ 0\leq L\cdot\beta\leq d+1. \]

Furthermore, let $v_i$ denote the primitive ray generators of $\Sigma$ and $D_i$ the corresponding torus-invariant divisors. By Proposition~\ref{prop:torich2}, we identify \[ H_2(X_\Sigma,\R) \cong \left\{ (a_i)_i\in\R^k \;\middle|\; \sum_i a_i v_i=0 \right\}, \] where $a_i=D_i\cdot\beta$ for the class $\beta\in H_2(X_\Sigma,\R)$. 

\begin{lemma} \label{lemma: classes}
    For each $\sigma\in\Sigma$ such that $Y\cap\overline{O(\sigma)}\neq\varnothing$, define 
    \[
    J_\sigma=
    \left\{
    \rho_i\in\Sigma(1)\setminus\sigma(1)
    \;\middle|\;
    \sigma+\rho_i\in\Sigma,\quad
    Y\cap \overline{O(\sigma+\rho_i)}\neq\varnothing
    \right\}
    \]
    where $\sigma(1)$ is the set of rays of $\sigma$, and the rational polyhedral cone
    \[
    \mathcal C_\sigma=
    \left\{
    \beta=(a_i)\in H_2(X_\Sigma,\mathbb R)
    \;\middle|\;
    \begin{array}{ll}
    a_i\geq0 & \text{if }\rho_i\in J_\sigma,\\
    a_i=0 & \text{if }
    \rho_i\notin\sigma(1)\cup J_\sigma
    \end{array}
    \right\}.
    \]
    Then every ambient class of an effective curve cycle on $Y$ belongs to \[ \mathcal C = \sum_{\substack{\sigma\in\Sigma\\ Y\cap\overline{O(\sigma)}\neq\varnothing}} \mathcal C_\sigma. \] 
\end{lemma}
\begin{proof}
    Given an irreducible curve $C \subset Y$, choose $\sigma \in \Sigma$ such that $\overline{O(\sigma)}$ is the smallest toric orbit closure containing $C$. Given a ray $\rho_i$ such that $\rho_i \not\in \sigma(1)$, the Orbit-Cone correspondence gives 
    \[
    D_i\cap \overline{O(\sigma)}=
    \begin{cases}
    \overline{O(\sigma+\rho_i)},
    & \sigma+\rho_i\in\Sigma,\\
    \varnothing,
    & \text{otherwise}.
    \end{cases}
    \]
    Thus, if $\rho_i \in J_\sigma$, we have
    \[
     D_i \cap \overline{O(\sigma)} =  \overline{O(\sigma + \rho_i)} \subseteq  \overline{O(\sigma)}
    \]
    and by the minimality of $\overline{O(\sigma)}$, $C$ can not be contained in $\overline{O(\sigma + \rho_i)}$, and thus $D_i \cdot C \geq 0$. 
    On the other hand, if $\rho_i \not \in \sigma(1) \cup J_\sigma$, either $D_i \cap \overline{O(\sigma)} = \varnothing$, and thus $a_i =0$, or $Y \cap \overline{O(\sigma + \rho_i)} = \varnothing$, and then of course $a_i =0$. Thus, $\iota_*[C]\in\mathcal C_\sigma$. Since every effective curve cycle is a nonnegative integral combination of irreducible curves, its ambient class belongs to $\mathcal C$. 
\end{proof}

\begin{remark}
The definition of $J_\sigma$ involves the incidence of $Y$ with toric orbit closures. For a general section of the globally generated bundle, this incidence can be determined from the corresponding intersection classes on the orbit closures.
\end{remark}

\begin{figure}[ht]
\centering
\begin{tikzpicture}[scale=0.82,>=stealth]
  \fill[black!8] (0,0)--(4,0)--(0,4)--cycle;
  \draw[->] (0,0)--(5,0) node[right] {$w_1$};
  \draw[->] (0,0)--(0,5) node[above] {$w_2$};
  \draw[very thick] (4,0)--(0,4);
  \node[below left] at (0,0) {$0$};
   \node[below] at (4,0) {$\frac{M}{L\cdot w_1}w_1$};
  \node[left] at (0,4) {$\frac{M}{L\cdot w_2}w_2$};
  \node at (1.0,1.15) {$P$};
  \node[above right] at (2,2) {$Q$};
\end{tikzpicture}
\caption{For $M=d+1$, $\mathcal C=\mathbb R_{\geq0}w_1+ \mathbb R_{\geq0}w_2$, and $L\cdot w_i>0$, the cross-section $Q=\{\beta\in\mathcal C\mid L\cdot\beta=M\}$ is the segment joining the two labeled points, while $P=\{\beta\in\mathcal C\mid 0\leq L\cdot\beta\leq M\}$ is its convex hull with the origin.}
\end{figure}
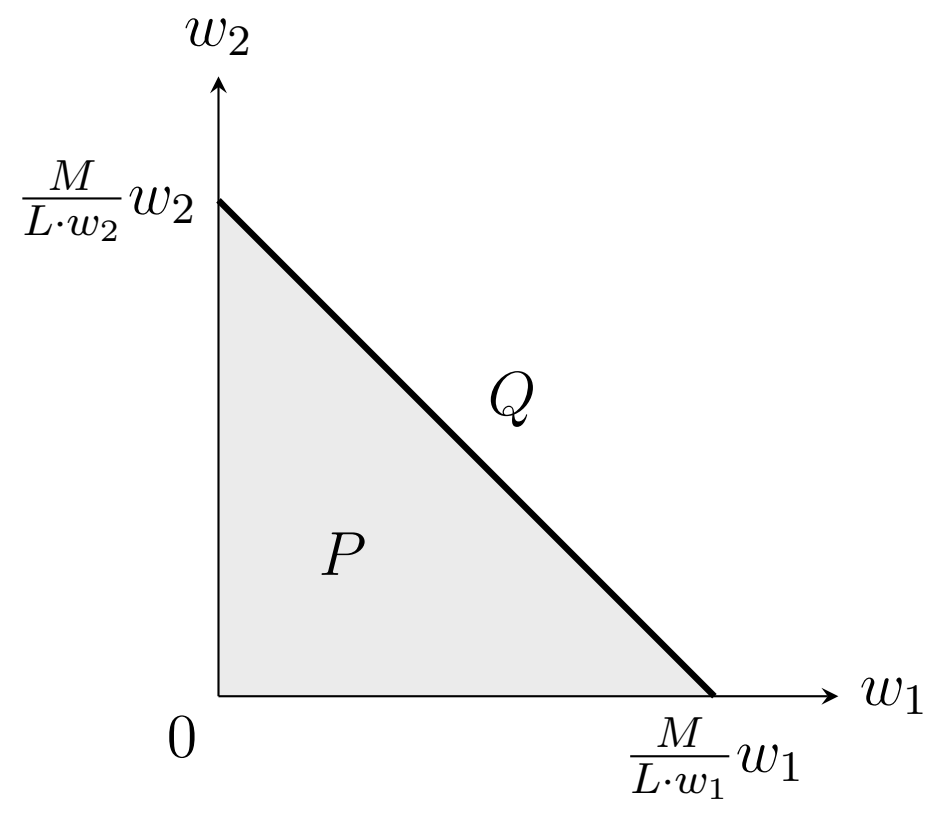

\begin{proposition} \label{prop:finite enum}
     Suppose that $L$ is strictly positive on every nonzero element of $\mathcal C_\sigma$, for every $\sigma\in\Sigma$ such that $Y\cap\overline{O(\sigma)}\neq\varnothing$. Then the set 
    \[
    \mathcal B = \{ \beta \in \mathcal C \cap H_2(X_\Sigma, \Z) \, \, | \, \, 0 \leq L \cdot \beta \leq d+1 \}
    \]
    is finite and contains the ambient image of every effective curve class that can contribute to the small quantum multiplication by $c_1(TY)$.
\end{proposition}

\begin{proof}
  By hypothesis, $L$ is strictly positive on every nonzero element of each $\mathcal C_\sigma$. Since $\mathcal C=\sum_\sigma\mathcal C_\sigma$, it follows that $L$ is strictly positive on $\mathcal C\setminus\{0\}$. Indeed, if $\beta=\sum_\sigma\beta_\sigma\neq0$ with $\beta_\sigma\in\mathcal C_\sigma$, then every $L\cdot\beta_\sigma$ is nonnegative and at least one is positive.

The cone $\mathcal C$ is rational polyhedral, being the sum of finitely many rational polyhedral cones. Consequently,
\[
P=
\left\{
\beta\in\mathcal C
\;\middle|\;
0\leq L\cdot\beta\leq d+1
\right\}
\]
is a bounded rational polytope. Indeed, if $P$ were unbounded, normalizing an unbounded sequence of points in $P$ would produce a nonzero element of $\mathcal C$ on which $L$ vanishes, contradicting strict positivity.

Therefore, $P$ contains only finitely many integral points. By the dimension axiom, a nonzero invariant $\langle c_1(TY),\alpha,\gamma\rangle^Y_{0,\beta'}$ must satisfy 
    \[
    1+\deg_{\mathbb C}\alpha+
    \deg_{\mathbb C}\gamma
    =d+c_1(TY)\cdot\beta'.
    \]
    But since $\deg_{\mathbb C}\alpha, \deg_{\mathbb C}\gamma \leq d$, we have $c_1(TY) \cdot \beta' \leq d+1$. By adjunction and the projection formula, $c_1(TY)\cdot\beta'=L\cdot\iota_*\beta'$. Therefore, the ambient image of every contributing class belongs to $\mathcal B$.
\end{proof}

We can now use this proposition to give an explicit search region for our software. Fix a basis of $H_2(X_\Sigma, \Z)$ and let $\beta=(b_1, \ldots, b_r)$ be the coordinates of $\beta$ on this basis and notice that
\[
\mathcal B = P \cap H_2(X_\Sigma, \Z )
\]
Since $P$ is a bounded rational polytope and the coordinate functions $b_i$ are linear, their maxima and minima are attained at the vertices of $P$. Thus, setting \[ m_i=\left\lceil\min_{\beta\in P}b_i\right\rceil, \qquad M_i=\left\lfloor\max_{\beta\in P}b_i\right\rfloor, \] every integral point of $P$ satisfies $b_i\in[m_i,M_i]\cap\Z$. This gives a finite search box. The points belonging to $\mathcal B$ are precisely those integral tuples in this box that also satisfy the defining inequalities of $P$. We now compute its vertices explicitly. 

Take the cross-section $Q= \{ \beta \in \mathcal C \, | \, L \cdot \beta = d+1\}$. It is easy to see that $P = \text{Conv}(\{0\} \cup Q )$, and thus we simply need to compute $\text{Vert}(Q)$, which are determined by the extremal rays of $\mathcal C$.

We say that a ray $\rho$ of $\mathcal C$ is \emph{extremal} if
\begin{itemize}
    \item $\rho=\R_{\geq0}w$ for some nonzero $w\in\mathcal C$.
    \item Whenever $u,v\in\mathcal C$ satisfy $u+v\in\rho$,
    then $u,v\in\rho$.
\end{itemize}
These are the rays that cross $Q$ at its vertices. By definition, a non-extremal ray will cross at a point that could be written as a convex combination of other points on $Q$ and thus could not be a vertex.  

Hence, if $w_1,\ldots,w_t$ are nonzero generators of the extremal rays of $\mathcal C$, then
\[
\operatorname{Vert}(Q)
= \left\{
\frac{d+1}{L\cdot w_i}\,w_i
\;\middle|\;
i=1,\ldots,t
\right\}.
\]
Consequently,
\[
\operatorname{Vert}(P)=\{0\}\cup\operatorname{Vert}(Q).
\]

\bigskip

\end{document}